\documentclass[10pt]{article}

\usepackage{amsmath}
\usepackage{amsfonts}
\usepackage{amssymb}
\usepackage{framed}
\usepackage{mathtools}
\usepackage{enumerate}
\usepackage{mathrsfs}
\usepackage[hidelinks]{hyperref}
\usepackage{enumitem}
\usepackage{mathrsfs}
\usepackage{scrextend}
\usepackage{amsthm}
\usepackage{bbold}
\usepackage{bbm}
\usepackage{thmtools}
\usepackage{thm-restate}
\usepackage{mathabx}
\usepackage[margin=1.5in]{geometry}
\usepackage{rotating}

\usepackage{xcolor}
\hypersetup{
    colorlinks,
    linkcolor={red!50!black},
    citecolor={blue!50!black},
    urlcolor={blue!80!black}
}

\usepackage{tikz-cd}

\theoremstyle{plain}
\newtheorem{thm}{Theorem}[section]
\newtheorem{lem}[thm]{Lemma}
\newtheorem{prop}[thm]{Proposition}
\newtheorem{cor}[thm]{Corollary}

\theoremstyle{definition}
\newtheorem{defn}[thm]{Definition}

\newtheorem{sit}[thm]{Situation}

\theoremstyle{remark}
\newtheorem*{rem}{Remark}

\newtheorem*{notn}{Notation}

\tikzset{
  symbol/.style={
    draw=none,
    every to/.append style={
      edge node={node [sloped, allow upside down, auto=false]{$#1$}}}
  },
    labl/.style={anchor=south, rotate=90, inner sep=.5mm}
}

\newcommand{\Spec}{\textrm{Spec} \hspace{0.15em} }

\newcommand\restr[2]{{
	\left.\kern-\nulldelimiterspace
	#1
	\vphantom{\big|}
	\right|_{#2}
	}}
\newcommand{\an}[1]{#1^{\textrm{an}}}
\newcommand{\ad}[1]{#1^{\textrm{ad}}}

\newcommand{\ep}{\varepsilon}
\newcommand{\Aut}{\textrm{Aut}}

\newcommand{\Hom}{\textrm{Hom}}

\newcommand{\ch}[1]{\widecheck{{#1}}}

\newcommand{\NL}{\textrm{NL}}

\newcommand{\GL}{\textrm{GL}}
\newcommand{\Spa}{\textrm{Spa}}

\newcommand{\Tot}{\textrm{Tot}}

\usepackage{amsmath,calligra,mathrsfs}

\usepackage[cal=boondoxo]{mathalfa}

\DeclareMathOperator{\sheafhom}{\mathcal{H \kern -1pt o \kern -2pt m}}
\DeclareMathOperator{\sheafend}{\mathcal{E \kern -1pt n \kern -2pt d}}
\DeclareMathOperator{\sheafaut}{\mathcal{A \kern -1pt u \kern -2pt t}}

\tikzset{
  trim node/.default=1cm,
  trim node/.style={
    overlay,
    append after command={
      ([xshift={+#1}]\tikzlastnode.north west)
      ([xshift={+-#1}]\tikzlastnode.south east)}},
  down and trim/.default=1cm,
  down and trim/.style={
    yshift=-(\pgfmatrixcurrentcolumn-1)*1.5\baselineskip,
    trim node={#1}},
  downup and trim/.default=1cm,
  downup and trim/.style={
    yshift=iseven(\pgfmatrixcurrentcolumn) ? -1.5\baselineskip : 0pt,
    trim node={#1}},
  -|/.style={to path={-|(\tikztotarget)\tikztonodes}},
  |-/.style={to path={|-(\tikztotarget)\tikztonodes}},
  -| sl/.style={-|, xslant=-1},
  |- sl/.style={|-, xslant= 1},
  center picture/.style={
    trim left=(current bounding box.center),
    trim right=(current bounding box.center)}}

\title{Geometric G-functions and Atypicality}
\author{David Urbanik}

\begin{document}

\maketitle

\begin{abstract}
We describe a general method for giving $p$-adic interpretations of $G$-functions arising from degenerating periods of smooth projective algebraic varieties. Using this, we are able to implement a strategy due to Andr\'e for bounding heights of moduli points where period functions acquire unusual algebraic relations. This leads to new results on Galois lower bounds for special moduli, and new cases of the Zilber-Pink conjecture. In particular, we establish the first Galois-orbit lower bounds on CM moduli in non-Shimura settings.

As a more technical contribution, we introduce a streamlined Pila-Zannier strategy capable of handling Zilber-Pink-type atypical intersection problems in arbitrary dimension and for arbitrary smooth projective families.
\end{abstract}

\tableofcontents

\section{Introduction}

Let $g \geq 2$, and $\mathcal{M}_{g}$ be the moduli stack of genus $g$ curves with universal family $\mathcal{C} \to \mathcal{M}_{g}$. Let $\overline{\mathcal{M}}_{g}$ be its compactification which parameterizes stable curves, as constructed by Deligne-Mumford in \cite{zbMATH03289080}, and let $\mathcal{B} \subset \overline{\mathcal{M}}_{g} \setminus \mathcal{M}_{g}$ be the locus of stable curves $\mathcal{C}_{x} = C_{1} \cup \cdots \cup C_{\ell}$ with smooth components for which 
\begin{equation}
\label{genuscond}
p_{g}(C_{1}) + \cdots + p_{g}(C_{\ell}) \leq g-2 ,
\end{equation}
where $p_{g}(C_{i})$ denotes the genus of the curve $C_{i}$. This condition is equivalent to $\delta - \ell \geq 1$, where $\delta$ is the number of nodal singularities of $C$. 

Our first result concerns the Jacobians of points of $\mathcal{M}_{g}$. We recall that every abelian variety $A$ over a field $k$ admits a unique isogeny decomposition $A \approx A_{1} \times_{k} \cdots \times_{k} A_{\ell}$ into simple factors, where the relation $\approx$ is given by the existence of an isogeny. Write $\mathcal{S} \subset \mathcal{M}_{g}(\mathbb{C})$ for the set of $x \in \mathcal{M}_{g}(\mathbb{C})$ for which the Jacobian $J(\mathcal{C}_{x})$ admits an isogeny factor with complex multiplication. 

\begin{thm}
\label{modulicurvesthm}
Let $S \subset \mathcal{M}_{g,\mathbb{C}}$ with $g \geq 2$ be an irreducible Hodge-generic curve whose closure $\overline{S} \subset \overline{\mathcal{M}}_{g}$ intersects $\mathcal{B}$. Then $S(\mathbb{C}) \cap \mathcal{S}$ is finite.
\end{thm}

Because there exist non-trivial families of curves for which the Jacobian above a generic point has a CM factor, some kind of Hodge-genericity hypothesis on $S$ like the one given in \autoref{modulicurvesthm} is clearly necessary, but the role of the hypothesis concerning the intersection of $\overline{S}$ with the boundary locus $\mathcal{B}$ is less clear. This will turn out to be an artifact of our method, which proceeds by analyzing cohomological periods near normal crossing degenerations. For instance, essentially the same argument will also give the following: 

\begin{defn}
\label{vanishingcycdef}
Given a normal crossing singularity with equation $x_{1} \cdots x_{r} = 0$, we say a cycle is induced by this singularity if it is obtained as the product of $k$ simple loops around components $x_{i_{1}} = 0, \hdots, x_{i_{k}} = 0$ for some subset $\{ i_{1}, \hdots, i_{k} \} \subset \{ 1, \hdots, r \}$.
\end{defn}

\begin{thm}
\label{abfamthm}
Let $f : X \to S$ be a degenerating Hodge-generic family of abelian varieties of dimension $g \geq 3$ over an irreducible complex algebraic curve. Suppose that at some point $s_{0} \in S(\mathbb{C})$ the fibre $X_{s_{0}}$ has simple normal crossing singularities which induce at least two vanishing cycles in cohomological degree one which are linearly independent in a nearby fibre. Then the set $\mathcal{S} \subset S(\mathbb{C})$ of points $s \in S(\mathbb{C})$ for which $X_{s}$ admits an isogeny factor with complex multiplication is finite. 
\end{thm}

\begin{rem}
A moduli interpretation of \autoref{abfamthm} in the style of \autoref{modulicurvesthm} is also possible. However compactified moduli spaces of abelian varieties are typically constructed using semiabelian schemes instead of singular normal crossing varieties as degeneration points, and although one can always construct a normal crossing compactification of a semiabelian variety, even in families \cite{zbMATH00049143}, there is no canonical way to do so. Thus to understand which curves in $\mathcal{A}_{g}$ are covered by \autoref{abfamthm} one has to undertake an analysis of the different normal crossing compactifications of the fibres of boundary strata in $\overline{\mathcal{A}}_{g}$, and as this will lead us too far astray we do not attempt this here. 
\end{rem}

The methods underlying the proof of both these theorems apply to any smooth projective family over a curve admitting a simple normal crossing degeneration. The key input is ultimately a lower bound on the degree of the field of definition of points in the set $\mathcal{S}$ in terms of a Weil height on $S$. The crucial intermediate theorem is \autoref{genheightthm}. To state it we recall the notion of complex multiplication in the Hodge-theoretic setting.

\begin{defn}
\label{CMdef}
We say that a $\mathbb{Q}$-Hodge structure $V$ of algebro-geometric origin has complex multiplication if for every simple $\mathbb{Q}$-Hodge summand $W \subset V$, the algebra of endomorphisms of $W$ generated by algebraic cycles is a field with $\mathbb{Q}$-vector space dimension equal to $\dim_{\mathbb{Q}} W$. 
\end{defn}

For the relationship between \autoref{CMdef} and the more common description in terms of abelian Mumford-Tate groups see \cite[V.2, V.3]{GGK}; under the Hodge conjecture, CM Hodge structures of algebro-geometric origin in our sense have abelian Mumford-Tate groups, as follows from \cite[V.2, V.3]{GGK} and the fact that the Mumford-Tate group of a direct sum of Hodge structures is contained in the product of the Mumford-Tate groups of the summands. We adopt this more geometric definition since we will explicitly study CM endomorphisms in this paper.

\begin{thm}
\label{genheightthm}
Let $K \subset \mathbb{C}$ be a number field, and suppose that $f : X \to S$ is a projective $K$-algebraic family of geometrically connected varieties over the curve $S$ with smooth generic fibre, and $S$ smooth. Assume that the primitive local subsystem $\mathbb{V}' \subset R^{w} \an{f'}_{*} \mathbb{Q}$ is absolutely simple, where $f'$ is base-change of $f$ to the locus $S' \subset S$ above which the fibres are smooth. Suppose additionally that at some point $s_{0} \in S(K)$ the fibre $X_{s_{0}}$ has simple normal crossing singularities which induce at least two primitive degree $w$ vanishing cycles which are linearly independent in a nearby fibre. Write $\mathcal{S} \subset S(\mathbb{C})$ for the set of points $x \in S(\mathbb{C})$ at which:
\begin{itemize}
\item[-] the Hodge conjecture for endomorphisms of the fibre $\mathbb{V}'_{x}$ holds; and
\item[-] the fibre $\mathbb{V}'_{x}$ has a $\mathbb{Q}$-Hodge summand with complex multiplication.
\end{itemize}
Then $\mathcal{S} \subset S(\overline{\mathbb{Q}})$, and for any logarithmic Weil height $\theta : S(\overline{\mathbb{Q}}) \to \mathbb{R}$ there exists constants $\kappa, a \in \mathbb{R}_{> 0}$ such that 
\[ \theta(\xi) \leq \kappa \, [K(\xi) : K]^{a} \]
for all $\xi \in \mathcal{S}$.
\end{thm}

Note that \autoref{genheightthm} already implies:

\begin{cor}
In the situation of \autoref{genheightthm}, for any constant $d > 0$, the number of points of $\mathcal{S}$ lying inside a number field of degree at most $d$ is finite.
\end{cor}

\noindent Indeed, \autoref{genheightthm} gives an absolute height bound in this case, so this is just the Northcott property. 

\vspace{0.5em}

\autoref{genheightthm} produces unconditional results even in higher weight settings beyond the Shimura case. Indeed, the following is a formal consequence of \autoref{genheightthm}:

\begin{thm}
\label{genCMpointthm}
Let $f, K, \mathbb{V}', f', s_{0}, X_{s_{0}}$ and $w$ be as in \autoref{genheightthm}, and write $\mathcal{S} \subset S(\mathbb{C})$ for the set of points $x \in S(\mathbb{C})$ for which the Hodge structure $\mathbb{V}'_{x}$ has complex multiplication. Then for any logarithmic Weil height $\theta : S(\overline{\mathbb{Q}}) \to \mathbb{R}$ there exists constants $\kappa, a \in \mathbb{R}_{> 0}$ such that 
\[ \theta(\xi) \leq \kappa \, [K(\xi) : K]^{a} \]
for all $\xi \in \mathcal{S}$. In particular, for each $d > 0$, the number of points of $\mathcal{S}$ lying in a number field of degree at most $d$ is finite.
\end{thm}

\noindent Note that the condition that we have two distinct vanishing cycles can be verified in explicit cases by computing the limiting mixed Hodge structure associated to the degeneration point $s_{0}$; for instance, Andr\'e computes in \cite[IX, \S4.4]{zbMATH00041964} the number of such vanishing cycles for cohomology in middle degree in terms of a cohomological invariant of a dual graph associated to the singularities in the special fibre. We will use this ourselves when we prove \autoref{modulicurvesthm} and \autoref{hypsufcor} below. For instance, one can easily verify the hypotheses of \autoref{genCMpointthm} to show the following:

\begin{cor}
\label{hypsufcor}
The conclusion of \autoref{genCMpointthm} holds when $f'$ is a Hodge-generic family of smooth projective hypersurfaces of degree $d$ in $\mathbb{P}^{n+1}$, with $d > n+2$, $X$ is smooth, and where the fibre $X_{s_{0}}$ is a union of $d$ hyperplanes in general position.
\end{cor}

\begin{rem}
The general position assumption is much too strong, one really only needs the hyperplane arrangement to not be overly degenerate; see the proof in \S\ref{hyperplanesec} for details.
\end{rem}

\noindent As far as we are aware, results like this beyond the setting of Shimura-type families have not appeared in the literature before. (Recent work in \cite{papas2022height} gives results under additional arithmetic hypotheses and conjectures near the degeneration point $s_{0}$; we discuss related work in the next section.)

\subsection{Strategy and Previous Work}

\subsubsection{The Pila-Zannier Strategy}
\label{pilazanniersec}

The results we prove fall broadly into the theory of unlikely intersections, and in particular work on the Zilber-Pink conjecture. Let us begin by illustrating the sort of phenomenon one wants to understand with reference to \autoref{modulicurvesthm} above. The set of points $\mathcal{S} \subset \mathcal{M}_{g}$ where the Jacobian acquires a CM isogeny factor is infinite, and in particular contains the complex points of infinitely many subvarieties of $\mathcal{M}_{g}$ whose dimension is a linear function of $g$. This can be seen, for instance, by intersecting the image of the Torelli map with Hecke translates of loci like $\mathcal{A}_{g-1} \times \{ y \} \subset \mathcal{A}_{g}$, with $y$ a point corresponding to a CM elliptic curve (see \cite[Thm. 1.6(i)]{khelifa2023existence} and \cite{eterovic2022likely}). What one is tasked with doing is ruling out the possibility that infinitely many of these subvarieties might intersect some curve in $\mathcal{M}_{g}$; that only finitely many should intersect is a consequence of the Zilber-Pink conjecture. Beyond the special case of the (now resolved) Andr\'e-Oort conjecture, much less is known about Zilber-Pink-type problems, even when considering families over a curve.

The most common approach to Zilber-Pink-type problems is something often referred to as the ``Pila-Zannier strategy''. So as to limit the discussion to the situation which will be relevant to us, we discuss this strategy only in the situation where one has a parameter space $S$ defined over a number field $K \subset \mathbb{C}$, and a set of ``special'' points $\mathcal{S} \subset S(\overline{\mathbb{Q}})$ which one believes is not Zariski-dense. For illustrative purposes we will divide the approach into three steps:

\begin{itemize}[leftmargin=+0.75in]
\item[\textbf{Step (1):}] Prove a bound of the form
\[ \theta(\xi) \leq \kappa [K(\xi) : K]^{a} , \]
for all $\xi \in \mathcal{S}$, where $\theta : S(\overline{\mathbb{Q}}) \to \mathbb{R}$ is a (logarithmic) height function and the constants $\kappa, a > 0$ are independent of $\xi$.
\item[\textbf{Step (2):}] Estimate the heights of some algebraic numbers with uniformly bounded field of definition constructed from the Betti-de Rham periods of $X_{\xi}$ in terms of a polynomial in the height $\theta(\xi)$ and the degree $[K(\xi) : K]$.
\item[\textbf{Step (3):}] Deduce from (1) and (2) that some transcendental analytic set $\mathcal{A}$ contains lots of points defined over fields of uniformly bounded degree, and with the number of points of height at most $N$ growing faster than $N^{\ep}$ for some $\ep > 0$. Use this to contradict a functional transcendence principle.
\end{itemize}

In the particular case of a special point (we follow here the terminology of \cite{BKU})  problem of Andr\'e-Oort type, where $S$ may be viewed as a subvariety $S \subset \Gamma \backslash D$ of a Shimura variety and $\mathcal{S}$ is a set of CM moduli, the transcendental set $\mathcal{A}$ may be taken to be the inverse image of $S$ inside a fundamental set $\mathcal{F} \subset D$ for the action of $\Gamma$. The algebraic numbers whose heights are estimated in terms of $\theta(\xi)$ and $[K(\xi) : K]$ are coordinate-values of lifts $\widetilde{\xi} \in \mathcal{F}$ of the points $\xi$, where the coordinate functions are algebraic functions coming from the embedding $D \subset \ch{D}$ into an ambient algebraic flag variety. If one believes that $\mathcal{S}$ is infinite, then Step (3) proves there is some $\ep > 0$ such that the number of points $\widetilde{\xi}$ of height at most $N$ grows faster than $N^{\ep}$. A theorem of Pila-Wilkie \cite[Theorem 1.6]{zbMATH05680945} shows this can only happen on account of algebraic subvarieties of $\ch{D}$ passing through the inverse image of $S$ in $\mathcal{F}$. The Ax-Schanuel Theorem \cite{AXSCHAN} characterizes such algebraic intersections with $\mathcal{F}$ completely, and can then be used to show that all but finitely many points of $\mathcal{S}$ lie in the intersection of $S$ with special subvarieties of $\Gamma \backslash D$. 

In the more general setting where the points in $\mathcal{S}$ are not CM, but merely have smaller-than-expected Mumford-Tate groups, the strategy becomes more complicated. Instead of taking $\mathcal{A}$ to be the transcendental set $\mathcal{F} \subset D$, we instead consider a transcendental subset $\mathcal{A}$ of a moduli space of flag subvarieties of $\ch{D}$. These flag subvarieties give rise to the points of $\mathcal{S}$ through their intersection with the preimage $\mathcal{F}_{S} \subset \mathcal{F}$ of $S$ inside $\mathcal{F}$. By controlling the heights of the moduli of these subvarieties, we are able to produce, via the Pila-Wilkie theorem, an algebraic family of flag subvarieties intersecting $\mathcal{F}_{S}$ in a positive-dimensional locus of unexpectedly-large dimension. The Ax-Schanuel theorem can once again be used to conclude.

\vspace{0.5em}

Step (1) is fully understood for Andr\'e-Oort-type problems in the Shimura setting, but open in general. Step (2) may be understood to a great extent in the case of abelian families (and more generally, abelian-type Shimura varieties) using isogeny-degree estimates due to Masser-W\"ustholtz. We will give a very general approach to Step (3) in \S\ref{pilazansec7} capable of handling special moduli in arbitrary smooth projective families; see \S\ref{pilazannatypintro} for further discussion. We note that the only reason we do not obtain the finiteness of all special points in $S(\overline{\mathbb{Q}})$ in the setting of \autoref{genCMpointthm} and \autoref{hypsufcor} is only due to the fact that Step (2) has yet to be understood beyond the Shimura setting.

\subsubsection{Height Bounds from $G$-functions}

Our most significant contribution will be to demonstrate, in broad generality, a new way of achieving Step (1) when $S$ is a curve and the family $X \to S$ has a suitable degeneration. The basic idea dates back to Andr\'e's work \cite{zbMATH00041964}, where he developed a method for giving height bounds on special moduli based on a transcendence principle due to Bombieri \cite[\S11]{zbMATH03721021}. The point is that one studies periods near a degeneration point $s_{0}$ of a smooth projective family which, when expanded in a uniformizing parameter $s$ at $s_{0}$, are given by power series $h_{1}, \hdots, h_{m}$ with coefficients in a number field $K$. Bombieri's theorem, which goes by both the name ``Principle of Global Relations'' and ``Hasse Principle for $G$-functions'', gives a way of producing height bounds on special moduli $\xi$ provided one can produce $K$-algebraic relations on the values $h_{1}(\xi), \hdots, h_{m}(\xi)$ which hold at every place of $K$ where $\xi$ lies within the radius of convergence of the $h_{1}, \hdots, h_{m}$. 

In his work \cite{zbMATH00041964}, Andr\'e gave a way of producing such relations which hold at infinite places, but needed to work in a setting where integrality assumptions on both the special moduli and the family $f : X \to S$ ensured that $\xi$ would be bounded away from $s_{0}$ at the finite places. The reason this was necessary, as is discussed by Andr\'e in \cite[pg.8-10]{zbMATH00041964}, is the lack of a suitable way of interpreting the power series $h_{1}, \hdots, h_{m}$ at the finite places; in particular, one would like to be able to produce $K$-algebraic relations on the values of $h_{1}(\xi), \hdots, h_{m}(\xi)$ holding at \emph{all} finite places. Follow up work, both by Andr\'e \cite{zbMATH00764346} \cite{zbMATH00903629}, and more recently by Daw and Orr \cite{daw2022zilber}, and by Papas \cite{papas2022height}, has either similar limitations, or is only able to resolve these issues for families coming from elliptic curves. (We note here that one needs \emph{both} a $p$-adic interpretation of these functions, \emph{and} a way of producing $K$-algebraic relations on the values $h_{1}(\xi), \hdots, h_{m}(\xi)$ at finite places; it is not difficult to interpret the functions $h_{1}, \hdots, h_{m}$ using rigid-analytic de Rham cohomology, but this alone is insufficient.) Consequently, one is presently unable to use the $G$-function method to produce height bounds in many general settings, and in particular to resolve Step (1) of the Pila-Zannier strategy for Zilber-Pink-type problems where the $G$-function method could be expected to apply.

Our primary contribution will be to resolve this issue in broad generality, and show how one can use a combination of rigid analytic de Rham cohomology and $p$-adic Hodge theory to interpret Andr\'e's $G$-functions $p$-adically, and moreover use the Galois-invariance of the $p$-adic Hodge comparison to produce $K$-algebraic relations on the values of $h_{1}(\xi), \hdots, h_{m}(\xi)$ at special moduli. To explain how we will do this, we now give a more detailed description of Andr\'e's construction.

\subsection{The $G$-function Method}
\label{Gfuncmethodintro}

Suppose that $f : X \to S$ is a projective family of relative dimension $n = \nu - 1$ defined over a number field $K \subset \overline{\mathbb{Q}} \subset \mathbb{C}$, with $X$ and $S$ both smooth, $S$ a geometrically irreducible curve, and which has geometrically-connected fibres. Denote by $f' : X' \to S'$ its base-change to the open locus $S' \subset S$ above which the fibres are smooth, and suppose that $S \setminus S' = \{ s_{0} \}$ with $s_{0} \in S(K)$. We fix some subset $\mathcal{S} \subset S(\overline{\mathbb{Q}})$ of special moduli; this means that for each point $\xi \in \mathcal{S}$ some self-product
\[ X^{n}_{\xi} = \underbrace{X_{\xi} \times \cdots \times X_{\xi}}_{n \textrm{ times}} \]
carries an algebraic cycle class which does not spread out to the generic point of $S$. We also fix a (logarithmic) Weil height $\theta : S(\overline{\mathbb{Q}}) \to \mathbb{R}$ from which one obtains a function $\mathcal{S} \to \mathbb{R}$ which we also denote by $\theta$. 

We now recall Andr\'e's setup from \cite{zbMATH00041964}. He considers the case where monodromy around $s_{0}$ acts via a unipotent linear transformation, and for which the fibre $X_{0} \subset X$ at $s_{0}$ degenerates via a reduced normal crossing; this latter condition means that:
\begin{itemize}
\item[-] there is an affine open subset $U \subset X$ and functions $z_{1}, \hdots, z_{\nu}$ on $U$ whose differentials trivialize $\Omega^{1}_{U}$; and
\item[-] the equation $z_{1}^{e_{1}} \cdots z^{e_{\nu}}_{\nu} = s$ defines the graph of $\restr{f}{U}$ near $s_{0}$, where $s$ is a uniformizing function on $S$ at $s_{0}$, and $e_{j} \in \{ 0, 1 \}$ for all $1 \leq j \leq \nu$.
\end{itemize}
After reordering we obtain an integer $\mu$ such that that $e_{j} = 0$ for $j > \mu$ and $e_{j} = 1$ otherwise. On $U$ one can then fix a point $q$ in the locus $z_{1} = \cdots = z_{\mu} = 0$, hence in the singular locus of $X_{0}$, and consider, in an analytic neighbourhood of $s_{0}$, complex analytic functions of $s$ given by 
\begin{equation}
\label{Pdef}
 h(s) = \frac{1}{(2 \pi i)^{\mu-1}} \int_{\gamma_{s}} \iota_{s}^{*}\left( g^{q} \frac{dz_{2} \cdots dz_{\mu}}{z_{2} \cdots z_{\mu}} \right) , 
\end{equation}
where
\begin{itemize}
\item[-] $\iota_{s} : X_{s} \cap U \hookrightarrow U$ is the inclusion of the fibre above $s$; 
\item[-] $\gamma_{s}$ is a ``vanishing cycle'' in the fibre $U_{s}$ obtained as a product of $\mu-1$ simple loops near $q$ around the divisors $z_{j} = 0$ for $j = 2, \hdots, \mu$; and
\item[-] $g^{q}$ is a holomorphic function chosen so that $g^{q} \frac{dz_{2} \cdots dz_{\mu}}{z_{2} \cdots z_{\mu}}$ represents the restriction of a relative class in the algebraic de Rham cohomology of $X$ over $S$ and whose power series representation in the coordinates $z_{1}, \hdots, z_{\mu}$ at $q$ has coefficients in $K$.
\end{itemize}
As explained by Andr\'e in \cite[IX,\S4]{zbMATH00041964}, the functions $h(s)$ are described by power series in $s$ with coefficients in $K$ when expanded around $s_{0}$, and give, in a punctured neighbourhood around $s_{0}$, a relative period of the Betti-algebraic de Rham comparison $R^{\mu-1} f'_{*} \mathbb{Z}(\mu-1) \otimes \mathcal{O}_{\an{S'}} \simeq \an{(R^{\mu-1} f'_{*} \Omega^{\bullet}_{X'/S'})_{\mathbb{C}}}$. In particular, this power series representation of $h$ is a $G$-function in the sense of \cite[I]{zbMATH00041964}.

\begin{rem}
Andr\'e actually assumes that $\mu = \nu$, and assumes that $X$ is covered by neighbourhoods $U$ of the above type. As it is not substantially more difficult, we will work in greater generality.
\end{rem}

To explicate the relationship between these $G$-functions and the projective family $f$ Andr\'e classifies, at least in degree $n$, the period functions $h$ of this form in terms of the monodromy around the degeneration point. To explain what we mean, let us fix an integer $w$ and denote by $\mathbb{V}'$ the variation of Hodge structure with underlying local system $R^{w} f'_{*} \mathbb{Z}(w) / tor.$, and let $\mathcal{H}' = R^{w} f'_{*} \Omega^{\bullet}_{X'/S'}$ be the associated algebraic de Rham cohomology vector bundle. The vector bundle $\mathcal{H}'$ has a so-called \emph{canonical extension} $\mathcal{H}$ to a vector bundle over $S$ which we recall in \S\ref{canonextsec}. The sections of $\mathbb{V}'$ that extend to sections of $\an{\mathcal{H}_{\mathbb{C}}}$ under the comparison isomorphism span a local subsystem $\mathbb{M} \subset \restr{\mathbb{V}'_{\mathbb{Q}}}{B}$, where $B \subset S(\mathbb{C})$ is a small analytic ball centered at $s_{0}$. Using the duality between cohomology and homology in the case of $\mathbb{Q}$ coefficients, $\mathbb{M}$ corresponds to a subsystem $\mathbb{M}^{*} \subset \restr{\mathbb{V}'^{*}_{\mathbb{Q}}}{B}$ of monodromy-invariant homology classes. We then obtain a natural pairing
\begin{equation}
\label{poinpairing}
\restr{\an{\mathcal{H}'}}{B} \otimes \mathbb{M}^{*} \to \mathcal{O}_{B} ,
\end{equation}
and sections in the image of this pairing we refer to as \emph{non-degenerating} (relative) periods at $s_{0}$; the functions $h(s)$ described above were of this type. Andr\'e gives a description of the image in the case when $w = \nu-1$, and characterizes the span in the image of (\ref{poinpairing}) by $G$-functions of the form (\ref{Pdef}) using the monodromy-weight filtration. 

In what follows we write $\mathbf{M}_{m \times n}(A)$ for the space of $m \times n$ matrices with values in the ring $A$; we will also write $\mathbf{M}_{m}(A)$ for $\mathbf{M}_{m \times m}(A)$. Andr\'e then uses the non-degenerating periods at $s_{0}$ to give a method for bounding the heights of the points in the set $\mathcal{S}$, based on a so-called \emph{Hasse principle for $G$-functions} \cite[VII, \S5]{zbMATH00041964}, itself a variant of a theorem of Bombieri, which may be stated as follows:

\begin{notn}
Given a number field $L$, we write $\Sigma_{L}$ for the set of places of $L$.
\end{notn}

\begin{thm}[Hasse Principle]
\label{hasseprincip}
Suppose that $G = (G_{1}, \hdots, G_{m})^{t} \in \mathbf{M}_{m \times 1}(K[[t]])$ satisfies the differential system $\frac{d}{dt} G = \Gamma G$ for some $\Gamma \in \mathbf{M}_{m}(K(t))$, write $G_{i} = \sum_{j = 0}^{\infty} G_{ij} t^{j}$, and suppose that
\[ \limsup_{n \to \infty} \left( \frac{1}{n} \sum_{v \in \Sigma_{K}} \max_{i \leq m, j \leq n} \log^{+}|G_{ij}|_{v} \right) < \infty \]
where $\log^{+}(a) = \log \textrm{max}\{ 1, a \}$. Let $\Xi(G, \delta)$ denote the set of $\xi \in \overline{\mathbb{Q}}$ satisfying the following property: there exists a homogeneous polynomial $Q \in K[y_{1}, \hdots, y_{m}]$ of degree at most $\delta$ such that: 
\begin{itemize}
\item[(1)] there does not exist $\widetilde{Q} \in K[t][y_{1}, \hdots, y_{m}]$, homogeneous in $y_{1}, \hdots, y_{m}$, specializing to $Q$ at $\xi$ such that the relation $\widetilde{Q}(t,G_{1}, \hdots, G_{m}) = 0$ holds on the level of formal power series; and
\item[(2)] for all $v \in \Sigma_{K(\xi)}$ for which $|\xi|_{v} < 1$, either
\begin{itemize}
\item[(i)] at least one of the series $G_{i}$ does not have $v$-adic convergence radius greater than $|\xi|_{v}$; or
\item[(ii)] the relation $Q(G_{1}(\xi), \hdots, G_{m}(\xi)) = 0$ holds $v$-adically at $\xi$.
\end{itemize}
\end{itemize}
Then there exists an exponent $e \in \mathbb{Z}_{\geq 0}$ such that 
\begin{equation}
\label{hasseprineq}
 \sup_{\xi \in \Xi(G,\delta)} \theta(\xi) = O(\delta^{e}) ,
\end{equation}
where $e$ and the implied constant depend only on $G$, and $\theta$ is the standard logarithmic Weil height function.
\end{thm}

\begin{rem}
In practice one does not actually have to find a single $\Gamma$ for which all the $G_{i}$ belong to a single common solution vector $G$ for the differential operator $\frac{d}{dt} - \Gamma$, and one could instead assume that each $G_{i}$ is a component of such a solution vector for a differential operator $\frac{d}{dt} - \Gamma_{i}$ depending on $i$. (Or what is equivalent, that each $G_{i}$ is a $G$-function.) See \cite[Thm. 2.5]{daw2022zilber} for a proof of this.
\end{rem}

\begin{defn}
In the context of an application of the Hasse principle, a point $\xi \in \overline{\mathbb{Q}}$ and a place $v$ of $K(\xi)$, we will say that $v$ is \emph{relevant} for $\xi$ if $|\xi|_{v} < 1$ and the series $G_{1}, \hdots, G_{m}$ all have $v$-adic convergence radius greater than $|\xi|_{v}$; in particular, a place being relevant means that one has to demonstrate condition (2)(ii) holds if one wants to apply the Hasse principle. 
\end{defn}

\begin{rem}
The constants implicit in the term $O(\delta^{e})$ can be made explicit, and even (in principle) effective, see the footnote in \cite[pg.129]{zbMATH00041964} and the corresponding discussion. 
\end{rem}

\noindent By taking $G$ to be a vector consisting of $G$-functions at $s_{0}$ arising from a basis of $\mathcal{H}'$ and sections of $\mathbb{M}$ Andr\'e uses this principle to bound the height of elements in certain subsets of $\mathcal{S}$; here the polynomials $Q$ are relations on periods coming from algebraic (or absolute Hodge) cycles associated to the cohomology groups of the fibre $X_{\xi}$. We note that if one has different relations $Q_{1}$ and $Q_{2}$ at two different places $\iota_{1}, \iota_{2} \in \Sigma_{K(\xi)}$, one can always consider the product relation $Q = Q_{1} Q_{2}$ which will hold at both $\iota_{1}$ and $\iota_{2}$. The difficulty is ultimately to show that the infinite set $\Sigma_{K(\xi)}$ can be partitioned into finitely many subsets $\Sigma_{\xi,1}, \hdots \Sigma_{\xi,j}$ such that a common relation $Q_{i}$ holds on the values $G_{1}(\xi), \hdots, G_{m}(\xi)$ for each place $\iota \in \Sigma_{\xi,i}$ which is relevant for $\xi$. Then to resolve Step (1) above, one has to ensure that the degree of $Q = \prod_{i} Q_{i}$ is bounded by a polynomial in $[K(\xi) : K]$, and then (\ref{hasseprineq}) leads to an inequality of the type required in Step (1).

\subsection{$p$-adic Interpretations of $G$-functions}

Our analysis will start by giving a purely algebraic-de-Rham description of Andr\'e's G-functions; a similar description already appears in the proof of \cite[IX, \S4, Theorem 2]{zbMATH00041964}. To do this, we once again fix an affine open subset $U \subset X$, with coordinates $z_{1}, \hdots, z_{\nu}$ inducing an \'etale map $U \to \mathbb{A}^{\nu}$, and such that the map to $S$ is given by $s \mapsto z_{1} \hdots z_{\mu}$ for some $1 \leq \mu \leq \nu$; we also set $w = \mu - 1$. The vector bundle $\mathcal{H}'$ extends canonically to a $K$-algebraic vector bundle $\mathcal{H}$ over $S$, as we review in \S\ref{canonextsec}. From our description of $\mathcal{H}$ we will see that any section $\omega$ of $\mathcal{H}$ then admits a restriction $\omega_{U}$ to a relative de Rham sheaf on $U$. Fixing a point $q \in U(K)$ in the locus $z_{1} = \cdots = z_{\mu} = 0$ mapping to $s_{0}$, we will further restrict $\omega_{U}$ to a formal neighbourhood of $q$ to obtain a unique representation
\[ i^{*}_{q} \omega_{U} = h^{q} \frac{dz_{2} \cdots dz_{\mu}}{z_{2} \cdots z_{\mu}} ,  \]
in a formal de Rham complex at $q$, with $h^{q} \in K[[s]]$. This $K$-algebraic power series $h^{q}$ is obtained without leaving the algebro-geometric category, and its analytification agrees with the function $h$ in (\ref{Pdef}) above. Our goal is to provide a $p$-adic interpretation of the same object, for which we begin with a $p$-adic interpretation of vanishing cycles.

\subsubsection{$p$-adic Vanishing Cycles}

Let $v$ be a finite place of $K$ and write $K_{v}$ for the associated completion. For a connected adic space $Y$ defined over $\Spa(K_{v}, \mathcal{O}_{K_{v}})$, we recall in \S\ref{fundgpsec} the definition of the \'etale fundamental group $\pi^{\textrm{\'et}}_{1}(Y, y)$ at a geometric point $y \in Y$. We consider the case where $Y = \Delta^{\circ} = \textrm{Spa}(K_{v}\langle T, T^{-1} \rangle, \mathcal{O}_{K_{v}}\langle T, T^{-1} \rangle)$ is the rigid-analytic torus, and try to describe an element in $\pi^{\textrm{\'et}}_{1}(\Delta^{\circ},y)$ giving a pro-$p$ analogue of a rigid-analytic ``loop'' around the puncture. Unfortunately, the space $\Delta^{\circ}$ admits more rigid-analytic coverings than just those of Kummer type, making a na\"ive approach difficult. Our idea is basically to define this ``loop'' on just those coverings of Kummer type, and then to choose an extension to $\pi^{\textrm{\'et}}_{1}(\Delta^{\circ},y)$ which will be compatible with the formalism of $p$-adic Hodge theory. This is not so easy to do in the (possibly non-abelian) setting of fundamental groups, so we instead work dually, viewing each element of $\pi^{\textrm{\'et}}_{1}(\Delta^{\circ},y)$ through its induced functional on first-degree cohomology via the isomorphism $H^{1}(\Delta^{\circ}_{\textrm{p\'et}}, \hat{\mathbb{Z}}_{p}(1)) \simeq \Hom_{\textrm{cont}}(\pi^{\textrm{\'et}}_{1}(\Delta^{\circ},y), \mathbb{Z}_{p}(1))$. (The notation $(-)_{\textrm{p\'et}}$ denotes the pro-\'etale site introduced by Scholze in \cite{zbMATH06209107}, which we review in \S\ref{proetsec}.) 
 
After fixing a compatible system of $p$'th roots of unity we obtain a functional $\alpha^{*} : H^{1}(\mathbb{G}_{m,\textrm{p\'et}}, \hat{\mathbb{Z}}_{p}(1)) \to \mathbb{Z}_{p}(1)$, where $\mathbb{G}_{m}$ is the multiplicative group, which we then extend to a functional $\hat{\alpha}^{*} : H^{1}(\Delta^{\circ}_{\textrm{p\'et}}, \hat{\mathbb{Z}}_{p}(1)) \otimes \mathbb{Q}_{p} \to \mathbb{Q}_{p}(1)$ compatible with pullback by the map $H^{1}(\mathbb{G}_{m,\textrm{p\'et}}, \hat{\mathbb{Z}}_{p}(1)) \otimes \mathbb{Q}_{p} \to H^{1}(\Delta^{\circ}_{\textrm{p\'et}}, \hat{\mathbb{Z}}_{p}(1)) \otimes \mathbb{Q}_{p}$ induced by the embedding $\Delta^{\circ} \hookrightarrow \mathbb{G}_{m}$ of adic spaces. In the case of a more general space $\Delta^{a,b} = (\Delta^{\circ})^{a} \times \Delta^{b}$ embedding into $\mathbb{G}_{m}^{a} \times \mathbb{A}^{b}$, where $\Delta$ is the open $v$-adic disk of radius $1$, one extends both $\alpha^{*}$ and $\hat{\alpha}^{*}$ to functionals $\alpha^{*}_{a,b}$ and $\hat{\alpha}^{*}_{a,b}$ on the cohomology in degree $a$ using the cup product.

To relate this to the function $h^{q}$ above, we fix a rigid analytic disk $\mathcal{D}$ around $s_{0}$, and a family of neighbourhods $\Delta^{w,\nu-\mu}_{s} \subset X_{s}$ for points $s \in \mathcal{D}$, where $\Delta^{w,\nu-\mu}_{s} \simeq (\Delta^{\circ})^{w} \times \Delta^{\nu-\mu}$. Using our fixed compatible system of $p$'th roots of unity, one obtains a fundamental $p$-adic period $t \in B_{\textrm{dR}}$. We then pull back $\widehat{\alpha}^{*}_{a,b}$ to a functional 
\[ \hat{\gamma}^{*}_{s} : H^{w}(X_{s,\overline{K_{v}},\textrm{p\'et}}, \hat{\mathbb{Z}}_{p}(w)) \to \mathbb{Z}_{p}(w) \] 
and let $\hat{\gamma}^{*}_{s,B_{\textrm{dR}}}$ be its scalar extension. We then obtain the following relationship with the functions $h^{q}$ above:
\begin{prop}
\label{Gfuncinterp1}
For each point $s \in \mathcal{D}(K_{v})$ one has that
\begin{equation}
\label{hfuncid}
h^{q}(s) = \frac{1}{t^{w}} \, \hat{\gamma}^{*}_{s,B_{\textrm{dR}}}\left( \rho^{-1} \left( \omega \right)\right) ,
\end{equation}
where 
\[ \rho : H^{w}(X_{s, \overline{K_{v}}, \textrm{p\'et}}, \widehat{\mathbb{Z}}_{p}(w)) \otimes B_{\textrm{dR}} \to H^{w}_{\textrm{dR}}(X_{s}) \otimes B_{\textrm{dR}} \]
is the $p$-adic period isomorphism
\end{prop}

\begin{rem}
To make sense of the equality (\ref{hfuncid}) we have fixed an embedding $K \hookrightarrow K_{v}$, and the result holds for all such choices.
\end{rem}

\subsubsection{Producing Relations on $p$-adic Periods}

After fixing a frame $\omega_{1}, \hdots, \omega_{m}$ of $\restr{\mathcal{H}'}{\mathcal{D}}$ with corresponding representations
\[ i^{*}_{q} \omega_{j} = h_{j} \frac{dz_{2} \cdots dz_{\mu}}{z_{2} \cdots z_{\mu}} , \]
our result \autoref{Gfuncinterp1} can be interpreted as follows:

\begin{thm}
The values $h_{1}(s), \hdots, h_{m}(s)$ give the vector 
\[ \frac{1}{t^{w}} (\hat{\gamma}^{*}_{s,\textrm{B}_{dR}} \circ \rho^{-1}) \in \Hom(H^{\mu-1}_{\textrm{dR}}(X_{s}) \otimes B_{\textrm{dR}}, B_{\textrm{dR}}) \]
in the dual coordinates induced by $\omega_{1,s}, \hdots, \omega_{m,s}$.
\end{thm}

To apply this fact to produce $K$-algebraic relations on the values of $h_{1}, \hdots, h_{m}$ at a special point we now consider the situation where $X_{s}$ admits a non-trivial algebra $E$ of algebraic self-correspondences. First, we observe the following, which is almost immediate from the construction of $\hat{\gamma}^{*}_{s}$:

\begin{prop}
Let $v$ be the place of $K$ at which the functional $\hat{\gamma}^{*}_{s}$ is defined, and $G_{K_{v}} = \textrm{Gal}(\overline{K_{v}}/K_{v})$ the associated local Galois group. Then $G_{K_{v}}$ acts on $\hat{\gamma}^{*}_{s}$ through the character $\chi_{\textrm{cycl}}^{w}$, where $\chi_{\textrm{cycl}} : G_{K_{v}} \to \mathbb{Z}^{*}_{p}$ is the usual cyclotomic character. 
\end{prop}

\noindent The point is that the only choice not invariant under the Galois action made in the construction of $\hat{\gamma}^{*}_{s}$ is our choice of a non-trivial compatible system of $p$-power roots of unity corresponding to a $p$-adic ``loop'' inside the torus $\Delta^{\circ}$, and if we ``integrate'' around $w$ such loops then $G_{K}$ acts on the ``integrals'' through the $w$'th power of $\chi_{\textrm{cycl}}$. We observe that this simple fact is already enough to produce non-trivial algebraic relations on the de Rham coordinates of $\hat{\gamma}^{*}_{s}$ in the presence of an $L$-algebraic correspondence $\tau : X_{s} \dashrightarrow X_{s}$. Indeed, the cohomology class $[\tau]$ in both de Rham and \'etale cohomology is fixed by a finite index subgroup of $G_{K_{v}}$, hence the functionals
\begin{equation}
\label{pullbackseq}
\hat{\gamma}^{*}_{s}, \, \hat{\gamma}^{*}_{s} \circ [\tau],  \, \hat{\gamma}^{*}_{s} \circ [\tau]^{2},  \, \hat{\gamma}^{*}_{s} \circ [\tau]^{3}  \, \hdots
\end{equation}
all lie in a subspace of functionals on $H^{w}(X_{s, \overline{K_{v}}, \textrm{p\'et}}, \hat{\mathbb{Z}}_{p}(w))$ on which a finite index subgroup of $G_{K}$ acts by $\chi_{\textrm{cycl}}^{w}$. If one expresses the functionals in (\ref{pullbackseq}) in a fixed $L$-algebraic de Rham basis, then the coordinates of each element of the sequence (\ref{pullbackseq}) are $L$-linear combinations of the coordinates of $\hat{\gamma}^{*}_{s}$. The Hodge-Tate comparison provides a simple way of bounding the dimension of the $\chi^{w}_{\textrm{cycl}}$-character space in terms of the de Rham Hodge numbers, and hence one may obtain an $L$-algebraic relation on the coordinates of $\hat{\gamma}^{*}_{s}$ simply by taking determinants of minors of matrices constructed from the vectors in the sequence (\ref{pullbackseq}). To make the relation $K$-algebraic one takes the product of this relation with all its Galois conjugates. (We will instead use a slightly more subtle idea to construct $K$-algebraic relations as products of linear relations, but the above gives the essential idea.)

\vspace{0.5em}

Let us remark that, in the cases we consider, we sometimes don't expect relations like the ones we construct on the coordinates of $\hat{\gamma}_{s}$ to actually exist for points $s$ which are $p$-adically close to $s_{0}$. The reason is that the there can be constraints on which special moduli are allowed to be close to degeneration points in $p$-adic metrics: this is the case for curves and abelian varieties with complex multiplication, all of which have potentially good reduction. In such a case, our construction of non-trivial relations on the coordinates of $\hat{\gamma}^{*}_{s}$ under the (in this case counterfactual) assumption that $s$ is a special modulus $p$-adically close to $s_{0}$ is to be interpreted as a sort of integrality constraint on $s$. On the other hand, there are genuine cases where special moduli do come arbitrarily $p$-adically close to boundary divisors, and in these cases the algebraic relations on special moduli reflect the constraints placed on algebraic cycles by $p$-adic vanishing cycle data.

\subsection{Applications to Height Bounds}

The generality in which one can now construct algebraic relations on Andr\'e's $G$-functions at finite places eliminates a broad class of obstructions to applying the $G$-function method to problems in arithmetic geometry. Indeed one can now show, in quite general settings, that polynomial height bounds on special moduli follow as soon as one can establish $K$-algebraic constraints at the infinite places. 

We now give a sample result of this type, assuming the setup of \S\ref{Gfuncmethodintro} (also recalled in \S\ref{notconvtpt1} below). We also assume the Hodge conjecture holds for endomorphisms appearing in the fibres of $\mathbb{V}'$. We then define:
\begin{align*}
\mathcal{S} &:= \{ \xi \in S(\mathbb{C}) : \textrm{there exists a CM Hodge summand } W \subset \mathbb{V}'_{\xi} \} 
\\
\mathcal{S}_{\textrm{split}} &:= \left\{ \xi \in S(\mathbb{C}) : \substack{\textrm{there exists a Hodge summand } W \subset \mathbb{V}'_{\xi} \\ \textrm{with CM by a field }E\textrm{ such that }\mathbb{V}'_{E}\textrm{ is }E\textrm{-simple } } \right\} .
\end{align*}

\noindent We note that if $\mathbb{V}'$ is absolutely simple, then $\mathcal{S} = \mathcal{S}_{\textrm{split}}$. This will be the case in all the applications we outlined at the beginning of this introduction.

\begin{prop}
\label{onlycontrolinfplaces}
Suppose that $\mathbb{V}'$ is $\mathbb{Q}$-simple, and that the $G$-functions $(h_{1}, \hdots, h_{m})$ associated to the point $q \in X_{s_{0}}$ give the coordinates of a family of non-zero functionals on cohomology. Then after replacing $K$ with a finite extension, there exists
\begin{itemize}
\item[(i)] a finite covering $c : C \to S$, and a parameter $s$ on $C$ with simple zeros and vanishing exactly on the set $c^{-1}(s_{0})$; and
\item[(ii)] for all but finitely many $\xi \in c^{-1}(\mathcal{S})$ a $K(\xi)$-algebraic relation on the values 
\[ h_{1}(s(\xi)), \hdots, h_{m}(s(\xi)) , \]
not induced by a functional relation on $s, h_{1}, \hdots, h_{m}$, and which holds at all finite places relevant for $s(\xi)$.
\end{itemize}
Moreover, the degree of the relation in (ii) is bounded independently of $\xi$.
\end{prop}

\begin{rem}
In \autoref{onlycontrolinfplaces}(ii) we actually mean to replace the original $G$-functions with the ones computed in terms of the parameter $s$; we give a more precise description in \S\ref{heightboundsec}. 
\end{rem}

Using Andr\'e's original strategy for constructing relations at the infinite places, this then leads to the following theorem, which reinterprets \autoref{genheightthm} above:

\begin{thm}
\label{sndsinggivesheights}
Assume there exists two order $w$ normal crossing points $q, q' \in X_{s_{0}}$ whose associated tuples $(h_{1}, \hdots, h_{m})$ and $(h'_{1}, \hdots, h'_{m})$ of $G$-functions correspond to linearly independent functionals in cohomology. Then for any logarithmic Weil height $\theta : S(\overline{\mathbb{Q}}) \to \mathbb{R}_{> 0}$ there exists constants $\kappa, a \in \mathbb{R}_{> 0}$ such that 
\[ \theta(\xi) \leq \kappa \, [K(\xi) : K]^{a} \]
for all $\xi \in \mathcal{S}_{\textrm{split}}$.
\end{thm}

\subsection{Pila-Zannier for General Atypicality}
\label{pilazannatypintro}

Lastly, with reference to our discussion of the Pila-Zannier strategy in \S\ref{pilazanniersec}, let us elaborate on our approach to Step (3), a good portion of which we actually carry out in a general Hodge-theoretic setting. In this setting we have a set $\mathcal{S} \subset S(\mathbb{C})$ of special points, and these induce points $\widetilde{\xi}$ inside a definable period image $\mathcal{I} \subset D$, where $D$ is a period domain for polarized Hodge structures appearing in the fibres of $\mathbb{V}$. If the points in $\mathcal{S}$ one is studying are not CM points, but merely points above which the fibre $X_{s}$ acquires some extra algebraic cycle, the points $\widetilde{\xi}$ are no longer $\overline{\mathbb{Q}}$-algebraic, and merely have lower-than-expected transcendence degree. This can be understood as arising from the intersection between $\mathcal{I}$ and a $\overline{\mathbb{Q}}$-algebraic flag subvariety $\ch{D}_{M} \subset \ch{D}$ determined by the Mumford-Tate group $M$ of $\widetilde{\xi}$. 

Our observation, which is related to ideas appearing in \cite{daw2018applications}, is that one can obtain results in this more general setting by applying Pila-Zannier-type reasoning to the moduli of the varieties $\ch{D}_{M}$. More specifically, one can reduce to the case where one considers only varieties $\ch{D}_{M}$ for which the associated Mumford-Tate groups $M$ lie inside a single $\mathbf{G}_{S}(\mathbb{C})$-orbit for the generic Mumford-Tate group $\mathbf{G}_{S}$ of the variation $\mathbb{V}$, where the action on Mumford-Tate groups is by conjugacy. The situation one is then tasked with dealing with is the situation where there are many $\overline{\mathbb{Q}}$-algebraic translates $g \cdot \ch{D}_{M}$ of the variety $\ch{D}_{M}$ which intersect $\mathcal{I}$ atypically. One can understand the elements $g$ that arise in terms of heights of Hodge tensors defining the associated Mumford-Tate groups $g M g^{-1}$, and use this to bound both the heights of such $g$ and the degree of their field of definition. The Pila-Wilkie theorem then produces, under appropriate bounds on the heights of some Hodge tensors associated to points of $\mathcal{I}$, an algebraic family of subvarieties of $\ch{D}$ which intersect $\mathcal{I}$ atypically, and from this one can run the usual functional transcendence arguments. We do not need any constraints on $S$; in particular, we do not use that $S$ is a curve.

As an application of this, we prove the following general result, which we state here somewhat informally (see \S\ref{pilazansec7} for the relevant definitions and precise statements).

\begin{thm}
\label{informalstatement}
Suppose that $f : X \to S$ is a family of abelian varieties whose derived Mumford-Tate group is $\mathbb{Q}$-simple, and $\mathcal{S} \subset S(\mathbb{C})$ is the subset of points in the zero-dimensional Hodge locus which are defined by, and atypical for, a central Hodge idempotent. Then if there exists constants $\kappa, a \in \mathbb{R}_{> 0}$ such that
\[ \theta(\xi) \leq \kappa \, [K(\xi) : K]^{a} \]
for all $\xi \in \mathcal{S}$, with $\theta$ some logarithmic Weil height, then $\mathcal{S}$ is finite.
\end{thm}

\noindent As alluded to in \S\ref{pilazanniersec}, the fact that $f$ is an abelian family here is only relevant because we are unable to carry out Step (2) without it. We note that we will obtain both \autoref{modulicurvesthm} and \autoref{abfamthm} by applying \autoref{sndsinggivesheights} to the appropriate family and then reducing to \autoref{informalstatement}.

\subsection{Additional Related Work}
\label{relworksec}

As we have discussed in great detail, the $G$-function method for bounding heights on special moduli was introduced in Andr\'e's book \cite{zbMATH00041964}, but was limited by the lack of $p$-adic interpretations of these functions. Some subsequent works (see \cite{zbMATH00764346}, \cite{zbMATH00903629} and \cite{zbMATH00085903}) gave $p$-adic interpretations of $G$-functions in some special cases, but were limited either by integrality assumptions or applied only to families of elliptic curves.

A second more recent revival of the $G$-function method has been initiated by Daw and Orr with a focus on applications to Zilber-Pink. Their first such works \cite{zbMATH07481643} and \cite{zbMATH07608391} produced relations only at infinite places and used integrality assumptions at finite places, following the arguments in Andr\'e's book \cite{zbMATH00041964}. Their most recent work \cite{daw2022zilber} --- which was a key inspiration for this paper --- uses the $p$-adic Tate uniformization for elliptic curves to produce relations at the finite places and thereby extend the applicability of the method. Following the release of this paper Daw and Orr have written a follow-up to \cite{daw2022zilber} in the more general setting of families of multiplicatively degenerating abelian varieties \cite{daw2023large}. 

The $G$-function method has also been recently taken up by Papas, again with Zilber-Pink in mind. His thesis work \cite{papas2022height} gives a general strategy for applying the $G$-function method beyond the setting of abelian families, and carries out substantially the task of constructing relations at infinite places for various types of endomorphism algebras associated to special moduli. (His thesis work has since been updated to \cite{papas2023unlikely}, which incorporates results from this work. We will later cite the first version \cite{papas2022height} for a certain result concerning absolute Hodge endomorphisms which will be useful for us.) Subsequently, Papas has also studied Zilber-Pink applications of the $G$-function method on curves in $\mathcal{A}_{g}$ \cite{papas2023cases} and in $Y(1)^{n}$ \cite{papas2024effective} \cite{papas2024zilberpink}.

Finally, with respect to the Pila-Zannier strategy for Zilber-Pink-type atypical intersections, Daw and Ren in \cite{daw2018applications} give an approach for the special case of subvarieties of Shimura varieties. The basic idea is essentially the same as ours in that one tries to argue that having ``many'' special points in $S$ will allow one to produce some low-dimension algebraic variety which interacts exceptionally with an analytic period image in order to contradict an Ax-Schanuel principle. Our results are similar, except that we are able to work in a general algebro-geometric setting beyond the case of Shimura varieties, and some aspects of our approach seem simpler to us.

\subsection{Acknowledgements}

We thank Mark Kisin for helpful discussions and suggestions, and for helping streamline the introduction. We thank Jacob Tsimerman, Chris Daw, Martin Orr, and Georgios Papas for comments on a draft of this manuscript. The author would in particular like to thank Chris Daw and Martin Orr for their invitation to discuss his work with them during the summer of 2023, and Martin Orr for alerting the author to some errors in an earlier version of this manuscript.

We thank Donu Arapura for a MathOverflow comment which suggested to the author the idea of considering hypersurface degenerations to hyperplanes as an application of \autoref{genheightthm}; we thank Thomas Scanlon for pointing us to the reference \cite{zbMATH06256440}. Finally, we thank the referees for their numerous helpful comments, suggestions, and corrections.

\subsection{Notational Conventions}
\label{notconvtpt1}

The paper uses roughly two tranches of standard notation. The first tranche consists largely of the notation that has just been introduced in the introduction and will be used throughout \S\ref{cohomprelimsec}, \S\ref{cohomcompsec}, \S\ref{realgfuncsec}, and \S\ref{algrelsec}. The second tranche which will be introduced at the beginning of \S\ref{heightboundsec}, and will then then be in force throughout the rest of the paper. To aid the reader, we collect the notation in the first tranche here.

\begin{itemize}
\item[-] $K$ is a number field, $K_{v}$ is a non-Archimedean completion of $K$ above the prime $p$, $k$ is a discretely valued characteristic zero local field with residue field $\kappa$, and $\widehat{k}$ is the completion of its algebraic closure;
\item[-] $f : X \to S$ is a projective $K$-algebraic family of geometrically connected varieties of dimension $n$, with smooth generic fibre, with $S$ a smooth curve, and $f' : X' \to S'$ is the base-change of $f$ to the locus $S' \subset S$ above which the fibres are smooth;
\item[-] $S \setminus S' = \{ s_{0} \}$, with $s_{0} \in S(K)$ an isolated singular point;
\item[-] $E = X \setminus X' = f^{-1}(S \setminus S')$ is a simple normal crossing divisor;
\item[-] $q \in E(K)$ is a point where $E$ has $\mu$ branches;
\item[-] $z_{1}, \hdots, z_{\nu}$ is an \'etale coordinate system with $\nu = n+1$ defined in some open neighbourhood $U \subset X$ of $q$, such that $E \cap U$ is defined by $z_{1} \cdots z_{\mu} = 0$;
\item[-] $s = z_{1} \cdots z_{\mu}$ is a local uniformizing parameter at $s_{0}$;
\item[-] $w = \mu - 1$ is a fixed integer denoting a cohomological degree;
\item[-] $\mathbb{V}' = R^{w} \an{f}_{*} \mathbb{Z} / \textrm{tor}.$ is the variation of Hodge structure in degree $w$ and whose fibres have dimension $m$, and such that the monodromy around $s_{0}$ is unipotent;
\item[-] $\mathcal{H}' = R^{w} f'_{*} \Omega^{\bullet}_{X'/S'}$ is the associated de Rham cohomology vector bundle; 
\item[-] $\mathcal{H}$ is a vector bundle on $S$ extending $\mathcal{H}'$, which we introduce in \S\ref{canonextsec} below;
\item[-] $\mathcal{H}$ is a free $\mathcal{O}_{S}$-module with basis $\omega_{1}, \hdots, \omega_{m}$; and
\item[-] $h_{1}, \hdots, h_{m}$ are $G$-functions (power series in the parameter $s$).
\end{itemize}

Given an algebraic variety $Y$ over a subfield of $\mathbb{C}$, we write $Y^{\textrm{an}}$ for its analytification. If $Y$ is instead defined over a characteristic zero non-archimedean local field $L$, we write $Y^{\textrm{ad}}$ for the associated adic space. Given a space $Y$ (e.g., a scheme, adic space, or complex analytic variety), a subscript of the form $(-)_{\textrm{\'et}}, (-)_{\textrm{p\'et}}, (-)_{\textrm{ad}}$, etc., will refer to a site associated to that space.

\section{Cohomological Preliminaries}
\label{cohomprelimsec}

We continue with the notation and setup of \S\ref{notconvtpt1}. 

\subsection{A Model for the Canonical Extension}
\label{canonextsec}

We begin by describing an explicit model for the canonical extension $\mathcal{H}$ of $\mathcal{H}'$ referenced in the introduction, following Steenbrink \cite{zbMATH03473825}. We define the de Rham complex $\Omega^{\bullet}_{X}(\log E)$ of algebraic differentials with logarithmic poles along $E$ as follows: for an open set $V \subset X$ the sections of $\Omega^{p}_{X}(\log E)$ over $V$ are the algebraic forms $\omega$ on $V \setminus E$ such that $\omega$ and $d \omega$ have at most a simple pole along $E$. If one chooses local coordinates $(z_{1}, \hdots, z_{\nu})$ around a point $q \in E$ so that $E$ is defined by $z_{1} \cdots z_{\mu} = 0$ for some $1 \leq \mu \leq \nu$, then the stalk $\Omega^{1}_{X}(\log E)_{(q)}$ is a free module over $\mathcal{O}_{X,(q)}$ with generators $dz_{1}/z_{1}, \hdots, dz_{\mu}/z_{\mu}, z_{\mu+1}, \hdots, z_{\nu}$ and $\Omega^{p}_{X}(\log E) = \bigwedge^{p} \Omega^{1}_{X}(\log E)$. We further define $\Omega^{p}_{X/S}(\log E)$ as the $p$'th exterior power of the quotient $\Omega^{1}_{X}(\log E)/f^{*} \Omega^{1}_{S}(\log \{ s_{0} \})$, with $\Omega^{1}_{S}(\log \{ s_{0} \})$ defined analogously via differentials with at most a logarithmic pole at $s_{0}$.

The following is proven in \cite[2.18]{zbMATH03473825} (note that it makes no difference whether one uses the algebraic or analytic site, c.f. \S\ref{cohcohomsec} below):

\begin{prop}
For all $w \geq 0$, the sheaf $R^{w}f_{*}(\Omega^{\bullet}_{X/S}(\log E))$ is locally free on $S$ and for all $s \in S$ the canonical map
\[ R^{w} f_{*}(\Omega^{\bullet}_{X/S}(\log E)) \otimes_{\mathcal{O}_{S}} (\mathcal{O}_{S,s}/\mathfrak{m}_{S,s}) \to H^{w}(X_{s}, \Omega^{\bullet}_{X/S}(\log E) \otimes_{\mathcal{O}_{X}} \mathcal{O}_{X_{s}}) \]
is an isomorphism. \qed
\end{prop}

We may therefore take $\mathcal{H} = R^{w} f_{*} \Omega^{\bullet}_{X/S}(\log E)$. Let us now consider the setup in the introduction, where $U \subset X$ was a fixed affine Zariski open subset with coordinates $z_{1}, \hdots, z_{\nu}$ trivializing $\Omega^{1}_{U}$ and such that the map $U \to S$ takes the form $s \mapsto z_{1} \cdots z_{\mu}$ for some $1 \leq \mu \leq \nu$. As before, we set $w = \mu-1$. We have a natural map $R^{w} f_{*} \Omega^{\bullet}_{X/S}(\log E) \to R^{w} f_{*} \Omega^{\bullet}_{U/S}(\log (E \cap U))$. For $S$ affine, the cohomology module $R^{w} f_{*} \Omega^{\bullet}_{U/S}(\log (E \cap U))$ may be identified with $R^{w} \Gamma \, \Omega^{\bullet}_{U/S}(\log (E \cap U))$, as follows from the Leray spectral sequence. Moreover, because $U$ is affine, this can in turn be identified with the cohomology in degree $w$ of the complex
\[ 0 \to \mathcal{O}_{U} \to \Omega^{1}_{U/S}(\log (E \cap U)) \to \cdots \to \Omega^{n}_{U/S}(\log (U \cap E)) , \]
viewed as a module over $\mathcal{O}_{S}$. Restricting to the completed stalk at a point $q$ in the locus $z_{1} = \cdots = z_{\mu} = 0$ mapping to $s_{0}$ one obtains a complex of $\widehat{\mathcal{O}}_{S,(s_{0})}$-modules
\[ 0 \to \widehat{\mathcal{O}}_{U,(q)} \to \widehat{\Omega}^{1}_{U/S}(\log (E \cap U))_{(q)} \to \cdots \to \widehat{\Omega}^{n}_{U/S}(\log (E \cap U))_{(q)} , \]
and by composition a restriction map
\[ \eta : \Gamma(R^{w} f_{*} \Omega^{\bullet}_{X/S}(\log E), S) \to \textrm{Cohom}^{w} \left[ \widehat{\Omega}^{\bullet}_{U/S}(\log(E \cap U))_{(q)} \right] , \]
where we denote by $\textrm{Cohom}^{w}$ the na\"ive cohomology in degree $w$. This map will be used in the proof of \autoref{bigsupergfuncthm} to construct $G$-functions.

The following lemma describes the form of the elements in 
\[ \textrm{Cohom}^{w} \left[ \widehat{\Omega}^{\bullet}_{U/S}(\log(E \cap U))_{(q)} \right] , \]
and also the analogous modules obtained by considering convergent power series in the complex and rigid-analytic topologies.

\begin{lem}
\label{uniquereplem}
Suppose that $A$ is any of the rings 
\[ \{ K[[z_{1}, \hdots, z_{\nu}]], \mathbb{C} \{ z_{1}, \hdots, z_{\nu} \}, \mathbb{C} \{ z_{1}, \hdots, z_{\nu} \}_{\sqrt{\mu} r^{1/\mu}},  k\langle z_{1}, \hdots, z_{\nu} \rangle_{r^{1/\mu}} \} \]
which are (in order) formal power series over the characteristic zero field $K$, germs of complex analytic power series, complex analytic power series convergent on the open ball of radius $\sqrt{\mu} r^{1/\mu}$, and power series convergent in the open non-archimedian ball of radius $r^{1/\mu}$, respectively. Let $B$ be the corresponding ring in $\{ K[[s]], \mathbb{C} \{ s \}, \mathbb{C} \{ s \}_{r}, k\langle s \rangle_{r} \}$ and consider the map $B \to A$ given by $s \mapsto z_{1} \cdots z_{\mu}$. Consider the complex
\begin{equation}
\label{stdrelcomp}
 0 \to A \to \Omega^{1}_{A/B}(\log E) \to \cdots \to \bigwedge^{\nu} \Omega^{1}_{A/B}(\log E) \to 0 , 
\end{equation}
where $\Omega^{1}_{A/B}(\log E)$ is the quotient of $\Omega^{1}_{A}(\log E)$ and $A \otimes_{B} \Omega^{1}_{B}(\log \{ 0 \})$, with
\[ \Omega^{1}_{A}(\log E) = A \frac{dz_{1}}{z_{1}} \oplus \cdots \oplus A \frac{dz_{\mu}}{z_{\mu}} \oplus A dz_{\mu+1} \oplus \cdots \oplus A dz_{\nu} , \]
and $\Omega^{1}_{B}(\log \{ 0 \}) = B \frac{ds}{s}$. Then every element $\alpha$ of $\textrm{Cohom}^{w} \left[ \Omega^{\bullet}_{A/B}(\log E) \right]$ admits a unique representation of the form 
\[ h \frac{dz_{2} \cdots dz_{\mu}}{z_{2} \cdots z_{\mu}} \]
where $h$ is a uniquely determined element of $B$.
\end{lem}

\begin{rem}
We note that the fact that the complex open ball of radius $\sqrt{\mu} r^{1/\mu}$ maps into the ball of radius $r$ follows from the arithmetic-geometric mean inequality:
\[ |z_{1}| \cdots |z_{\mu}| \leq \left( \frac{|z_{1}| + \cdots + |z_{\mu}|}{\mu} \right)^{\mu} \leq \left( \frac{\sqrt{|z_{1}|^2 + \cdots + |z_{\nu}|^2}}{\sqrt{\mu}} \right)^{\mu} < r \]
\end{rem}

\begin{proof}
In the complex analytic setting the entire cohomology of the complex (\ref{stdrelcomp}) is described in \cite[1.13]{zbMATH03473825}, and the same proof works in general. For the convenience of the reader we give some details. The relation $s = z_{1} \cdots z_{\mu}$ induces the relation $\frac{dz_{1}}{z_{1}} + \cdots + \frac{dz_{\mu}}{z_{\mu}} = 0$ in the complex (\ref{stdrelcomp}), which gives a natural presentation of the complex in terms of the forms $\frac{dz_{2}}{z_{2}}, \hdots, \frac{dz_{\mu}}{z_{\mu}}, dz_{\mu+1}, \hdots, dz_{\nu}$ only. The complex then reduces to a Kozul-type complex $L^{\bullet}$ generated by
\begin{align*}
A \frac{dz_{2}}{z_{2}} \oplus \cdots \oplus A \frac{dz_{\mu}}{z_{\mu}} \oplus A \, dz_{\mu+1} \oplus \cdots \oplus A \, dz_{\nu},
\end{align*}
where the differential operators for $L^{\bullet}$ are given by $D_{i} = z_{i} \partial_{i} - z_{1} \partial_{1}$ for $2 \leq i \leq \mu$ and $D_{i} = \partial_{i}$ for $\mu+1 \leq i \leq \nu$. 

Suppose that an element 
\[ \beta = g \frac{dz_{i_{1}}}{z_{i_{1}}^{e_{i_{1}}}} \wedge \cdots \wedge \frac{dz_{i_{r}}}{z_{i_{1}}^{e_{i_{r}}}} \]
in the complex $L^{\bullet}$ lies in the kernel of the differential, with $g \in A$ a monomial, and where the exponents $e_{i_{j}}$ satisfy the property that $e_{i_{j}} \in \{ 0, 1 \}$ and $e_{i_{j}} = 0$ if $i_{j} > \mu$. Then from the construction of the Kozul complex we must have that $D_{j}(g) = 0$ for each $j$ not appearing in the set $\{ i_{1}, \hdots, i_{r} \}$. If we have $i_{k} > \mu$ for some $k$ (and hence $e_{i_{k}} = 0$), and the variable $z_{i_{k}}$ occurs with exponent $a \geq 0$ in the monomial $g$, we compute that
\begin{align*}
& \, d \left( (-1)^{i_{k}} \frac{z_{i_{k}}}{a+1} g \frac{dz_{i_{1}}}{z_{i_{1}}^{e_{i_{1}}}} \wedge \cdots \wedge \widehat{dz_{i_{k}}} \wedge \cdots \wedge \frac{dz_{i_{r}}}{z_{i_{1}}^{e_{i_{r}}}} \right) \\ = & \, \beta + (-1)^{i_{k}} \frac{z_{i_{k}}}{a+1} \underbrace{\sum_{j} \, D_{j}(g) \, dz_{j} \wedge \left( \frac{dz_{i_{1}}}{z_{i_{1}}^{e_{i_{1}}}} \wedge \cdots \wedge \frac{dz_{i_{r}}}{z_{i_{1}}^{e_{i_{r}}}} \right)}_{=0}
\end{align*}
which shows that $\beta$ can be integrated. Thus, non-trivial contributions to cohomology appear only when $\{ i_{1}, \hdots, i_{r} \} \subset \{ 2, \hdots, \mu \}$. In degree $r = w$ this means that $\{ i_{1}, \hdots, i_{w} \} = \{ 2, \hdots, \mu \}$ and that $D_{j}(g) = 0$ for $j > \mu$, meaning we can assume our class $\alpha$ is of the form $h \frac{dz_{2} \cdots dz_{\mu}}{z_{2} \cdots z_{\mu}}$ with $h$ depending on $z_{1}, \hdots, z_{\mu}$ only. One checks that a monomial $m$ in $h$ cannot be integrated if and only if $D_{j}(m) = 0$ for $2 \leq j \leq \mu$. But this means that $h$ consists only of terms like $c \cdot (z_{1} \cdots z_{\mu})^{e}$, hence $g$ lies in the image of $B$. 

The uniqueness claim is checked directly from the construction of the complex, as no non-zero elements of the specified form lie in the image of the differential.
\end{proof}

\subsection{\v{C}ech cohomological recollections}
\label{cechrecsec}

\subsubsection{\v{C}ech cohomology of complexes}
\label{cechcompsec}

We now develop the general formalism of the \v{C}ech double complex associated to a complex $(\mathcal{F}^{\bullet}, d_{\mathcal{F}})$ of sheaves on a site $\mathcal{C}$, valued in an abelian category $\mathcal{A}$, generalizing the case of sheaves on a topological space which appears in \cite[\href{https://stacks.math.columbia.edu/tag/01ED}{Section 01ED}]{stacks-project} and \cite[\href{https://stacks.math.columbia.edu/tag/01FP}{Section 01FP}]{stacks-project}. 

We assume that $\mathcal{C}$ has a final object $X$, and we let $\mathcal{U} = \{ c_{i} : U_{i} \to X \}_{i \in I}$ be a covering of $X$. We first consider the case where $\mathcal{F}^{\bullet} = \mathcal{F}$ consists of a single sheaf. We define 
\[ C^{p}(\mathcal{U}, \mathcal{F}) = \prod_{(i_{0}, \hdots, i_{p}) \in I^{p+1}} \mathcal{F}(U_{i_{0}} \times_{X} \cdots \times_{X} U_{i_{p}}) . \]
Given $s \in C^{p}(\mathcal{U}, \mathcal{F})$ we will write $s_{i_{0} \hdots i_{p}}$ for its value in the factor $\mathcal{F}(U_{i_{0}} \times_{X} \hdots \times_{X} U_{i_{p}})$, and we define the differential
\[ \delta : C^{p}(\mathcal{U}, \mathcal{F}) \to C^{p+1}(\mathcal{U}, \mathcal{F})  \]
by the formula
\[ \delta(s)_{i_{0} \hdots i_{p+1}} = \sum_{j = 0}^{p+1}(-1)^{j} \restr{s_{i_{0} \hdots \widehat{i_{j}} \hdots i_{p+1}}}{U_{i_{0}} \times_{X} \cdots \times_{X} U_{i_{p+1}}} \]
where restriction comes from the natural fibre product projection. One checks that $(C^{\bullet}(\mathcal{U}, \mathcal{F}), \delta)$ is a complex. 

The formation of the complex $C^{\bullet}(\mathcal{U}, \mathcal{F})$ is functorial in $\mathcal{F}$, so given a complex $(\mathcal{F}^{\bullet}, d_{\mathcal{F}})$ one naturally obtains a double complex $C^{\bullet}(\mathcal{U}, \mathcal{F}^{\bullet})$. We write $(L^{\bullet}(\mathcal{U}, \mathcal{F}^{\bullet}), d)$ for the associated total complex, with terms
\[ L^{n}(\mathcal{U}, \mathcal{F}^{\bullet}) = \bigoplus_{p+q=n} \prod_{i_{0} \hdots i_{p}} \mathcal{F}^{q}(U_{i_{0} \hdots i_{p}}) \]
and with differential of an element $\alpha$ of degree $n$ given by $d = \delta + (-1)^{p+1} d_{\mathcal{F}}$. Finally, we write $\ch{H}^{p}(\mathcal{U}, \mathcal{F}^{\bullet})$ for the cohomology groups of this complex. 

We now compare the \v{C}ech cohomology to sheaf cohomology. We denote the cohomology of a complex of sheaves $\mathcal{F}^{\bullet}$ computed on an object $V$ of $\mathcal{C}$ using covers of $\mathcal{C}$ by $H^{\bullet}(V_{\mathcal{C}}, \mathcal{F}^{\bullet})$. Then we have the following generalization of \cite[III, Ch.0, 12.4.6]{zbMATH03193326}:

\begin{prop}
\label{specseqprop}
Let $\mathcal{F}^{\bullet}$ be a bounded below complex of sheaves on $\mathcal{C}$, valued in an abelian category $\mathcal{A}$, and let $\mathcal{U} = \{ c_{i} : U_{i} \to X \}_{i \in I}$ be a covering of $X$. Then there exists a spectral sequence abutting to $H^{\bullet}(X_{\mathcal{C}}, \mathcal{F}^{\bullet})$ whose second page is given by
\[ E^{pq}_{2} = \textrm{Cohom}^{p}(L^{\bullet}(\mathcal{U}, J^{q}(\mathcal{F}^{\bullet}))) , \]
where $J^{q}(\mathcal{F}^{\bullet})$ denotes the complex of presheaves whose $j$'th term is given by $[V \mapsto H^{q}(V_{\mathcal{C}}, \mathcal{F}^{j})]$. 
\end{prop}

\begin{proof}
One just has to check that all the steps in the argument in \cite[III, Ch.0, 12.4.6]{zbMATH03193326} generalize to this situation (c.f. \cite[\href{https://stacks.math.columbia.edu/tag/08BN}{Lemma 08BN}]{stacks-project}). We consider a Cartan-Eilenberg resolution $\mathcal{L}^{\bullet \bullet}$ of $\mathcal{F}^{\bullet}$ by injective sheaves, constructed as in \cite[\href{https://stacks.math.columbia.edu/tag/015I}{Lemma 015I}]{stacks-project}. From the functorality of \v{C}ech cohomology we obtain a tricomplex $C^{\bullet}(\mathcal{U}, \mathcal{L}^{\bullet \bullet}) = [C^{i}(\mathcal{U}, \mathcal{L}^{jk})]$ which we may regard as a bicomplex in degrees $i$ and $j+k$. Because the sheaves $\mathcal{L}^{jk}$, and hence the terms in the total complex of $\mathcal{L}^{\bullet\bullet}$, are all injective sheaves of abelian groups, the \v{C}ech complex $C^{\bullet}(\mathcal{U},\mathcal{L}^{\bullet \bullet})$ computes the cohomology of the total complex associated to $\mathcal{L}^{\bullet\bullet}$: this follows by combining \cite[\href{https://stacks.math.columbia.edu/tag/03AW}{Lemma 03AW}]{stacks-project}, which shows that the positive degree \v{C}ech cohomology on $\mathcal{U}$ of each injective sheaf is zero, with \cite[\href{https://stacks.math.columbia.edu/tag/0133}{Lemma 0133}]{stacks-project} applied to the map $\mathcal{L}^{\bullet \bullet} \to C^{\bullet}(\mathcal{U},\mathcal{L}^{\bullet \bullet})$, where we regard $\mathcal{L}^{\bullet \bullet}$ as a single complex in $j+k$ (its total complex), as stated. Because the total complex of $\mathcal{L}^{\bullet \bullet}$ computes the cohomology of $\mathcal{F}^{\bullet}$, it follows that the map $\mathcal{F}^{\bullet} \to C^{\bullet}(\mathcal{U},\mathcal{L}^{\bullet \bullet})$ induces an isomorphism on cohomology.

We now consider the tricomplex $C^{\bullet}(\mathcal{U}, \mathcal{L}^{\bullet \bullet})$ as a bicomplex in degrees $i+j$ and $k$. Then because $\mathcal{L}^{j,\bullet}$ is an injective resolution of $\mathcal{F}^{j}$ for all $j$, the degree $q$ cohomology of the complex $C^{i}(\mathcal{U}, \mathcal{L}^{j,\bullet})$ is then given by the \v{C}ech complex $C^{i}(\mathcal{U}, J^{q}(\mathcal{F}^{j}))$. Computing the second page then gives the result. 
\end{proof}

\begin{cor}
\label{Cechcomputescohomcor}
Suppose that for each $U'$ obtained as a fibre product of objects in the cover $\mathcal{U}$ and for each $k$ one has that $H^{q}(U'_{\mathcal{C}}, \mathcal{F}^{k}) = 0$ for all $q > 0$. Then $L^{\bullet}(\mathcal{U}, \mathcal{F}^{\bullet})$ computes the cohomology of the complex $\mathcal{F}^{\bullet}$.
\end{cor}

\begin{proof}
The assumption ensures that the cohomology presheaves $J^{q}$ appearing in \autoref{specseqprop} are zero on all open neighbourhoods obtained by taking fibre products of elements in $\mathcal{U}$, and hence the total complexes $L^{\bullet}(\mathcal{U}, J^{q}(\mathcal{F}^{\bullet}))$ are zero for $q > 0$. Thus the only non-zero terms in the spectral sequence appearing on the second page are those when $q = 0$, which correspond to the cohomology groups $E^{p0}_{2} = \textrm{Cohom}^{p}(L^{\bullet}(\mathcal{U}, \mathcal{F}^{\bullet}))$. Because the spectral sequence degenerates at the second page, it follows that $E^{p0}_{2}$ is nothing other than the cohomology $H^{p}(X_{\mathcal{C}}, \mathcal{F}^{\bullet})$.
\end{proof}

\subsubsection{Cup product in \v{C}ech cohomology}
\label{cechcupprodsec}

We also recall how to define the cup product on the \v{C}ech complex, following \cite[\href{https://stacks.math.columbia.edu/tag/01FP}{Section 01FP}]{stacks-project} in the setting of complexes of sheaves on topological spaces. Given two complexes of sheaves $\mathcal{F}^{\bullet}$ and $\mathcal{G}^{\bullet}$ of abelian groups on the site $\mathcal{C}$, we write $\Tot(\mathcal{F}^{\bullet} \otimes \mathcal{G}^{\bullet})$ for the complex with terms $\bigoplus_{p+q=n} \mathcal{F}^{p} \otimes \mathcal{G}^{q}$ and where the differential is given by
\[ d(\alpha \otimes \beta) = d(\alpha) \otimes \beta + (-1)^{\deg(\alpha)} \alpha \otimes d(\beta) . \]
Given a covering $\mathcal{U} = \{ c_{i} : U_{i} \to X \}_{i \in I}$, our cup product is then a map
\[ \cup : \Tot\left(\Tot(C^{\bullet}(\mathcal{U}, \mathcal{F}^{\bullet})) \otimes \Tot(C^{\bullet}(\mathcal{U}, \mathcal{G}^{\bullet}))\right) \to \Tot(C^{\bullet}(\mathcal{U}, \Tot(\mathcal{F}^{\bullet} \otimes \mathcal{G}^{\bullet}))) . \]
It is given by the rule 
\[ (\alpha \cup \beta)_{i_{0} \hdots i_{p}} = \sum_{r=0}^{p} \ep(\deg \alpha,\deg \beta,p,r) \alpha_{i_{0} \hdots i_{r}} \otimes \beta_{i_{r} \hdots i_{p}} , \]
where $\ep(n,m,p,r) = (-1)^{(p+r)n + rp + r}$. The associativity of the cup product as well as the identity
\[ d(\alpha \cup \beta) = d(\alpha) \cup \beta + (-1)^{\deg(\alpha)} \alpha \cup d(\beta) \]
may be proved by explicit calculation, exactly as is done in \cite[\href{https://stacks.math.columbia.edu/tag/01FP}{Section 01FP}]{stacks-project} in the setting of complexes on topological spaces. Moreover, the cup product is compatible with a graded commutative structure on the complex $\mathcal{F}^{\bullet}$, as we now explain, following \cite[\href{https://stacks.math.columbia.edu/tag/01FP}{Section 01FP}]{stacks-project}.

Suppose that we have a graded commutative multiplication map
\[ \wedge^{\bullet} : \Tot(\mathcal{F}^{\bullet} \otimes \mathcal{F}^{\bullet}) \to \mathcal{F}^{\bullet} . \]
This is defined to mean that given sections $s$ of $\mathcal{F}^{a}$ and $t$ of $\mathcal{F}^{b}$ we obtain a section $s \wedge t$ of $\mathcal{F}^{a+b}$ in such a way that $s \wedge t = (-1)^{ab} t \wedge s$, and that we have $d(s \wedge t) = d(s) \wedge t + (-1)^{a} s \wedge d(t)$. We may then consider the composition
\begin{align*}
\Tot\left(\Tot(C^{\bullet}(\mathcal{U}, \mathcal{F}^{\bullet})) \otimes \Tot(C^{\bullet}(\mathcal{U}, \mathcal{F}^{\bullet}))\right) &\xrightarrow{\cup} \Tot(C^{\bullet}(\mathcal{U}, \Tot(\mathcal{F}^{\bullet} \otimes \mathcal{F}^{\bullet}))) \\
&\xrightarrow{\wedge} \Tot(C^{\bullet}(\mathcal{U}, \mathcal{F}^{\bullet})) .
\end{align*}
It may be checked as in \cite[\href{https://stacks.math.columbia.edu/tag/01FP}{Section 01FP}]{stacks-project} that this induces a map on \v{C}ech cohomology
\[ H^{n}(\Tot(C^{\bullet}(\mathcal{U}, \mathcal{F}^{\bullet}))) \times H^{m}(\Tot(C^{\bullet}(\mathcal{U}, \mathcal{F}^{\bullet}))) \to H^{n+m}(\Tot(C^{\bullet}(\mathcal{U}, \mathcal{F}^{\bullet}))) . \]
In our situation of interest, this will reproduce the cup product on both \'etale cohomology and (algebraic) de Rham cohomology.

\subsection{The pro-\'etale site}
\label{proetsec}

Let us fix an adic space $X$ over $\Spa(k, \mathcal{O}_{k})$. We will assume that $X$ is locally noetherian. (This assumption will also continue to be in force in subsequent sections without further comment.) 

We will begin by defining some categories (and sites) associated to $X$. First, one has the \'etale site $X_{\textrm{\'et}}$, whose objects consist of \'etale maps $U \to X$ of adic spaces and morphisms between them. Next we consider the category $\textrm{pro}(X_{\textrm{\'et}})$: its objects consist of projective limits $\varprojlim_{i \in I} U_{i}$ of objects of $X_{\textrm{\'et}}$ and its morphisms are the natural morphisms of limit diagrams. A map of objects $U \to V$ in the category $\textrm{pro}(X_{\textrm{\'et}})$ is called \'etale if it is induced by an \'etale morphism of objects $U_{0} \to V_{0}$ in $X_{\textrm{\'et}}$. A map of objects $U \to V$ is called pro-\'etale if we have $U = \lim_{i} U_{i}$ in such a way so that $U \to V$ is given by an inverse limit $U_{i} \to V$ of objects \'etale over $V$, and such that $U_{i} \to U_{j}$ is finite \'etale and surjective for large $i > j$. The category $X_{\textrm{p\'et}}$ is then defined to be the full subcategory of $\textrm{pro}(X_{\textrm{\'et}})$ consisting of objects which are pro-\'etale over $X$. A covering $\mathcal{U}$ in $\textrm{pro}(X_{\textrm{\'et}})$ of an object $U$ is given by a family of pro-\'etale morphisms $\mathcal{U} = \{ f_{i} : U_{i} \to U \}$ such that $|U| = \bigcup_{i} f_{i}(|U_{i}|)$, where we give pro-objects the limit topology. By \cite[Lemma 3.10]{zbMATH06209107}, this defines a site.

If one instead starts with the category $X_{\textrm{f\'et}}$ of objects finite \'etale over $X$, one may carry out the analogous procedure to define a category $X_{\textrm{pf\'et}}$, which we call the ``pro-finite finite \'etale site''. It is naturally a subcategory of $X_{\textrm{p\'et}}$, and we have a natural map of sites $X_{\textrm{p\'et}} \to X_{\textrm{pf\'et}}$.

The pro-finite \'etale site can be used to compute \'etale cohomology with coefficients in $\mathbb{Z}_{p}$ as the cohomology of the sheaf $\hat{\mathbb{Z}}_{p}(U) = \Hom_{\textrm{cont}}(|U|,\mathbb{Z}_{p})$, where we consider continuous morphisms of the underlying topological spaces, and $\mathbb{Z}_{p}$ has the usual $p$-adic topology.

Next we introduce some important sheaves on $X_{\textrm{p\'et}}$, following \cite{zbMATH06209107}. The first is $\mathcal{O}_{X}$, the ``uncompleted structure sheaf'', which is the pullback $\gamma^{*} \mathcal{O}_{X_{\textrm{\'et}}}$ under the natural map $\gamma : X_{\textrm{p\'et}} \to X_{\textrm{\'et}}$ of sites. Likewise we have the subring of integral elements $\mathcal{O}^{+}_{X} = \gamma^{*} \mathcal{O}^{+}_{X_{\textrm{\'et}}}$. These sheaves can then be completed to obtain $\hat{\mathcal{O}}^{+}_{X} = \varprojlim \mathcal{O}^{+}_{X} / p^{n}$ and $\hat{\mathcal{O}}_{X} = \hat{\mathcal{O}}^{+}_{X}\left[\frac{1}{p}\right]$. Next we have the tilted integral structure sheaf, defined as $\hat{\mathcal{O}}^{+}_{X^{\flat}} = \varprojlim_{\Phi} \mathcal{O}^{+}_{X}/p$, with the inverse limit over Frobenius. If $k = \widehat{k}$ we set $\hat{\mathcal{O}}_{X^{\flat}} = \hat{\mathcal{O}}^{+}_{X^{\flat}} \otimes_{k^{\flat+}} k^{\flat}$; here we use the notion of the tilt $X^{\flat}$ of $X$ from Scholze's theory of perfectoid spaces \cite{zbMATH06320533}. We then define $\mathbb{A}_{\textrm{inf}} = W(\hat{\mathcal{O}}^{+}_{X^{\flat}})$ and $\mathbb{B}_{\textrm{inf}} = \mathbb{A}_{\textrm{inf}}\left[\frac{1}{p}\right]$. 

We have a natural map $\theta : \mathbb{A}_{\textrm{inf}} \to \widehat{\mathcal{O}}^{+}_{X}$ which extends to a map $\mathbb{B}_{\textrm{inf}} \to \hat{\mathcal{O}}_{X}$. To define it, we work locally, where the sheaf $\mathbb{A}_{\textrm{inf}}$ is represented by a ring $W(A^{\flat})$, and we wish to construct a map $W(A^{\flat}) \to A$. We may represent an element $x \in W(A^{\flat})$ via its Witt vector components as a sum $\sum_{i} p^{i}[x_{i}]$. We then define
\[ \theta\left(\sum_{i} p^{i}[x_{i}]\right) = \sum_{i} p^{i} x_{i}^{\sharp} , \]
where the operation $(-)^{\sharp}$ is defined on $y \in A^{\flat}$, represented by the sequence $(y_{1}, y_{2}, \hdots)$, by choosing lifts $\hat{y_{j}}$ for all $j$ and setting
\[ y^{\sharp} = \lim_{j \to \infty} \hat{y_{j}}^{p^{j}} . \]
We then define $\mathbb{B}^{+}_{\textrm{dR}} = \varprojlim \mathbb{B}_{\textrm{inf}} / (\ker \theta)^{n}$ and $\mathbb{B}_{\textrm{dR}} = \mathbb{B}^{+}_{\textrm{dR}}[t^{-1}]$, where $t$ is any element generating the kernel of $\theta$. 

Finally we define $\mathcal{O}\mathbb{B}_{\textrm{inf}} = \mathcal{O}_{X} \otimes_{W(\kappa)} \mathbb{B}_{\textrm{inf}}$. The map $\theta$ on $\mathbb{B}_{\textrm{inf}}$ extends to a map $\theta : \mathcal{O} \mathbb{B}_{\textrm{inf}} \to \hat{\mathcal{O}}_{X}$. One then defines the de Rham structure sheaf $\mathcal{O} \mathbb{B}^{+}_{\textrm{dR}}$ as in \cite{zbMATH06700126}, correcting the definition in \cite{zbMATH06209107}, and $\mathcal{O} \mathbb{B}_{\textrm{dR}} = \mathcal{O} \mathbb{B}^{+}_{\textrm{dR}}[t^{-1}]$, where $t$ is a generator of $\ker \theta$ (this makes sense locally, as is checked in \cite[\S6]{zbMATH06209107}). Lastly we define $\mathbb{\Omega}^{i}_{X} = \mathcal{O}\mathbb{B}_{\textrm{dR}} \otimes_{\mathcal{O}_{X}} \Omega^{i}_{X}$ as sheaves on $X_{\textrm{p\'et}}$.

\subsection{Coherent Cohomology on Various Sites}
\label{cohcohomsec}

An important sort of fact that we will use (often implicitly) throughout the paper is that it ``doesn't matter'' on which site one computes the cohomology of coherent objects associated to a space $X$. What is meant by this is that if one has two sites associated to $X$, say $X_{1}$ and $X_{2}$, with a natural map of ringed sites $\tau : X_{2} \to X_{1}$, and given a complex of coherent sheaves $\mathcal{F}^{\bullet}$ on $X_{1}$ the natural map $H^{i}(X_{1}, \mathcal{F}^{\bullet}) \to H^{i}(X_{2}, \tau^{*} \mathcal{F}^{\bullet})$ is an isomorphism for all $i$. Note that one typically only needs to check this when $\mathcal{F}^{\bullet} = \mathcal{F}$ is a single sheaf rather than a complex of such. The reason is that, in the situations of interest, the sites $X_{1}$ and $X_{2}$ will contain certain types of objects on which the coherent cohomology of any individual coherent sheaf vanishes (e.g., affine, affinoid, Stein, etc.), and using this fact for each $\mathcal{F}^{i}$ in the complex $\mathcal{F}^{\bullet}$ and an appropriate cover one learns that the ``same'' \v{C}ech complex computes cohomology on both $X_{1}$ and $X_{2}$, and the resulting fact is formal. 

These facts we will only need for complexes of differentials (possibly relative, possibly logarithmic) and for sufficiently nice spaces $X$. (Our spaces or maps of spaces will also often be proper, which makes things even easier.) Nevertheless, we give some of the required facts in greater generality.

\begin{prop}
\label{cohcohomprop}
For $\mathcal{F}^{\bullet}$ a complex of coherent sheaves on $X_{1}$ and $\tau : X_{2} \to X_{1}$ a map of sites, the natural map $H^{i}(X_{1}, \mathcal{F}^{\bullet}) \to H^{i}(X_{2}, \tau^{*} \mathcal{F}^{\bullet})$ is an isomorphism when
\begin{itemize}
\itemsep0em
\item[(i)] $X$ is a scheme, $X_{1} = X_{\textrm{Zar}}$, $X_{2} = X_{\textrm{\'et}}$;
\item[(ii)] $X$ is a smooth locally finite-type $\mathbb{C}$-scheme, $X_{1} = X_{\textrm{Zar}}$, $X_{2} = X_{\textrm{an}}$, and $\mathcal{F}^{\bullet} = \Omega^{\bullet}_{X}$;
\item[(iii)] $X$ is a locally finite-type proper $K_{v}$-scheme, $\overline{K_{v}}$-scheme or $\mathbb{C}_{p}$-scheme, $X_{1} = X_{\textrm{Zar}}$, $X_{2} = X_{\textrm{ad}}$, and $\mathcal{F}^{\bullet} = \Omega^{\bullet}_{X}$;
\item[(iv)] $X$ is a rigid space, $X_{1} = X_{\textrm{ad}}$ is the adic site, $X_{2} = X_{\textrm{\'et}}$;
\item[(v)] $X$ is an adic space, $X_{1} = X_{\textrm{\'et}}$, $X_{2} = X_{\textrm{p\'et}}$.
\end{itemize}
\end{prop}

\begin{proof}
For (i) see \cite[\href{https://stacks.math.columbia.edu/tag/03DW}{Proposition 03DW}]{stacks-project}; for (ii) see the introduction to \cite{zbMATH03234150}; (iii) is a consequence of rigid analytic GAGA, see for instance the appendix to \cite{zbMATH05145717}; for (iv) see \cite[Example 2.1.3]{zbMATH05610519}; for (v) see \cite[Corollary 3.17]{zbMATH06209107}.
\end{proof}

\subsection{The \'etale fundamental group and cohomology}
\label{fundgpsec}

We now describe the \'etale fundamental group of an adic space $X$ and its relation to the cohomology of $X$. We write $X_{\textrm{f\'et}}$ for the category of adic spaces $Y$ which are finite \'etale over $X$. Fixing a geometric point $\overline{x}$ of $X$, we obtain a natural fibre functor $F_{X,\overline{x}} : X_{\textrm{f\'et}} \to \textrm{Set}$, and as usual the group $\pi_{1}^{\textrm{\'et}}(X,\overline{x})$ is defined as the group of automorphisms of this functor. 

For any finite abelian group $\Lambda$, we now describe a natural isomorphism
\[ \textrm{Hom}(\pi^{\textrm{\'et}}_{1}(X,\overline{x}), \Lambda) \xrightarrow{\sim} H^{1}(X_{\textrm{p\'et}}, \Lambda) . \]
We note that to compute $H^{1}(X_{\textrm{p\'et}}, \Lambda)$ for a finite abelian group $\Lambda$ it suffices to use the usual \'etale site $X_{\textrm{\'et}}$, since the natural map $H^{1}(X_{\textrm{p\'et}}, \Lambda) \to H^{1}(X_{\textrm{\'et}}, \Lambda)$ induced by the map of sites $X_{\textrm{\'et}} \to X_{\textrm{p\'et}}$ is an isomorphism; this is due to  \cite[Corollary 3.17]{zbMATH06209107}, as mentioned above.

The description of this isomorphism is essentially identical to the case of schemes, for which \cite[I.\S11]{milneLEC} is a reference. We will give some details. In what follows we also denote by $\Lambda$ the constant sheaf on $X_{\textrm{\'et}}$ it defines, and we use multiplicative notation for group multiplication. A sheaf $\mathcal{L}$ of sets on $X_{\textrm{\'et}}$ on which $\Lambda$ acts is called a \emph{torsor} for $\Lambda$ if:
\begin{itemize}
\item[(i)] there exists a covering $\mathcal{U} = \{ c_{i} : U_{i} \to X \}_{i \in I}$ in $X_{\textrm{\'et}}$ such that $\mathcal{L}(U_{i}) \neq \varnothing$ for all $i$; and
\item[(ii)] for every object $U \to X$ in $X_{\textrm{\'et}}$ and $s \in \mathcal{L}(U)$ the map $\restr{\Lambda}{U} \to \restr{\mathcal{L}}{U}$ given by $g \mapsto g s$ is an isomorphism of sheaves over $U_{\textrm{\'et}}$.
\end{itemize}
A covering $\mathcal{U} = \{ c_{i} : U_{i} \to X \}_{i \in I}$ for which (i) holds is said to split $\mathcal{L}$. Supposing we have such a covering, we construct a \v{C}ech cocycle $[\mathcal{L}] \in H^{1}(\mathcal{U}, \Lambda)$ as follows. Choose some sections $s_{i} \in \mathcal{L}(U_{i})$ for each $i$. By (ii), on each ``intersection'' $U_{ij}$ arising from the cover $\mathcal{U}$ there exists a unique element $\lambda_{ij} \in \Lambda(U_{ij})$ such that $\lambda_{ij} \cdot \restr{s_{i}}{U_{ij}} = \restr{s_{j}}{U_{ij}}$. Then $(\lambda_{ij})_{I \times I}$ is a cocycle, and defines a class $[\mathcal{L}]$ in $H^{1}(\mathcal{U}, \Lambda)$. Moreover we have:

\begin{lem}
\label{torsorcechbij}
The map $\mathcal{L} \mapsto [\mathcal{L}]$ defines a bijection from the set of isomorphism classes of torsors for $\Lambda$ split by $\mathcal{U}$ to $H^{1}(\mathcal{U}, \Lambda)$.
\end{lem}

\begin{proof}
In the case of the \'etale site of a scheme this is \cite[I. Prop 11.1]{milneLEC}, and the proof is identical in our case.
\end{proof}

We now use the fact that there is a further bijection
\begin{equation}
\label{HomLambda}
\{ \textrm{isom. classes of }\Lambda\textrm{-torsors} \} \longleftrightarrow \Hom(\pi^{\textrm{\'et}}_{1}(X,\overline{x}), \Lambda) .
\end{equation}
This is true in a great deal of generality by the work of \cite{zbMATH03290088} (c.f. the discussion in \cite[\S9]{adictamesite}). We describe this correspondence in the special case where the $\Lambda$-torsor $\mathcal{L}$ is representable by a Galois covering $Y \to X$. More precisely, we assume that $\Lambda = \Aut_{X}(Y)$ and that $\mathcal{L}(U) = \Hom_{X}(U, Y)$ for every $U \in X_{\textrm{f\'et}}$, with the natural action of $\Lambda$ on $\mathcal{L}$. Using the fact $\pi^{\textrm{\'et}}_{1}(X,\overline{x}) = \textrm{Aut}(F_{X, \overline{x}})$, we may define the map $\pi^{\textrm{\'et}}_{1}(X,\overline{x}) \to \Aut_{X}(Y)$ by sending $\eta \in \textrm{Aut}(F_{X, \overline{x}})$ to the automorphism $\alpha \in \Aut_{X}(Y)$ for which $\eta(\overline{y}) = \alpha(\overline{y})$ for all $\overline{y} \in F_{X,\overline{x}}(Y)$; that such an element exists follows from the assumption that $Y \to X$ be a Galois cover.

In the situation where $\Lambda = \varprojlim \Lambda_{n}$ is a pro-finite group, one can take the limit of both sides of (\ref{HomLambda}) to obtain a bijection 
\begin{equation}
\label{HomLambda2}
\{ \textrm{isom. classes of }\Lambda\textrm{-torsors} \} \longleftrightarrow \Hom_{\textrm{cont}}(\pi^{\textrm{\'et}}_{1}(X,\overline{x}), \Lambda) ,
\end{equation}
with a similar explicit description in the case of a torsor coming from a limit of Galois coverings.

\section{Cohomological Computations}
\label{cohomcompsec}

\subsection{Basic \v{C}ech Computations}
\label{Cechcalcsec}

Let $k$ be a discretely-valued characteristic zero non-archimedian local field, with ring of integers $\mathcal{O}_{k}$, residue field $\kappa$, and completed closure $\widehat{k}$. Write $\Delta^{\circ}$ for the adic space $\Spa(k\langle T^{\pm 1} \rangle, \mathcal{O}_{k} \langle T^{\pm 1} \rangle)$ over $\Spa(k, \mathcal{O}_{k})$, which can be thought of as a rigid-analytic annulus. We consider the natural cover of $\Delta^{\circ}_{\widehat{k}}$ inside $\Delta^{\circ}_{\widehat{k}, \textrm{pro\'et}}$ with the covering space modelled by the infinite tower
\[ \widetilde{\Delta}^{\circ} = \varprojlim_{T \mapsto T^{p}} \Delta^{\circ}_{\widehat{k}} . \]
The space $\widetilde{\Delta}^{\circ}$ is then perfectoid of the form 
\[\Spa(\widehat{k}\langle T^{\pm 1/p^{\infty}} \rangle, \mathcal{O}_{\widehat{k}} \langle T^{\pm 1/p^{\infty}} \rangle) = \varprojlim \Spa(\widehat{k}\langle T^{\pm 1/p^{j}} \rangle, \mathcal{O}_{\widehat{k}} \langle T^{\pm 1/p^{j}} \rangle) , \]
and the covering map $c : \widetilde{\Delta}^{\circ} \to \Delta^{\circ}_{\widehat{k}}$ with respect to these presentations is simply given by $T \mapsto T$.

Define $\mathbb{Z}_{p}(1) = \varprojlim_{j} \mu_{p^{j}}$, where the transition maps are given by $x \mapsto x^{p}$. We consider the self-product $\widetilde{\Delta}^{\circ 2} = \widetilde{\Delta}^{\circ} \times_{\Delta^{\circ}} \widetilde{\Delta}^{\circ}$, and observe that its connected components are naturally indexed by $\mathbb{Z}_{p}(1)$. Indeed, one has that
\[ \widetilde{\Delta}^{\circ 2} = \lim_{T \mapsto T^{p^{j}}} \underbrace{\Delta^{\circ}_{\widehat{k}} \times_{\Delta^{\circ}_{\widehat{k}}} \Delta^{\circ}_{\widehat{k}}}_{\Delta^{\circ}_{j}} . \]
Where the fibre product $\Delta^{\circ}_{j}$ may be modelled as
\[ \Delta^{\circ}_{j} = \Spa(\widehat{k}\langle T_{1}^{\pm 1}, T_{2}^{\pm 1} \rangle / (T_{1}^{p^{j}} - T_{2}^{p^{j}}), \mathcal{O}_{\widehat{k}}\langle T_{1}^{\pm 1}, T_{2}^{\pm 1} \rangle / (T_{1}^{p^{j}} - T_{2}^{p^{j}})) , \]
and the transition maps are given by $(T_{1}, T_{2}) \mapsto (T_{1}^{p}, T_{2}^{p})$. It is clear that the components of $\Delta^{\circ}_{j}$ are naturally identified with the group $\mu_{p^{j}}$ of $p^{j}$'th roots of unity, with $\zeta_{j} \in \mu_{p^{j}}$ identified with the component on which $T_{1} - \zeta_{j} T_{2} = 0$, and hence the idempotents of the coordinate ring of $\widetilde{\Delta}^{\circ 2}$ are identified with a compatible system of such roots and hence with $\mathbb{Z}_{p}(1)$.

Let $\mathcal{F}^{\bullet} = \mathbb{\Omega}^{\bullet}_{\Delta^{\circ}}$ be the de Rham complex defined in \S\ref{proetsec}. We let $\mathcal{U} = \{ c \}$ be our cover, and form the \v{C}ech complex $C^{\bullet}(\mathcal{U}, \mathcal{F}^{\bullet})$ and the associated total complex $L^{\bullet} = L^{\bullet}(\mathcal{F}^{\bullet})$ with differential $d + \delta$. We will consider the cocycle $s_{c} = \log (T/[T^{\flat}]) \in C^{0}(\mathcal{U}, \mathcal{O} \mathbb{B}_{\textrm{dR}})$. To be more precise, \cite[Cor. 6.13]{zbMATH06209107} shows that the sequence
\begin{equation}
\label{BdRres}
0 \to \mathbb{B}_{\textrm{dR}} \to \mathcal{O} \mathbb{B}_{\textrm{dR}} \xrightarrow{d} \mathbb{\Omega}^{1}_{\Delta^{\circ}} \to 0
\end{equation}
is exact, and does so by showing that a natural map $\restr{\mathcal{O}\mathbb{B}^{+}_{\textrm{dR}}}{\widetilde{\Delta}^{\circ}} \xrightarrow{\sim} \restr{\mathbb{B}^{+}_{\textrm{dR}}}{\widetilde{\Delta}^{\circ}}[[X]]$ sending $T \otimes 1$ to $[T^{\flat}] + X$ is an isomorphism; here $T^{\flat}$ is defined as in \cite[\S6]{zbMATH06209107}. One does this by showing that $\restr{\mathbb{B}^{+}_{\textrm{dR}}}{\widetilde{\Delta}^{\circ}}[[X]]$ admits the structure of an $\restr{\mathcal{O}_{\Delta^{\circ}}}{\widetilde{\Delta}^{\circ}}$-algebra, compatible with the one on the quotient 
\[ \mathbb{B}^{+}_{\textrm{dR}}[[X]]/(\ker \theta) = \hat{\mathcal{O}}_{\Delta^{\circ}} . \]
This then gives a natural map 
\begin{equation}
\label{sheafBdRtriviso}
\restr{\left(\mathcal{O}_{\Delta^{\circ}} \otimes_{W(\kappa)} W(\hat{\mathcal{O}}^{+}_{\Delta^{\circ \flat}})\right)}{\widetilde{\Delta}^{\circ}} \xrightarrow{\sim} \restr{\mathbb{B}^{+}_{\textrm{dR}}}{\widetilde{\Delta}^{\circ}}[[X]]
\end{equation}
inducing the inverse of the map $X \mapsto T \otimes 1 - 1 \otimes [T^{\flat}]$.

Using this description, one can define a section $\log (T/[T^{\flat}])$ of $\restr{\mathcal{O}\mathbb{B}^{+}_{\textrm{dR}}}{\widetilde{\Delta}^{\circ}}$ by applying the isomorphism (\ref{sheafBdRtriviso}), computing $\log(1 + X/[T^{\flat}])$ (using the power series expansion) and then using the isomorphism to translate the resulting expression back. The resulting function satisfies the property that $d(\log (T/[T^{\flat}])) = dT/T$, and that $\log(a (T/[T^{\flat}])) = \log(a) + \log(T/[T^{\flat}])$ for any non-zero $a \in B_{\textrm{dR}}$ for which $\log(a)$ is defined. 

We note that we have two natural maps $p_{i} : \widetilde{\Delta}^{\circ 2} \to \widetilde{\Delta}^{\circ}$ with $i \in \{ 1, 2 \}$. If one models $\widetilde{\Delta}^{\circ 2}$ as the space
\[ \widetilde{\Delta}^{\circ 2} = \Spa(\widehat{k}\langle T_{1}^{\pm p^{\infty}}, T_{2}^{\pm p^{\infty}} \rangle / (T_{1} - T_{2}), \mathcal{O}_{\widehat{k}}\langle T_{1}^{\pm p^{\infty}}, T_{2}^{\pm p^{\infty}} \rangle / (T_{1} - T_{2})) , \]
then these maps are given by $T \mapsto T_{i}$. The component $\widetilde{\Delta}^{\circ 2}_{\zeta_{\bullet}} \subset \widetilde{\Delta}^{\circ 2}$ corresponding to the sequence $\zeta_{\bullet} = (\zeta_{1}, \zeta_{2}, \hdots)$ is then given by imposing the infinitely many relations 
\[ T_{1}^{1/p} = \zeta_{1} T_{2}^{1/p}, \hspace{1em} T_{1}^{1/p^2} = \zeta_{2} T_{2}^{1/p^2},  \hspace{1em} \hdots \]
If we consider the restrictions (isomorphisms) $p_{i,\zeta_{\bullet}} : \widetilde{\Delta}^{\circ 2}_{\zeta_{\bullet}} \xrightarrow{\sim} \widetilde{\Delta}^{\circ}$ associated to this component, they are, on the level of the ring maps $r_{i,\zeta_{\bullet}}$, related by the fact that $r_{1,\zeta_{\bullet}}(T^{1/p^k}) = \zeta_{k} r_{2,\zeta_{\bullet}}(T^{1/p^k})$ for all $k \geq 1$. 

We now apply the differential $d + \delta$ to the cocycle $s_{c}$, which has its image inside $L^{1}(\mathcal{U}, \mathbb{\Omega}^{\bullet}_{\Delta^{\circ}}) = \mathbb{\Omega}^{1}_{\Delta^{\circ}}(\widetilde{\Delta}^{\circ}) \oplus \mathcal{O}\mathbb{B}_{\textrm{dR}}(\widetilde{\Delta}^{\circ 2})$. The result is the direct sum of $dT/T$, regarded as a differential form on $\widetilde{\Delta}^{\circ}$, and the difference $\log(T_{1}/[T_{1}^{\flat}]) - \log(T_{2}/[T_{2}^{\flat}])$. If we compute this latter difference on the component $\widetilde{\Delta}^{\circ 2}_{\zeta_{\bullet}}$ we obtain $-\log([\zeta_{\bullet}])$, coming from the fact that $[T_{1}^{\flat}] = [\zeta_{\bullet}][T_{2}^{\flat}]$. We thus have that $(d + \delta)(s_{c}) = dT/T - \mathbf{t}$, where $\mathbf{t}$ is the section of $\restr{\mathbb{B}_{\textrm{dR}}}{\widetilde{\Delta}^{\circ 2}}$ whose value on each component $\widetilde{\Delta}^{\circ 2}_{\zeta_{\bullet}}$ is $\log([\zeta_{\bullet}])$. Our conclusion is that

\begin{prop}
\label{tcalccor}
In the \v{C}ech complex associated to the sheaf $\mathbb{\Omega}^{\bullet}_{\Delta^{\circ}}$ and the cover $\mathcal{U} = \{ c \}$ the cycles $dT/T$ and $\mathbf{t}$ are cohomologous. \qed
\end{prop}

\begin{cor}
\label{dtnonzero}
The class of $dT/T$ is non-zero in the cohomology of the complex $\mathbb{\Omega}^{\bullet}_{\Delta^{\circ}}$.
\end{cor}

\begin{proof}
Using (\ref{BdRres}), we have a map $H^{1}(\Delta^{\circ}_{\textrm{p\'et}}, \mathbb{B}_{\textrm{dR}}) \xrightarrow{\sim} H^{1}(\Delta^{\circ}_{\textrm{p\'et}}, \mathbb{\Omega}^{\bullet}_{\Delta^{\circ}})$. In \v{C}ech cohomology, this map corresponds to the map $L^{\bullet}(\mathcal{U}, \mathbb{B}_{\textrm{dR}}) \to L^{\bullet}(\mathcal{U}, \mathbb{\Omega}^{\bullet}_{\Delta^{\circ}})$ of total \v{C}ech complexes (recall \S\ref{cechcompsec}), and $\mathbf{t}$ is also naturally a cocycle of $L^{\bullet}(\mathcal{U}, \mathbb{B}_{\textrm{dR}})$. By \autoref{tcalccor}, it therefore suffices to check the class of $\mathbf{t}$ is non-zero regarded as an element of the cohomology of $L^{\bullet}(\mathcal{U}, \mathbb{B}_{\textrm{dR}})$, which computes the cohomology of $H^{1}(\Delta^{\circ}_{\textrm{p\'et}}, \mathbb{B}_{\textrm{dR}})$ as a consequence of \cite[Thm. 6.5(ii)]{zbMATH06209107}, the fact that $c$ is an affinoid perfectoid cover, and \autoref{Cechcomputescohomcor}.

Thus we are asking whether there is an element $s$ of $\restr{\mathbb{B}_{\textrm{dR}}}{\widetilde{\Delta}^{\circ}}$ such that $\delta(s) = \restr{s}{p_{1,\zeta_{\bullet}}} - \restr{s}{p_{2,\zeta_{\bullet}}}$ agrees with a constant element $\log([\zeta_{\bullet}])$ on the component $\widetilde{\Delta}^{\circ 2}_{\zeta_{\bullet}}$ for each choice of $\zeta_{\bullet}$. Suppose there exists such an $s$, and choose $i$ such that $s$ is a section of $\textrm{Fil}^{i} \mathbb{B}_{\textrm{dR}}$. Necessarily $i \leq 1$, since $\mathbf{t}$ is a section of $\textrm{Fil}^{1} \mathbb{B}_{\textrm{dR}}$ and not of $\textrm{Fil}^{2} \mathbb{B}_{\textrm{dR}}$, and the restrictions preserve the filtration. We may then replace $s$ with an element of the associated graded quotient $\textrm{Fil}^{i} \mathbb{B}_{\textrm{dR}} / \textrm{Fil}^{2} \mathbb{B}_{\textrm{dR}} \cong \bigoplus_{i \leq j \leq 2} \xi^{j} \restr{\widehat{\mathcal{O}}_{\Delta^{\circ}}}{\widetilde{\Delta}^{\circ}}$ (see \cite[Cor. 6.4]{zbMATH06209107}), where $\xi$ is any generator of the kernel of Fontaine's map $\theta$ (recall \S\ref{proetsec}), and $\mathbf{t}$ with its image $\mathbf{i}$ in the same quotient. Setting $R = \restr{\widehat{\mathcal{O}}_{\Delta^{\circ}}}{\widetilde{\Delta}^{\circ}}$, $\mathbf{i}$ is identified with a collection of constant elements of $\xi^{j} R$ indexed by $\zeta_{\bullet}$ and $j$. By multiplying by $\xi^{-j}$, we can reduce to the same question where $\mathbf{i}$ is a collection of constant elements of $R$, and $s$ is likewise an element of $R$ (because the restrictions preserve the graded direct sum decomposition).

The two restrictions of $s$ are related, on the component corresponding to $\zeta_{\bullet}$, by an automorphism of $\widetilde{\Delta}^{\circ 2}_{\zeta_{\bullet}}$ inducing an automorphism of $\widehat{\mathcal{O}}_{\widetilde{\Delta}^{\circ 2}_{\zeta_{\bullet}}}$. But one easily checks that this automorphism, which is induced by scaling $p$'th roots of $T$ by the corresponding elements of $\zeta_{\bullet}$, cannot shift any function by a non-zero constant. 
\end{proof}

\subsection{Evaluation Functionals}
\label{alphaconstrsec}

We now consider the more general setting where we have a product $\Delta^{a,b} = (\Delta^{\circ})^{a} \times \Delta^{b}$, where $\Delta$ denotes the closed $v$-adic disk of radius $1$. Explicitly, this is given by
\[ \Delta^{a,b} = \Spa(k\langle T_{1}^{\pm 1}, \hdots, T_{a}^{\pm 1}, T_{a+1}, \hdots, T_{a+b} \rangle, \mathcal{O}_{k}\langle T_{1}^{\pm 1}, \hdots, T_{a}^{\pm 1}, T_{a+1}, \hdots, T_{a+b} \rangle) . \]

In this section we will view all spaces, including algebraic varieties, as adic spaces; in particular we consider the multiplicative group $\mathbb{G}_{m}$ and the affine line $\mathbb{A}^{1}$ as adic spaces over $\Spa(k, \mathcal{O}_{k})$. Pro-\'etale cohomology will be computed over the algebraic closure $\overline{k}$ of $k$. We wish to define certain ``evaluation functionals'' $\hat{\alpha}^{*}_{a,b} : H^{a}(\Delta^{a,b}_{\textrm{p\'et}}, \hat{\mathbb{Z}}_{p}(a)) \otimes \mathbb{Q}_{p} \to \mathbb{Q}_{p}(a)$ and study their relationship with the $p$-adic Hodge comparisons and our calculation in the previous section. 

\begin{defn}
By a Kummer-type cover of $V^{a,b} := \mathbb{G}^{a}_{m} \times \mathbb{A}^{b}$ we mean a finite-\'etale covering of $V^{a,b}$. By a Kummer-type cover of $\Delta^{a,b}$ we mean a cover obtained by pullback of a Kummer-type cover of $V^{a,b}$ under the natural map $\Delta^{a,b} \to V^{a,b}$.
\end{defn}

\begin{lem}
\label{Kpi1lem}
For each $n$ and $r$, the map $H^{n}(V^{a,b}_{\textrm{pf\'et}}, \widehat{\mathbb{Z}}_{p}(r)) \to H^{n}(V^{a,b}_{\textrm{p\'et}}, \widehat{\mathbb{Z}}_{p}(r))$ is an isomorphism: the \'etale cohomology can be computed with only finite \'etale coverings.
\end{lem}

\begin{proof}
It suffices to prove the analogous statement for the sites $V^{a,b}_{\textrm{f\'et}}$ and $V^{a,b}_{\textrm{\'et}}$, i.e., before passing to pro-sites. By \cite[Prop. 2.1.5(e)]{kpi1}, it suffices to show that $V^{a,b}$ is $K(\pi,1)$ (satisfies one of the equivalent conditions of \cite[Prop. 2.1.5]{kpi1}). Using \cite[Prop. 2.1.8(b)]{kpi1} we may reduce to the same statement with $V^{a,b}$ regarded as a variety over $\overline{\mathbb{Q}} \subset \overline{k}$, and then by base-change along $\overline{\mathbb{Q}} \hookrightarrow \mathbb{C}$, the same statement with $V^{a,b}$ regarded as a complex algebraic variety. Now the analytification of $(V^{a,b})^{\textrm{an}}$ has contractable universal cover, hence is topologically $K(\pi,1)$ by \cite[Prop. 2.1.1(b)]{kpi1}. The result then follows by $(2),(3) \implies (1)$ in \cite[Prop. 2.1.15]{kpi1}.
\end{proof}

\noindent We have a natural map $H^{a}(V^{a,b}_{\textrm{pf\'et}}, \hat{\mathbb{Z}}_{p}(a)) \to H^{a}(\Delta^{a,b}_{\textrm{p\'et}}, \hat{\mathbb{Z}}_{p}(a))$ induced by the map $\Delta^{a,b}_{\textrm{p\'et}} \to V^{a,b}_{\textrm{pf\'et}}$ of sites, and we will denote by $I^{a,b}$ its image. By \autoref{Kpi1lem} this agrees with the image of $H^{a}(V^{a,b}_{\textrm{p\'et}}, \hat{\mathbb{Z}}_{p}(a))$.

We will define $\hat{\alpha}^{*}_{a,b}$ as follows. We will first define a functional $\alpha^{*}_{a,b}$ on $H^{a}(V^{a,b}_{\textrm{pf\'et}}, \hat{\mathbb{Z}}_{p}(a))$, which induces a functional on $I^{a,b}$ and also $I^{a,b}_{\mathbb{Q}_{p}}$. We will then choose a splitting $H^{a}(\Delta^{a,b}_{\textrm{p\'et}}, \hat{\mathbb{Z}}_{p}(a)) \otimes \mathbb{Q}_{p} = I^{a,b}_{\mathbb{Q}_{p}} \oplus J^{a,b}$, and define $\hat{\alpha}^{*}_{a,b}$ by extending by zero. We will start with the case of $I^{1,0}$, for which we need:

\begin{lem}
The map $H^{1}(\mathbb{G}_{m,\textrm{p\'et}}, \hat{\mathbb{Z}}_{p}(1)) \to H^{1}(\Delta^{\circ}_{\textrm{p\'et}}, \hat{\mathbb{Z}}_{p}(1))$ is injective.
\end{lem}

\begin{proof}
From our discussion in \S\ref{fundgpsec} we may identify this map with the natural map $\Hom_{\textrm{cont}}(\pi^{1}_{\textrm{\'et}}(\mathbb{G}_{m}), \mathbb{Z}_{p}(1)) \to \Hom_{\textrm{cont}}(\pi^{1}_{\textrm{\'et}}(\Delta^{\circ}), \mathbb{Z}_{p}(1))$, so we are reduced to showing that $\pi^{1}_{\textrm{\'et}}(\Delta^{\circ}) \to \pi^{1}_{\textrm{\'et}}(\mathbb{G}_{m})$ is surjective. By \cite[\href{https://stacks.math.columbia.edu/tag/0BN6}{Lemma 0BN6}]{stacks-project}, we reduce to showing that every connected finite \'etale cover of $\mathbb{G}_{m}$ pulls back to a connected finite \'etale cover of $\Delta^{\circ}$, which is obvious. (Note that the rigid analytic finite \'etale coverings of $\mathbb{G}_{m}$ are just those of Kummer type as a consequence of the rigid analytic Riemann existence theorem \cite{zbMATH00224075}.)
\end{proof}

We now fix a distinguished system $\{ \zeta^{\star}_{i} \}_{i \geq 1}$ of $p$-power roots of unity of $\widehat{k}^{\flat +}$. This induces the following data:
\begin{itemize}
\item[-] A $p$-adic period $t = \log([\zeta^{\star}_{\bullet}]) \in B_{\textrm{dR}}$.
\item[-] Via the automorphisms $T \mapsto \zeta^{\star}_{i} T$, an element, denoted $\alpha$, of the pro-$p$ fundamental group $\pi^{\textrm{\'et}}_{1}(\mathbb{G}_{m})^{(p)}$; note that by the rigid analytic Riemann existence theorem \cite{zbMATH00224075} the covers $\mathbb{G}_{m} \to \mathbb{G}_{m}$ given by $T \mapsto T^{p^{i}}$ exhaust the connected finite \'etale coverings of $\mathbb{G}_{m}$ with degree a power of $p$, even on the adic finite \'etale site.
\item[-] A map $\alpha^{*} : H^{1}(\mathbb{G}_{m,\textrm{p\'et}}, \widehat{\mathbb{Z}}_{p}(1)) \to \mathbb{Z}_{p}(1)$. This uses the identification 
\[ H^{1}(\mathbb{G}_{m,\textrm{p\'et}}, \widehat{\mathbb{Z}}_{p}(1)) \simeq \Hom_{\textrm{cont}}(\pi^{\textrm{\'et}}_{1}(\mathbb{G}_{m})^{(p)}, \widehat{\mathbb{Z}}_{p}(1)) \]
and is defined by evaluation on $\alpha$.
\item[-] Using the identification $H^{1}(\mathbb{G}_{m,\textrm{p\'et}}, \widehat{\mathbb{Z}}_{p}(1)) \simeq H^{1}(\mathbb{G}_{m,\textrm{pf\'et}}, \widehat{\mathbb{Z}}_{p}(1))$, a map, also denoted $\alpha^{*}$, on the latter cohomology group.
\item[-] A map $\hat{\alpha}^{*} : I^{1,0} \to \mathbb{Z}_{p}(1)$, obtained by pulling back $\alpha^{*}$ along the map $I^{1,0} \to H^{1}(\mathbb{G}_{m,\textrm{pf\'et}}, \widehat{\mathbb{Z}}_{p}(1))$ (recall that $H^{1}(\mathbb{G}_{m,\textrm{pf\'et}}, \widehat{\mathbb{Z}}_{p}(1)) = H^{1}(\mathbb{G}_{m,\textrm{p\'et}}, \widehat{\mathbb{Z}}_{p}(1))$ by \autoref{Kpi1lem}), and its scalar extension $\hat{\alpha}^{*}_{\mathbb{Q}_{p}} : I^{1,0}_{\mathbb{Q}_{p}} \to \mathbb{Q}_{p}(1)$.  
\item[-] A linear functional
\[ \hat{\alpha}^{*}_{B_{\textrm{dR}}} : I^{1,0} \otimes B_{\textrm{dR}} \to B_{\textrm{dR}} , \]
which is defined by evaluating on $\alpha$ and extending scalars along the map $\mathbb{Z}_{p}(1) \hookrightarrow B_{\textrm{dR}}$ given by $\log([-])$.
\end{itemize}

From the coverings $T \mapsto T^{p^{j}}$ of $\Delta^{\circ}$ in the previous section we obtain torsors $\mathcal{L}_{j}$ on $\Delta^{\circ}_{\textrm{p\'et}}$ and a class $[\mathcal{L}_{\infty}] = (\lim_{j} [\mathcal{L}_{j}]) \in I^{1,0}$. To complete our definitions we will need the following facts.

\begin{lem}
\label{formalevallem}
Fix $a, b \geq 0$, and let $\hat{\alpha}^{*}_{i} : H^{1}(V^{a,b}_{\textrm{pf\'et}}, \widehat{\mathbb{Z}}_{p}(1)) \to \mathbb{Z}_{p}(1)$ be the pullback of $\hat{\alpha}^{*}$ along the $i$'th projection $H^{1}(V^{a,b}_{\textrm{pf\'et}}, \widehat{\mathbb{Z}}_{p}(1)) \to H^{1}(\mathbb{G}_{m,\textrm{pf\'et}}, \widehat{\mathbb{Z}}_{p}(1))$ induced by inclusion of factors. Then the map 
\[ \hat{\alpha}^{*}_{a,b} := \hat{\alpha}^{*}_{1} \otimes \cdots \otimes \hat{\alpha}^{*}_{a} : \underbrace{H^{1}(V^{a,b}_{\textrm{pf\'et}}, \widehat{\mathbb{Z}}_{p}(1)) \otimes \cdots \otimes H^{1}(V^{a,b}_{\textrm{pf\'et}}, \widehat{\mathbb{Z}}_{p}(1))}_{= H^{a}(V^{a,b}_{\textrm{pf\'et}}, \widehat{\mathbb{Z}}_{p}(a))} \to \mathbb{Z}_{p}(a) \]
takes the value $\{ \zeta^{\star}_{\bullet} \} \otimes \cdots \otimes \{ \zeta^{\star}_{\bullet} \}$ on the element $[\mathcal{L}_{\infty, 1}] \otimes \cdots \otimes [\mathcal{L}_{\infty, a}]$, where $[\mathcal{L}_{\infty, i}]$ is induced from $[\mathcal{L}_{\infty}]$ by the inclusion $H^{1}(\mathbb{G}_{m,\textrm{pf\'et}}, \widehat{\mathbb{Z}}_{p}(1)) \to H^{1}(V^{a,b}_{\textrm{pf\'et}}, \widehat{\mathbb{Z}}_{p}(1))$ coming from the projection onto the $i$'th factor.
\end{lem}

\begin{proof}
Immediate from the definitions and functoriality of cohomology.
\end{proof}

\begin{lem}
The map $H^{a}(V^{a,b}_{\textrm{pf\'et}}, \hat{\mathbb{Z}}_{p}(a)) \to H^{a}(\Delta^{a,b}_{\textrm{p\'et}}, \hat{\mathbb{Z}}_{p}(a))$ is injective.
\end{lem}

\begin{proof}
From the Kunneth formula, the cohomology group $H^{a}(V^{a,b}_{\textrm{pf\'et}}, \hat{\mathbb{Z}}_{p}(a))$ is generated by the class $[\mathcal{L}_{\infty, 1}] \cup \cdots \cup [\mathcal{L}_{\infty, a}]$, so it suffices to show the image of this class is non-zero. Using the functoriality of cup product and its compatibility with the differential graded structure of $\mathbb{\Omega}^{\bullet}_{\Delta^{a,b}}$ (as discussed in \S\ref{cechcupprodsec}), this will follow from \autoref{epsendslem} below.
\end{proof}

\begin{lem}
\label{epsendslem}
Denote by $\ep$ the natural comparison map
\[ \ep : H^{a}(\Delta^{a,b}_{\textrm{p\'et}}, \widehat{\mathbb{Z}}_{p}(a)) \otimes B_{\textrm{dR}} \to H^{a}(\Delta^{a,b}_{\textrm{p\'et}}, \mathbb{\Omega}^{\bullet}_{\Delta^{a,b}}) . \]
Then $\ep$ maps the element $\left( [\mathcal{L}_{1,\infty}] \cup \cdots \cup [\mathcal{L}_{a,\infty}] \right)$ to the element $dT_{1} / T_{1} \wedge \cdots \wedge dT_{a} / T_{a}$, and this element is non-zero in $H^{a}(\Delta^{a,b}_{\textrm{p\'et}}, \mathbb{\Omega}^{\bullet}_{\Delta^{a,b}})$. 
\end{lem}

\begin{proof}
If we compute cup product using \v{C}ech cohomology on both sides, and also on the cohomology of the sheaf $\mathbb{B}_{\textrm{dR}}$, the first part of the result (i.e., ignoring the non-zeroness of $dT_{1} / T_{1} \wedge \cdots \wedge dT_{a} / T_{a}$) will follow from the compatibility of cup product with the differential graded structure on $\mathbb{\Omega}^{\bullet}_{\Delta^{a,b}}$ (as discussed in \S\ref{cechcupprodsec}), as well as our result \autoref{formalevallem} above, as long as we can show that $[\mathcal{L}_{1,\infty}]$ maps to $\mathbf{t}_{1}$, where $\mathbf{t}_{1}$ is the class in $H^{1}(\Delta^{a,b}_{\textrm{p\'et}}, \mathbb{B}_{\textrm{dR}})$ obtained as the image of $\mathbf{t}$ under the map $H^{1}(\Delta^{\circ}_{\textrm{p\'et}}, \mathbb{B}_{\textrm{dR}}) \to H^{1}(\Delta^{a,b}_{\textrm{p\'et}}, \mathbb{B}_{\textrm{dR}})$ coming from projection. 

From functoriality it suffices to show that $[\mathcal{L}_{\infty}]$ maps to $\mathbf{t}$ under the natural map $H^{1}(\Delta^{\circ}_{\textrm{p\'et}}, \hat{\mathbb{Z}}_{p}(1)) \to H^{1}(\Delta^{\circ}_{\textrm{p\'et}}, \mathbb{B}_{\textrm{dR}})$. Working on the level of \v{C}ech complexes with respect to the perfectoid cover $c$ as in \S\ref{Cechcalcsec}, the gluing data for the torsor $[\mathcal{L}_{\infty}]$ assigns the system of compatible roots $\zeta_{\bullet}$ to the component $\widetilde{\Delta}^{\circ 2}_{\zeta_{\bullet}}$. As the map $\widehat{\mathbb{Z}}_{p}(1) \to \mathbb{B}_{\textrm{dR}}$ is induced by $(\zeta_{i})_{i \geq 1} \mapsto \log([\zeta_{\bullet}])$, one obtains the cycle $\mathbf{t}$ as desired.

For the non-zeroness of $dT_{1} / T_{1} \wedge \cdots \wedge dT_{a} / T_{a}$ we may argue as follows. We may first reduce to the case of $b = 0$ by using the natural map $H^{a}(\Delta^{a,b}_{\textrm{p\'et}}, \mathbb{B}_{\textrm{dR}}) \to H^{a}(\Delta^{a,0}_{\textrm{p\'et}}, \mathbb{B}_{\textrm{dR}})$ coming from inclusion. Because \autoref{tcalccor} and the cup-product compatibility implies that the class of $dT_{1} / T_{1} \wedge \cdots \wedge dT_{a} / T_{a}$ is cohomologous to the class of $\mathbf{t}_{1} \cup \cdots \cup \mathbf{t}_{a}$, we may reduce to the same statement for the latter. Using a Kunneth calculation to compute the cohomology of $\Delta^{a,0}$ this then reduces to showing that each $\mathbf{t}_{i}$ is non-zero, which is what we showed in the argument of \autoref{dtnonzero}. 
\end{proof}

We may now choose a splitting $H^{a}(\Delta^{a,b}_{\textrm{p\'et}}, \hat{\mathbb{Z}}_{p}(a)) \otimes \mathbb{Q}_{p} = I^{a,b}_{\mathbb{Q}_{p}} \oplus J^{a,b}$ such that $J^{a,b}$ contains the kernel of $H^{a}(\Delta^{a,b}_{\textrm{p\'et}}, \hat{\mathbb{Z}}_{p}(a)) \otimes \mathbb{Q}_{p} \to H^{a}(\Delta^{a,b}_{\mathbb{C}_{p}, \textrm{p\'et}}, \mathbb{B}_{\textrm{dR}})$. We then define $\hat{\alpha}^{*}_{a,b}$ on all of $H^{a}(\Delta^{a,b}_{\textrm{p\'et}}, \hat{\mathbb{Z}}_{p}(a)) \otimes \mathbb{Q}_{p}$ by extending by zero. 

\subsection{Extending to an ambient variety}
\label{funcextsec}

Fix a finite extension $k = K_{v}$ of $\mathbb{Q}_{p}$. We will now suppose that $Y$ is a smooth proper adic space over $\Spa(K_{v}, \mathcal{O}_{K_{v}})$, and compare two linear functionals associated to an embedding $\Delta^{a,b} \hookrightarrow Y$ of adic spaces defined over $K_{v}$. On the de Rham side we suppose we have a functional $\hat{\gamma}^{*}_{\textrm{dR}} : H^{a}_{\textrm{dR}}(Y) \to K_{v}$ which is defined on $\omega \in H^{a}_{\textrm{dR}}(Y)$ by $\restr{\omega}{\Delta^{a,b}} = \hat{\gamma}^{*}_{\textrm{dR}}(\omega) \, dT_{1}/T_{1} \wedge \cdots \wedge dT_{a} / T_{a}$ in cohomology. We compare such a functional to the pullback of $\hat{\alpha}^{*}_{a,b}$ to the cohomology of $Y_{\overline{K_{v}}}$. 

From the natural morphisms $\hat{\mathbb{Z}}_{p}(a) \to \mathbb{B}_{\textrm{dR}}$, $\mathbb{B}_{\textrm{dR}} \to \mathbb{\Omega}^{\bullet}_{(-)}$ and $\Omega^{\bullet}_{(-)} \to \mathbb{\Omega}^{\bullet}_{(-)}$ of sheaves on the pro-\'etale site, one obtains the following diagram, where all cohomology is computed on the pro-\'etale site:

\vspace{0.8em}

\begin{equation*}
\begin{tikzcd}[center picture]
H^{a}(Y_{\overline{K_{v}}}, \hat{\mathbb{Z}}_{p}(a)) \otimes B_{\textrm{dR}} \arrow[r, "\sim"] \arrow[d] & H^{a}(Y_{\mathbb{C}_{p}}, \mathbb{B}_{\textrm{dR}}) \arrow[r, "\sim"] \arrow[d] & H^{a}(Y_{\mathbb{C}_{p}}, \mathbb{\Omega}^{\bullet}_{Y}) \arrow[d] & \arrow[l, swap, "\sim"] H^{a}(Y, \Omega^{\bullet}_{Y}) \otimes B_{\textrm{dR}} \arrow[d] \\
H^{a}(\Delta^{a,b}_{\overline{K_{v}}}, \hat{\mathbb{Z}}_{p}(a)) \otimes B_{\textrm{dR}} \arrow[r,"\sigma"] & H^{a}(\Delta^{a,b}_{\mathbb{C}_{p}}, \mathbb{B}_{\textrm{dR}}) \arrow[r, "\sim"] &  H^{a}(\Delta^{a,b}_{\mathbb{C}_{p}}, \mathbb{\Omega}^{\bullet}_{\Delta^{a,b}}) & \arrow[l, swap, "\lambda"] H^{a}(\Delta^{a,b}, \Omega^{\bullet}_{\Delta^{a,b}}) \otimes B_{\textrm{dR}}
\end{tikzcd}
\end{equation*}

\vspace{0.8em}

All the squares in the diagram are commutative by general cohomological principles. That the middle horizontal rightward arrows are isomorphisms is \cite[6.13]{zbMATH06209107}. That the upper left horizontal arrow is an isomorphism is \cite[8.4]{zbMATH06209107}. That the upper right horizontal arrow is an isomorphism follows from the proof of \cite[7.11]{zbMATH06209107}. 

For later use we define a map $\hat{\gamma}^{*}_{\textrm{\'et}} : H^{a}_{\textrm{\'et}}(Y_{\overline{K_{v}}},\widehat{\mathbb{Z}}_{p}(a)) \otimes \mathbb{Q}_{p} \to \mathbb{Q}_{p}(a)$ obtained by composing the $\mathbb{Q}_{p}$-scalar extension of $H^{a}_{\textrm{\'et}}(Y_{\overline{K_{v}}},\widehat{\mathbb{Z}}_{p}(a)) \to H^{a}_{\textrm{\'et}}(\Delta^{a,b}_{\overline{K_{v}}}, \mathbb{Z}_{p}(a))$ with $\hat{\alpha}^{*}_{a,b}$. The next lemma shows that, when pulled back to $H^{a}(Y, \Omega^{\bullet}_{Y}) \otimes B_{\textrm{dR}}$ along the $p$-adic Hodge comparison, $\hat{\gamma}^{*}_{{\textrm{\'et}}, B_{\textrm{dR}}}$ agrees with $\hat{\gamma}^{*}_{\textrm{dR}, B_{\textrm{dR}}}$. 

\begin{lem}
\label{indepofpathlem}
Let $\hat{\gamma}^{*}_{a,b}$ be the pullback of $\hat{\alpha}^{*}_{a,b} : H^{a}(\Delta^{a,b}_{\textrm{p\'et}}, \widehat{\mathbb{Z}}_{p}(a)) \to \mathbb{Z}_{p}(a)$ to $H^{a}(Y_{\textrm{p\'et}}, \Omega^{\bullet}_{Y}) \otimes B_{\textrm{dR}}$ obtained using the leftmost arrow and the morphisms on the top row. Then if $\omega \in H^{a}(Y_{\textrm{p\'et}}, \Omega^{\bullet}_{Y})$, we have $\frac{1}{t^{a}} \hat{\gamma}^{*}_{a,b}(\omega) = \hat{\gamma}^{*}_{\textrm{dR}}(\omega) \otimes 1$.
\end{lem}

\begin{proof}
Invert all the arrows labelled isomorphisms in the above diagram. Then the commutativity of the diagram implies that the class $\omega$ has the same image in $H^{a}(\Delta^{a,b}_{\mathbb{C}_{p}}, \mathbb{B}_{\textrm{dR}})$ regardless of which path one takes. Travelling first along the right-most arrow, then $\lambda$, and then the bottom middle arrow, we obtain, using \autoref{tcalccor} and the compatibility of cup product with the graded differential structure discussed in \S\ref{cechcupprodsec}, the class $\hat{\gamma}^{*}_{\textrm{dR}}(\omega) \, \mathbf{t}_{1} \cup \cdots \cup \mathbf{t}_{a}$, where $\mathbf{t}_{i}$ corresponds to $dT_{i} / T_{i}$ analogously to \autoref{tcalccor}. Using \autoref{epsendslem}, this class is the image of $y := \left( [\mathcal{L}_{1,\infty}] \cup \cdots \cup [\mathcal{L}_{a,\infty}] \right) \otimes \hat{\gamma}^{*}_{\textrm{dR}}(\omega)$.

Now suppose instead we compute the image of $\omega$ in the bottom left corner of the diagram by travelling leftwards along the top row and then down along the leftmost arrow. From the commutativity of the diagram, the result lands inside $y + \ker \sigma \subset y + J^{a,b} \otimes B_{\textrm{dR}}$. Thus we get that $\hat{\gamma}^{*}_{a,b}(\omega) = \hat{\alpha}^{*}_{a,b}(y) = \hat{\gamma}^{*}_{\textrm{dR}}(\omega) \otimes t^{a} $, where we apply \autoref{formalevallem} and the fact that the embedding $\mathbb{Z}_{p}(a) \hookrightarrow B_{\textrm{dR}}$ sends $\{ \zeta^{\star}_{\bullet} \} \otimes \cdots \otimes \{ \zeta^{\star}_{\bullet} \}$ to $t^{a}$. 
\end{proof}

\section{Realizing $G$-functions}

\label{realgfuncsec}

We now give our main technical result, which will give a cohomological interpretation of Andr\'e's $G$-functions at all places of our fixed number field $K$. We invite the reader to recall the conventions in \S\ref{notconvtpt1}, which are in force throughout this section. We additionally assume for simplicity that $s$ is defined on all of $S$ and that $ds$ trivializes $\Omega^{1}_{S}$, which becomes true after removing finitely many points from $S$. This in particular implies that $S$ is affine and $\mathcal{H}$ may be identified with its module of global sections. Finally, write $\mathcal{H}_{U} = R^{w} f_{*} \Omega^{\bullet}_{U/S}(\log (E \cap U))$, and note that there is a natural restriction map $\mathcal{H} \to \mathcal{H}_{U}$. 

We have a commuting diagram 
\begin{equation}
\label{keysquare}
\begin{tikzcd}[column sep=20ex, row sep=10ex]
U \arrow[r, "g"] \arrow[d, swap, "\restr{f}{U}"] & \Spec K[x_{1}, \hdots, x_{\nu}] \arrow[d, "j"] \\
S \arrow[r, "u"] & \Spec K[t]
\end{tikzcd}.
\end{equation}
The map $g$ is defined by ${(x_{1}, \hdots, x_{\nu}) \mapsto (z_{1}, \hdots, z_{\nu})}$, the map $u$ by $t \mapsto s$, and the map $j$ by $t \mapsto x_{1} \cdots x_{\mu}$. Because $g$ is \'etale (the coordinates $z_{1}, \hdots, z_{\nu}$ were chosen so that $dz_{1}, \hdots, dz_{\nu}$ trivializes $\Omega^{1}_{U}$), its image is an open $K$-algebraic subvariety $V \subset \Spec K[x_{1}, \hdots, x_{\nu}]$. Write $T \subset \Spec K[t]$ for the image of $V$. 

\begin{defn}
By a scaling of the coordinates $(z_{1}, \hdots, z_{\nu})$ we mean coordinates $(\lambda z_{1}, \hdots, \lambda z_{\nu})$ for some $\lambda \in K^{\times}$. By a scaling of (\ref{keysquare}) we mean coordinates $(\lambda z_{1}, \hdots, \lambda z_{\nu})$ for some $\lambda \in K^{\times}$ and $\lambda^{\mu} s$ for some $\lambda \in K^{\times}$. 
\end{defn}

We note that the diagram (\ref{keysquare}) continues to commute if one replaces it with a scaling. 

\begin{lem}
\label{scalesoball}
Choose a $K$-point $q \in g^{-1}(0)$. After replacing $(z_{1}, \hdots, z_{\nu})$ with a scaling by $N^{-1}$, where $N \in \mathbb{Z}$, the following property holds: for any embedding $\iota : K \hookrightarrow K_{v}$ with $v$ a finite place of $K$, the map $g$ is invertible in the open ball of $v$-adic radius $1$ around $0 \in \Spec K[x_{1}, \hdots, x_{\nu}]$ onto a neighbourhood containing $q$. (In particular, this ball is contained inside $V$.) This property continues to hold if $N$ is replaced by some multiple $N'$ of $N$.
\end{lem}

\begin{proof}
The idea is that one can write down a formal inverse to the map of germs $(U,q) \to (V,0)$ and have this inverse converge at each finite place in the desired neighbourhood after scaling coordinates. More explicitly, let us begin by embedding the affine variety $U$ as a closed subvariety of $\Spec K[y_{1}, \hdots, y_{\sigma}]$ defined by polynomials $p_{1}, \hdots, p_{\ell} \in K[y_{1}, \hdots, y_{\sigma}]$. After translation we may identify $q$ with the origin in $\Spec K[y_{1}, \hdots, y_{\sigma}]$. The map $g$ is then given by component polynomials $g_{1}, \hdots, g_{\nu} \in K[y_{1}, \hdots, y_{\sigma}]$ with no constant terms. The formal inverse $A$ we wish to compute is then given by power series
\begin{equation}
\label{Adefs}
A_{i}(x_{1}, \hdots, x_{\nu}) = \sum_{\boldsymbol{J}} A_{i,\boldsymbol{J}} \overline{x}^{\boldsymbol{J}}
\end{equation}
for $1 \leq i \leq \sigma$, where $\boldsymbol{J}$ ranges over the set $\mathcal{C}$ of all appropriate compositions of integers $\geq 0$ and we use multi-index notation to exponentiate the vector $\overline{x} = (x_{1}, \hdots, x_{\nu})$. 

The fact that $g \circ A = \textrm{id}$ and $p \circ A = 0$ gives a system of linear equations for each coefficient appearing in each $A_{i}$ in terms of coefficients of the polynomials $g_{j}$ and $p_{k}$. Any solution with $A(0) = q$ to this system defines a formal inverse to the map $(U, q) \to (V,0)$ on the level of completed formal power series rings, and the fact that $U \to V$ is \'etale implies that there is a unique such inverse; in particular, this system of equations together with the condition $A(0) = q$ has a unique formal solution. Let us view the solution to this system as a formal function $A$, and suppose that this formal function converges in the open ball of $v$-adic radius $1$ around $0$ for all finite places of $K$ outside of a finite set $\Sigma$. Then it will suffice to scale the coordinates $(x_{1}, \hdots, x_{\nu})$ by multiplying each $x_{i}$ by a sufficiently large integer $N$ whose prime factors all lie above places of $\Sigma$: indeed, doing so does not affect the radius of convergence for finite places outside of $\Sigma$, and the radius of convergence of the resulting power series at a place $v \in \Sigma$ will increase by a factor of $1/|N|_{v}$ and hence be greater than $1$ as soon as $|N|_{v}$ is small enough. We are reduced to the following more formal fact:

\begin{lem}
Suppose that we have formal power series (\ref{Adefs}) with coefficients in a number field $K$ that are uniquely determined by the property that $B_{i} \circ A = C_{i}$ and an initial condition $A(0) = q$, where $B_{1}, \hdots, B_{c} \in K[y_{1}, \hdots, y_{\sigma}]$, $C_{1}, \hdots, C_{c} \in K[x_{1}, \hdots, x_{\nu}]$ are finitely many polynomials with coefficients in $K$. Then $A$ converges in the open $v$-adic ball of radius $1$ away from finitely many places $v$ of $K$.
\end{lem}

\begin{proof}
To understand the system of equations defined by these polynomials we recall the multivariate Fa\`a di Bruno formula \cite{zbMATH00870188}, which says that the derivatives of a composition $C = B \circ A$ of functions given by power series centred at zero are given by
\begin{equation}
\label{faadibruno}
 (\partial_{\boldsymbol{J}} C)(0) = \sum_{\substack{1 \leq |\boldsymbol{\lambda}| \leq |\boldsymbol{J}|}} (\partial_{\boldsymbol{\lambda}} B)(0) \sum_{s=1}^{|\boldsymbol{J}|} \sum_{\mathcal{C}_{s}(\boldsymbol{J}, \boldsymbol{\lambda})} \boldsymbol{J}! \prod_{j = 1}^{s} \frac{[\boldsymbol{A}_{\boldsymbol{\ell}_{j}}(0)]^{\boldsymbol{k}_{j}}}{(\boldsymbol{k}_{j}!)[\boldsymbol{\ell}_{j}!]^{|\boldsymbol{k}_{j}|}} , 
 \end{equation}
where we have made use of the following notation:
\begin{itemize}
\item[-] the vectors $\boldsymbol{\lambda}$ and $\boldsymbol{k}_{j}$ come from $\mathbb{Z}^{\sigma}_{\geq 0}$ and the vectors $\boldsymbol{J}$ and $\boldsymbol{\ell}_{j}$ come from $\mathbb{Z}^{\nu}_{\geq 0}$;
\item[-] for any vector $\boldsymbol{u} = (u_{1}, \hdots, u_{r}) \in (\mathbb{Z}_{\geq 0})^{r}$ we have $|\boldsymbol{u}| = u_{1} + \cdots + u_{r}$;
\item[-] we have
\[ \mathcal{C}_{s}(\boldsymbol{J}, \boldsymbol{\lambda}) = \left\{ (\boldsymbol{k}_{1}, \hdots, \boldsymbol{k}_{s} ; \boldsymbol{\ell}_{1}, \hdots, \boldsymbol{\ell}_{s}) : \substack{|\boldsymbol{k}_{i}| > 0 , \hspace{1em} \boldsymbol{0}  \prec \boldsymbol{\ell}_{1} \boldsymbol \prec \cdots \prec \boldsymbol{\ell}_{s} \\ \sum_{i=1}^{s} \boldsymbol{k}_{i} = \boldsymbol{\lambda} \textrm{ and } \sum_{i=1}^{s} |\boldsymbol{k}_{i}| \boldsymbol{\ell}_{i} = \boldsymbol{J} } \right\} ,  \]
\item[-] for vectors $\boldsymbol{u} = (u_{1}, \hdots, u_{r})$ and $\boldsymbol{u}' = (u'_{1}, \hdots, u'_{r})$, the symbol $\boldsymbol{u} \prec \boldsymbol{u}'$ means that one of the following conditions holds:
\begin{itemize}
\item[(i)] $|\boldsymbol{u}| < |\boldsymbol{u}'|$;
\item[(ii)] $|\boldsymbol{u}| = |\boldsymbol{u}'|$ and $u_{1} < u'_{1}$; or
\item[(iii)] $|\boldsymbol{u}| = |\boldsymbol{u}'|$, $u_{1} = u'_{1}, \hdots, u_{k} = u'_{k}$ and $u_{k+1} < u'_{k+1}$ for some $1 \leq k < r$;
\end{itemize}
\item[-] the notation $\boldsymbol{A}_{\boldsymbol{\ell}}$ for $\boldsymbol{\ell} = (\ell_{1}, \hdots, \ell_{\nu})$ means $(\partial_{\boldsymbol{\ell}} A_{1}, \hdots, \partial_{\boldsymbol{\ell}} A_{\nu})$; and
\item[-] for a vector $\boldsymbol{u} = (u_{1}, \hdots, u_{r})$, we have $\boldsymbol{u}! = u_{1}! \cdots u_{r}!$.
\end{itemize}
The terms on the right-hand side of the equation (\ref{faadibruno}) involving the components of $\mathbf{A}_{\boldsymbol{J}}(0)$ are then
\begin{equation}
\boldsymbol{J}! \sum_{i=1}^{\sigma} (\partial_{i} B)(0) \frac{A_{i, \boldsymbol{J}}}{\boldsymbol{J}!} .
\end{equation}

Now let us suppose that $B$ is a polynomial over $K$, and that $(\partial_{\boldsymbol{J}} C)(0) = 0$; this is the case for $B \in \{ B_{1}, \hdots, B_{c} \}$ and for $|\mathbf{J}|$ sufficiently large. Then at a finite place $v$, for all but finitely many $v$, the norms $|(\partial_{\boldsymbol{\lambda}} B)(0)|_{v}$ are $\leq 1$ if $|\boldsymbol{\lambda}| \leq \deg B$, and equal to $0$ if $|\boldsymbol{\lambda}| > \deg B$. In particular, the equation (\ref{faadibruno}) induces the following linear equation for $\mathbf{A}_{\boldsymbol{J}} / \boldsymbol{J}!$, at least when $|\boldsymbol{J}| > \textrm{max} \{ \deg B, \deg C \}$:

\begin{equation}
\label{lineq}
\sum_{i=1}^{\sigma} (\partial_{i} B)(0) \frac{A_{i, \boldsymbol{J}}}{\boldsymbol{J}!} =  - \sum_{\substack{2 \leq |\boldsymbol{\lambda}| \leq \deg B}} (\partial_{\boldsymbol{\lambda}} B)(0) \sum_{s=1}^{|\boldsymbol{J}|} \sum_{\mathcal{C}_{s}(\boldsymbol{J}, \boldsymbol{\lambda})} \left( \prod_{j = 1}^{s} \frac{1}{(\boldsymbol{k}_{j}!)} \right) \left( \prod_{j=1}^{s} \frac{[\boldsymbol{A}_{\boldsymbol{\ell}_{j}}(0)]^{\boldsymbol{k}_{j}}}{[\boldsymbol{\ell}_{j}!]^{|\boldsymbol{k}_{j}|}} \right) .
\end{equation}

\noindent We note that there are only finitely many possibilities for the coefficients $\prod_{j=1}^{s} \frac{1}{(\boldsymbol{k}_{j}!)}$ which appear in (\ref{lineq}), and these possibilities are independent of $\mathbf{J}$ and depend only on $\deg B$. Indeed, because only finitely many $\boldsymbol{\lambda}$ ever occur in all such terms, the equation $\sum_{i=1}^{s} \boldsymbol{k}_{i} = \boldsymbol{\lambda}$ together with the condition $|\boldsymbol{k}_{i}| > 0$ for all $i$ ensures that only finitely values for the tuple $(s, \boldsymbol{k}_{1}, \hdots, \boldsymbol{k}_{s})$ ever appear, and hence there are only finitely many $\prod_{j=1}^{s} \frac{1}{(\boldsymbol{k}_{j}!)}$ which appear. After excluding a further finite set of norms $|\cdot|_{v}$, we may assume that all these coefficients have norm $1$, and the non-zero coefficients of $B$ and its derivatives also have norm $1$. Letting $B$ range over the finitely many polynomials in the set $B \in \{ B_{1}, \hdots, B_{c} \}$, we have proven that:

\begin{quote}
For $|\mathbf{J}| > \max_{1 \leq i \leq c} \max\{ \deg B_{i}, \deg C_{i} \}$, and for all but finitely many places $v$, the vector $\frac{1}{\boldsymbol{J}!} \boldsymbol{A}_{\boldsymbol{J}}$ is the solution to a system of linear equations 
\begin{equation}
\label{mateqforcoeff}
\boldsymbol{M} \frac{\boldsymbol{A}_{\boldsymbol{J}}}{\boldsymbol{J}!} = \boldsymbol{N} ,
\end{equation}
where $\boldsymbol{M}$ is a $c \times \sigma$ matrix, independent of $\boldsymbol{J}$, whose non-zero entries all have unit $v$-norm, and $\boldsymbol{N}$ is a vector with $c$ entries whose norms are at most 
\begin{equation}
\label{Alnorms}
\max_{s, \mathcal{C}(\boldsymbol{J}, \boldsymbol{\lambda})} \left\| \prod_{j=1}^{s} \frac{[\boldsymbol{A}_{\boldsymbol{\ell}_{j}}(0)]^{\boldsymbol{k}_{j}}}{[\boldsymbol{\ell}_{j}!]^{|\boldsymbol{k}_{j}|}} \right\|_{v} .
\end{equation}
\end{quote}

Let $M = \max_{1 \leq i \leq c} \max\{ \deg B_{i}, \deg C_{i} \}$. It is a formal fact that any collection of vectors $\mathbf{A}'_{\mathbf{J}}$ indexed by $\mathbf{J}$ which satisfy $\mathbf{A}'_{\mathbf{J}} = \mathbf{A}_{\mathbf{J}}$ for $|\mathbf{J}| \leq M$ and satisfy (\ref{lineq}) for $|\mathbf{J}| > M$ determine a formal power series $A'$ such that $A'(0) = q$ and $B_{i} \circ A' = C_{i}$ for all $i$. Since by hypothesis such an $A'$ is unique, there must be a unique solution to the infinite dimensional linear system defined by (\ref{lineq}) for $|\mathbf{J}| > M$ and $\mathbf{A}'_{\mathbf{J}} = \mathbf{A}_{\mathbf{J}}$ for $|\mathbf{J}| \leq M$. Then if $\boldsymbol{M}$ does not have rank $\sigma$, one can construct infinitely many solutions to (\ref{lineq}) for any fixed $\mathbf{J}$, and so recursively construct infinitely many solutions to the equations $B_{i} \circ A' = C_{i}$ of formal power series.

All this is to say that we may assume $c \geq \sigma$ or else $\mathbf{M}$ does not have rank $\sigma$, and then $c = \sigma$ by choosing a linearly-independent subset of the rows of $\mathbf{M}$. We now use this to show that, after possibly throwing out a further finite set of places $v$, one has $\| A_{i,\boldsymbol{J}} / \boldsymbol{J}! \|_{v} \leq 1$ for all $i$ and all $\boldsymbol{J}$. We prove this by induction, starting from the case where $|\mathbf{J}| = M$; we note that the base cases with smaller $|\mathbf{J}|$ can be assumed after removing a further finite set of places. Removing a further finite set of places to ensure that $\| \det(\boldsymbol{M}) \|_{v} = 1$, we may use Cramer's rule and the equation (\ref{mateqforcoeff}) to compute the entries of $A_{i,\boldsymbol{J}} / \boldsymbol{J}!$ as quotients $\det(\boldsymbol{M}') / \det(\boldsymbol{M})$, where $\boldsymbol{M}'$ is a matrix obtained from $\boldsymbol{M}$ by replacing a column with the vector $\boldsymbol{N}$. By induction the bound (\ref{Alnorms}) is at most $1$, so it follows that $\det(\boldsymbol{M}') / \det(\boldsymbol{M})$, and hence the entries of $A_{i,\boldsymbol{J}} / \boldsymbol{J}!$, have $v$-adic norm at most $1$.
\end{proof}
\end{proof}

\begin{thm}
\label{bigsupergfuncthm}
Scale coordinates as in \autoref{scalesoball}, and fix a point $q$ in the common vanishing locus $z_{1} = \cdots = z_{\mu} = 0$ with image $s_{0} \in S$. Then $q$ induces a $K$-linear map, compatible with base change along a finite extension $L/K$, 
\[ \Gamma : \mathcal{H}(S) \to K[[t]] , \]
whose image consists of $G$-functions. These $G$-functions satisfy the following two ``realization'' properties: 
\begin{itemize}
\item[(i)] Fix an embedding $\iota : K \hookrightarrow \mathbb{C}$, suppose that $\{ \omega_{i} \}_{i \in I}$ is a subset of $\mathcal{H}(S)$, and that $R > 0$ is a real number such that $\Gamma(\omega_{i})_{\iota}$ has radius of convergence at least $R$ for all $i \in I$. Denote by $\mathcal{D}_{R} \subset \an{S}_{\iota}$ the component containing $s_{0}$ of the complex analytic neighbourhood defined by $|s| < R$. Then for each point $s_{1} \in \mathcal{D}_{R} \setminus \{ s_{0} \}$, there exists a linear functional
\[ \gamma^{*}_{1} : H^{w}(\an{X}_{s_{1}}, \mathbb{Z}(w)) \to \mathbb{Z}(w) , \]
such that if $\rho$ is the natural isomorphism
\[H^{w}(\an{X}_{s_{1}}, \mathbb{Z}(w)) \otimes \mathbb{C} \xrightarrow{\sim} H^{w}_{\textrm{dR}}(\an{X}_{s_{1}}) , \]
then we have
\begin{equation}
\label{compangammaprop}
\frac{1}{(2 \pi i)^{w}} (\gamma^{*}_{1,\mathbb{C}} \circ \rho^{-1})(\omega_{i, s_{1}}) = (\Gamma(\omega_{i})_{\iota})(u(s_{1})) 
\end{equation}
for all $i \in I$.

\item[(ii)] Fix an embedding $\iota : K \hookrightarrow K_{v}$ for some finite place $v$ above the prime $p$, suppose that $\{ \omega_{i} \}_{i \in I}$ is a subset of $\mathcal{H}(S)$, and that $1 \geq R > 0$ is a real number such that $\Gamma(\omega_{i})_{\iota}$ has radius of convergence at least $R$ for all $i \in I$. Denote by $\mathcal{D}_{R} \subset \ad{S}_{\iota}$ the component containing $s_{0}$ of the adic neighbourhood defined by $|s| < R$. Then for each point $s_{1} \in \mathcal{D}_{R}(K_{v}) \setminus \{ s_{0} \}$ there exists, after possibly replacing $K_{v}$ with a finite extension, a neighbourhood $\Delta^{w,\nu-\mu}_{s_{1}} \subset \ad{X}_{s_{1}}$ such that the pair $(\Delta^{w,\nu-\mu}_{s_{1}}, \ad{X}_{s_{1}})$ satisfies the same hypotheses as the pair $(\Delta^{a,b}, Y)$ in \S\ref{funcextsec}. In particular, denoting by $\hat{\gamma}^{*}_{1, B_{\textrm{dR}}} = \hat{\gamma}^{*}_{a,b}$ the functional obtained by scalar extending the functional $\hat{\gamma}^{*}_{1} := \hat{\gamma}^{*}_{\textrm{\'et}}$ defined in \S\ref{funcextsec}, we have
\begin{equation}
\label{adicgammaprop}
\frac{1}{t^{w}} (\hat{\gamma}^{*}_{1, B_{\textrm{dR}}} \circ \rho^{-1})(\omega_{i,s_{1}}) = (\Gamma(\omega_{i})_{\iota})(u(s_{1})) .
\end{equation}
for all $i \in I$, where $\rho$ is the $p$-adic period isomorphism associated to $X_{s_{1}}$. 
\end{itemize}
\end{thm}

\begin{rem}
The $p$-adic period $t$ in the statement of (ii) is not to be confused with the coordinate $t$ in the diagram (\ref{keysquare}). 
\end{rem}

The remainder of this section is devoted to the proof of \autoref{bigsupergfuncthm}. We begin by constructing the map $\Gamma$. From functoriality, the sections $\{ \omega_{i} \}_{i \in I}$ all have restrictions to the sheaf $R^{w} f_{U,*} \Omega^{\bullet}_{U/S}(\log (E \cap U))$, which as we explained in \S\ref{canonextsec} is represented by the cohomology in degree $w$ of the complex 
\begin{equation}
\label{algUScomp}
 0 \to \mathcal{O}_{U} \to \Omega^{1}_{U/S}(\log (E \cap U)) \to \cdots \to \Omega^{n}_{U/S}(\log (U \cap E)) . 
\end{equation}
We may then further restrict to a formal neighbourhood of $q$ and consider the complex
\begin{equation}
\label{formalqcomp}
 0 \to \widehat{\mathcal{O}}_{U,(q)} \to \widehat{\Omega}^{1}_{U/S}(\log (E \cap U))_{(q)} \to \cdots \to \widehat{\Omega}^{n}_{U/S}(\log (E \cap U))_{(q)} , 
\end{equation}
and obtain a $K$-linear map 
\[ \eta : \mathcal{H}(S) \to \textrm{Cohom}^{w} \left[ \widehat{\Omega}^{\bullet}_{U/S}(\log(E \cap U))_{(q)} \right] .  \]
As we saw in \autoref{uniquereplem}, the target of $\eta$ is naturally a $1$-dimensional free module over $\widehat{\mathcal{O}}_{S,(s_{0})}$, and so we may define $\Gamma$ as the composition
\[ \mathcal{H}(S) \xrightarrow{\eta} \textrm{Cohom}^{w} \left[ \widehat{\Omega}^{\bullet}_{U/S}(\log(E \cap U))_{(q)} \right] \xrightarrow{\sim} \widehat{\mathcal{O}}_{S,s_{0}} \xrightarrow{u} \widehat{\mathcal{O}}_{\Spec K[t], 0} = K[[t]] . \]

Before turning to the proof of (i), we briefly explain why the image consists of $G$-functions. The point is that, within the degree-$w$ cohomology of the complex (\ref{algUScomp}), each relative form $\omega$ is represented by $h \frac{dz_{2} \cdots dz_{\mu}}{z_{2} \cdots z_{\mu}}$, where $h$ is a function algebraic over a rational function field; in particular, the power series $h^{q}$ in the coordinates $z_{1}, \hdots, z_{\mu}$ representing $h$ at $q$ is algebraic over a rational function field. The calculation of \autoref{uniquereplem}, which we will see again in the proof of (i), computes $\Gamma(\omega)$ as the $\mu$\emph{-diagonal} of $h^{q}$, i.e., the power series in one variable $t = z_{1} \cdots z_{\mu}$ obtained by keeping all terms $a \, z_{1}^{e_{1}} \cdots z_{\mu}^{e_{\mu}}$ with $e_{1} = \cdots = e_{\mu}$ and discarding the others. It is known \cite[I, \S4.2]{zbMATH00041964} that any function obtained in this way is a $G$-function, which Andr\'e himself uses in his proof in \cite[IX, \S4.4]{zbMATH00041964}.

\vspace{0.5em}

\begin{proof}[Proof of (i):]

We work entirely in the complex analytic category, and view the diagram (\ref{keysquare}) in the complex analytic category using base-change along $\iota : K \hookrightarrow \mathbb{C}$ and analytifying. Because we do not assume that $\mathcal{D}_{R}$ is simply-connected, we will first prove the theorem in a small neighbourhood of the degeneration, and then extend the result to $\mathcal{D}_{R}$ by analytic continuation. In particular, we start by fixing a sufficiently small analytic disk $\mathcal{D} \subset \mathcal{D}_{R}$ around $s_{0}$ and a small analytic disk $\mathcal{U} \subset f^{-1}(\mathcal{D})$ around $q$ such that the map $g$ maps $\mathcal{U}$ isomorphically to the open ball $\mathcal{B}_{\sqrt{\mu} r^{1/\mu}} \subset (\mathbb{A}^{n})^{\textrm{an}}$ of radius $\sqrt{\mu} r^{1/\mu}$, and the map $\mathcal{U} \to \mathcal{D}$ is identified with the map $\mathcal{B}_{\sqrt{\mu} r^{1/\mu}} \to \Delta$ given by $t \mapsto x_{1} \cdots x_{\mu}$, with $\Delta \subset (\mathbb{A}^{1})^{\textrm{an}}$ the open ball of radius $r$. We then obtain a complex
\begin{equation}
\label{anacomp}
0 \to \mathcal{O}_{\mathcal{U}} \to \Omega^{1}_{\mathcal{U}/\mathcal{D}}(\log (E \cap \mathcal{U})) \to \cdots \to \Omega^{n}_{\mathcal{U}/\mathcal{D}}(\log (\mathcal{U} \cap E)) ,
\end{equation}
which is identified with the complex appearing in \autoref{uniquereplem}. Thus, applying \autoref{uniquereplem}, one obtains, for each $i$, analytic functions $P_{i}$ on $\mathcal{D}$ such that $P_{i} \frac{dz_{2} \cdots dz_{\mu}}{z_{2} \cdots z_{\mu}}$ represents the restriction of $\omega_{i, \mathbb{C}}$ to the cohomology group $\textrm{Cohom}^{w} \left[ \Omega^{\bullet}_{\mathcal{U}/\mathcal{D}}(\log (E \cap \mathcal{U})) \right]$. As the $P_{i}$ are obtained using the same calculations that produced $\Gamma(\omega_{i})_{\iota}$, these functions agree at the formal level.

We now give a geometric interpretation of the $P_{i}$ following \cite[IX, \S4.4]{zbMATH00041964}, as follows. We consider the fibre $\mathcal{U}_{t_{1}} \subset \mathcal{U} = \mathcal{B}_{\sqrt{\mu} r^{1/\mu}}$ defined by $x_{1} \cdots x_{\mu} = t_{1}$ for some fixed point $t_{1} \in \mathcal{D} = \Delta$. We then construct an injective map $\ep_{t_{1}} : (S^{1})^{w} = ([0,2\pi]/\{ 0, 2 \pi \})^{w} \to \mathcal{U}_{t_{1}}$ by the formula
\begin{align*}
x_{1} &= R_{1}^{-w/\mu} t_{1} \exp(-i(\alpha_{1} + \cdots + \alpha_{w})) \\
x_{j} &= R_{1}^{1/\mu} \exp(i \alpha_{j-1}) & 2 \leq j \leq \mu \\
x_{j} &= 0 & \mu+1 \leq j \leq \nu ,
\end{align*}
where each $\alpha_{j} \in [0,2\pi]$ and $R_{1} = |t_{1}|$. We then consider the family of integrals
\[ P'_{i}(t) = \frac{1}{(2\pi i)^{w}} \int_{\ep_{t}} \restr{\omega_{i}}{\mathcal{U}_{t}} . \]
Representing the restriction $\restr{\omega_{i}}{\mathcal{U}_{t}}$ as a function of the form $h_{i} \frac{dx_{2} \cdots dx_{\mu}}{x_{2} \cdots x_{\mu}}$ with $h$ a power series in $x_{1}, x_{2}, \hdots, x_{\mu}$, one computes using the residue formula that $P'_{i}(t)$ is the $\mu$-\emph{diagonal} of $h$, i.e., the function whose power series is obtained from that of $h$ by substituting in $a \, t^{e_{1}}$ for all terms $a \, x_{1}^{e_{1}} \cdots x^{e_{\mu}}_{\mu}$ with $e_{1} = \cdots = e_{\mu}$, and ignoring all other terms. This is compatible with the calculation in the proof of \autoref{uniquereplem}, and we have that $P'_{i} = P_{i}$.

By taking the image of the cycles $\ep_{t}$ inside the fibres $X_{t}$, this calculation realizes the functions $P_{i}$ (and hence the functions $\Gamma(\omega_{i})_{\iota}$) as functions inside the image of the integration pairing
\[ R_{w} f_{*} \mathbb{Z}(-w) \otimes [R^{w} f_{*} \Omega^{\bullet}_{X'/S'}]^{\textrm{an}} \to \mathcal{O}_{S^{\textrm{an}}} \]
restricted to the neighbourhood $\mathcal{D}$. By analytic continuation, the cycles $\ep_{t}$ extend to a (possibly multi-valued, since we do not assume $\mathcal{D}_{R}$ is simply-connected) section $\widetilde{\ep}$ of $\restr{R_{w} f_{*} \mathbb{Z}}{\mathcal{D}_{R}}$, and hence produce a (a priori possibly multi-valued) function $\widetilde{P}_{i}$ inside $\restr{\mathcal{O}_{S^{\textrm{an}}}}{\mathcal{D}_{R}}$ after pairing with $\omega_{i}$. But since $\widetilde{P}_{i}$ agrees with the analytification of $\Gamma(\omega_{i})_{\iota}$ near $s_{0}$ and its power series representation converges on $\mathcal{D}_{R}$, the analytic function $\widetilde{P}_{i}$ is single-valued, and gives an analytic realization of $\Gamma(\omega_{i})_{\iota}$ on all of $\mathcal{D}_{R}$. 

To complete the proof of (i), it suffices to define, for each $s_{1} \in \mathcal{D}_{R}$, the functional $\gamma^{*}_{1}$. This we define as the evaluation functional 
\[ \gamma^{*}_{1,\mathbb{C}} : H^{w}(\an{X}_{s_{1}}, \mathbb{C}) \to \mathbb{C}, \hspace{1em} \left[ \int_{(-)} \omega  \mapsto \int_{\widetilde{\ep}_{s_{1}}} \omega \right] . \]
To make sense of this definition, we are using the canonical isomorphism $H^{w}(\an{X}_{s_{1}}, \mathbb{C}) \simeq H^{w}_{\textrm{top-dR}}(\an{X}_{s_{1}}) \otimes \mathbb{C}$ (with \emph{topological} de Rham cohomology) to represent each element of $H^{w}(\an{X}_{s_{1}}, \mathbb{C})$ as an integration functional, and then defining $\gamma^{*}_{1,\mathbb{C}}$ by evaluating this functional on $\widetilde{\ep}_{s_{1}}$. Note that because the section $\widetilde{\ep}$ may in principle be multi-valued, this is also true of the function $\gamma^{*}_{1,\mathbb{C}}$. However the equality (\ref{compangammaprop}), which amounts to the above observation that the relative period with $\omega_{i}$ is given by $\Gamma(\omega_{i})$, holds regardless of which choice we make. The desired equality (\ref{compangammaprop}) therefore holds, and the proof is complete. 
\end{proof}

\vspace{0.5em}

\begin{proof}[Proof of (ii):]
As before, we will regard all spaces as analytic spaces over $K_{v}$ using the fixed embedding $\iota : K \hookrightarrow K_{v}$.

We wish to mirror the argument we made in the complex analytic case. The argument is in many respects simplified by our choice of coordinates, which ensure that the property in \autoref{scalesoball} holds over the neighbourhood $\mathcal{D}_{R}$ since $R \leq 1$. In particular if one writes $\mathcal{D} \subset S$ for the component of the neighbourhood $|s| < 1$ containing $s_{0}$, and one writes $\mathcal{U} \subset f^{-1}(\mathcal{D})$ for the neighbourhood around $q$ defined by $|(z_{1}, \hdots, z_{\nu})|_{\iota} < 1$, then we are guaranteed that the family $\mathcal{U} \to \mathcal{D}$ is isomorphic to the family $\mathcal{V} \to \mathcal{E}$ obtained by restricting $j$ to the open ball of radius $1$ around $0$. 

Now let us construct a subspace isomorphic to $\Delta^{w,\nu-\mu}$ inside $X_{s_{1}} \cap \mathcal{U}$. Letting $t_{1} = u(s_{1})$, we may instead construct such a subspace inside $\mathcal{V}_{t_{1}}$. Write $R_{1} = |t_{1}|_{v} < R \leq 1$. Passing to a finite extension of $K_{v}$ if necessary, fix a $\mu$-root $t^{1/\mu}_{1}$ of $t_{1}$, the choice of which is unimportant. The fibre $\mathcal{V}_{t_{1}}$ is defined by the equation $x_{1} \cdots x_{\mu} = t_{1}$ inside the open ball defined by $\textrm{max}_{i} |x_{i}| < 1$, and we may embed $\Delta^{w,\nu-\mu}$ inside this neighbourhood via the map $x_{i} \mapsto t_{1}^{1/\mu} T_{i-1}$ for $2 \leq i \leq \nu$, and $x_{1} \mapsto t^{1/\mu}_{1} / ( T_{1} \cdots T_{w} )$. This embedding identifies $\Delta^{w,\nu-\mu}$ with the closed neighbourhood $\mathcal{V}^{w,\nu-\mu}_{t_{1}}$ inside $\mathcal{V}_{t_{1}}$ defined by 
\[ |x_{i}| = R_{1}^{1/\mu} \textrm{ for } 1 \leq i \leq \mu, \hspace{1em} \textrm{and} \hspace{1em} |x_{i}| \leq R_{1}^{1/\mu} \textrm{ for } i > \mu . \]
We may then define $\hat{\gamma}^{*}_{1}$ by pulling back the functional $\hat{\alpha}^{*}_{w,\nu-\mu}$ on the cohomology of $\Delta^{w,\nu-\mu}$ constructed in \S\ref{alphaconstrsec}, and then extend this to the entire fibre $X_{s_{1}}$ as in \S\ref{funcextsec}. Define $\hat{\gamma}^{*}_{\textrm{dR}} : H^{w}_{\textrm{dR}}(X^{\textrm{ad}}_{s_{1}}) \to K_{v}$ by the property that $\hat{\gamma}^{*}_{\textrm{dR}}(\omega) \frac{dx_{2} \cdots dx_{\mu}}{x_{2} \cdots x_{\mu}} = \restr{\omega}{\mathcal{V}^{w,\nu-\mu}_{t_{1}}}$, analogously to \S\ref{funcextsec}.

It now suffices to verify the desired equality. As before, we may restrict each $\omega_{i}$ to the cohomology in degree $w$ of the complex $\Omega^{\bullet}_{\mathcal{U}_{R}/\mathcal{D}_{R}}$ to obtain, via \autoref{uniquereplem}, sections $h_{i} \frac{dz_{2} \cdots dz_{\mu}}{z_{2} \cdots z_{\mu}}$, with $h_{i}$ a function on $\mathcal{D}_{R}$. We first observe that the $h_{i}$ agree with the functions $\Gamma(\omega_{i})_{\iota} \circ u$ on $\mathcal{D}_{R}$: this is true by construction at the formal level, and that implies that these functions agree in a small enough ball around $s_{0}$, so this follows by uniqueness of analytic continuation of rigid-analytic functions on affinoid balls (c.f. \cite[Prop. 1.6.24]{mehmeti2019rhoatching}). Now the coherent cohomology group $H^{w}(\mathcal{U}_{s_{1}, \textrm{p\'et}}, \Omega^{\bullet}_{\mathcal{U}_{s_{1}}})$ is computed by the cohomology of the na\"ive de Rham complex $\Omega^{\bullet}_{\mathcal{U}_{s_{1}}}$. The natural restriction map 
\[ H^{w}(\mathcal{U}_{\textrm{p\'et}}, \Omega^{\bullet}_{\mathcal{U}/\mathcal{D}}(\log (E \cap \mathcal{U}))) \to H^{w}(\mathcal{U}_{s_{1}, \textrm{p\'et}}, \Omega^{\bullet}_{\mathcal{U}_{s_{1}}}) \]
is then represented by the natural map 
\[ \textrm{Cohom}^{w}\left[ \Omega^{\bullet}_{\mathcal{U}/\mathcal{D}}(\log (E \cap \mathcal{U})) \right] \to \textrm{Cohom}^{w}\left[ \Omega^{\bullet}_{\mathcal{U}_{s_{1}}} \right] , \]
and this map sends $h_{i} \frac{dz_{2} \cdots dz_{\mu}}{z_{2} \cdots z_{\mu}}$ to $h_{i}(s_{1}) \frac{dz_{2} \cdots dz_{\mu}}{z_{2} \cdots z_{\mu}}$. If we then evaluate this element with $\gamma_{1,B_{\textrm{dR}}}^{*} \circ \rho^{-1}$, one sees by \autoref{indepofpathlem} that the result is simply $t^{w} h_{i}(s_{1})$, which completes the proof. 
\end{proof}

\section{Algebraic Relations on Functionals}
\label{algrelsec}

We now introduce a new way to obtain relations on Andr\'e's $G$-functions at finite places, which uses the explicitly $p$-adic nature of our construction. The first key observation is the following:

\begin{lem}
Let $Y$ be a proper algebraic variety over the finite extension $K_{v}$ of $\mathbb{Q}_{p}$, and let $\hat{\gamma}^{*} : H^{a}(Y_{\overline{K_{v}}, \textrm{p\'et}}, \hat{\mathbb{Z}}_{p}(a)) \to \mathbb{Z}_{p}(a)$ be the map $\hat{\gamma}^{*} = \hat{\gamma}^{*}_{\textrm{\'et}}$ constructed in \S\ref{funcextsec}. Then the Galois group $G_{K_{v}} = \textrm{Gal}(\overline{K_{v}}/K_{v})$ acts on $\hat{\gamma}^{*}$ through $\chi^{a}_{\textrm{cycl}}$, where $\chi_{\textrm{cycl}}$ is the cyclotomic character. 
\end{lem}

\begin{proof}
Using \autoref{indepofpathlem} we have the formula $\hat{\gamma}^{*}_{B_{\textrm{dR}}} = \hat{\gamma}^{*}_{\textrm{dR}} \otimes t^{a}$ after scalar extension. But $\hat{\gamma}^{*}_{\textrm{dR}}$ is $\textrm{Gal}(\overline{K_{v}}/K_{v})$-invariant, and $\textrm{Gal}(\overline{K_{v}}/K_{v})$ acts on $t^{a}$ through $\chi^{a}_{\textrm{cycl}}$, so the result follows. 
\end{proof}

We now adopt a slightly more abstract perspective to produce relations on de Rham coordinates of $\hat{\gamma}^{*}$. The situation of interest is as follows.


\begin{sit}
\label{relprodsit}
We have a $B_{\textrm{dR}}$-admissible $G_{K_{v}}$-representation of weight $-a$ on a $\mathbb{Q}_{p}$-vector space $V_{\textrm{\'et}}$ such that $V_{\textrm{dR}, K_{v}} := (V_{\textrm{\'et}} \otimes_{\mathbb{Q}_{p}} B_{\textrm{dR}})^{G_{K_{v}}}$ admits a $K$-structure $V_{\textrm{dR}} \otimes_{K} K_{v} \simeq V_{\textrm{dR}, K_{v}}$, where $K \subset K_{v}$ is a finite extension of $\mathbb{Q}$. We also have
\begin{itemize}
\item[-] a $G_{K_{v}}$-invariant endomorphism $\tau : V_{\textrm{\'et}} \to V_{\textrm{\'et}}$, whose de Rham realization is defined over $K$; 
\item[-] a functional $\hat{\gamma}^{*} : V_{\textrm{\'et}} \to \mathbb{Q}_{p}(a)$ on which $G_{K_{v}}$ acts through the $a$'th power of $\chi_{\textrm{cycl}}$; and
\item[-] the dimension $k$ of the $G_{K_{v}}$-invariant subspace of
\[ \Hom(V_{\textrm{\'et}}, \mathbb{Q}_{p}) \simeq \Hom(V_{\textrm{\'et}}, \mathbb{Q}_{p}(a)) \otimes_{\mathbb{Q}_{p}} \mathbb{Q}_{p}(-a) \]
is less than the degree of the minimal polynomial of $\tau$.
\end{itemize}
We have a fixed $K$-basis $\omega_{1}, \hdots, \omega_{m}$ for $V_{\textrm{dR}}$. 
\end{sit}

In the situation above, the most typical application will be with $Y$ a variety defined over a number field $K$, and $V_{\textrm{\'et}} = H^{a}(Y_{\overline{K}, \textrm{\'et}}, \mathbb{Q}_{p}(a))$, $V_{\textrm{dR}} = H^{w}(Y, \Omega^{\bullet}_{Y})$ its \'etale and de Rham cohomology groups, although the extra generality will be useful to handle summands appearing in such cohomology groups as well.

\begin{prop}
\label{relprop}
Suppose we are in the situation of \autoref{relprodsit}. Then there is a non-zero $K$-algebraic relation on the dual coordinates of $\hat{\gamma}^{*}_{B_{\textrm{dR}}}$ of degree equal to $k+1$. 
\end{prop}

\begin{proof}
It suffices to show that the set
\begin{equation}
\label{lindepvecs}
 \hat{\gamma}^{*}_{B_{\textrm{dR}}}, \hat{\gamma}^{*}_{B_{\textrm{dR}}} \circ \tau, \hdots, \hat{\gamma}^{*}_{B_{\textrm{dR}}} \circ \tau^{k} 
\end{equation}
is linearly dependent. Indeed, if we evaluate the vectors in this set on the basis $\omega_{1}, \hdots, \omega_{m}$ we will obtain a sequence of vectors $v_{1}, \hdots, v_{k+1} \in B_{\textrm{dR}}^{m}$ such that $v_{i}$ is obtained by applying a non-zero linear transformation to $v_{1}$ whose matrix entries lie in $K$. But we may then get a $K$-algebraic relation on the coordinates of $v_{1}$ by taking the determinant of a square submatrix of $[v_{1} | \cdots | v_{k+1}]$. This relation has degree $k+1$ and depends only on the coordinates of $\tau$. That this relation is non-zero follows from the definition of minimal polynomial and the fact that $k$ is assumed smaller than its degree in \autoref{relprodsit}.

To show this linear dependence, we observe that $\tau$ is invariant under $G_{K_{v}}$, and therefore the vectors (\ref{lindepvecs}) are all obtained from scalar extension of vectors in $\Hom(V_{\textrm{\'et}}, \mathbb{Q}_{p}(a))$. Twisting by $\mathbb{Q}_{p}(-a)$, the twisted vectors are all invariant under the action of $G_{K_{v}}$ on $\Hom(V_{\textrm{\'et}}, \mathbb{Q}_{p})$. But the space of $G_{K_{v}}$-invariants on $\Hom(V_{\textrm{\'et}}, \mathbb{Q}_{p})$ has dimension at most $k$ by assumption. 
\end{proof}

\begin{prop}
\label{charpolyrelversion}
Suppose we are in the setting of \autoref{relprodsit}, and that the characteristic polynomial $P$ of $\tau$ is equal to its minimal polynomial. Then there is a non-zero $K$-algebraic relation of degree at most $(\deg P)!$ on the dual coordinates of $\hat{\gamma}^{*}_{B_{\textrm{dR}}}$ which is a product of linear relations, with each linear factor defined over $K F$, with $F$ the splitting field of $P$. The relation depends only on the coordinates of $\tau$ in the fixed basis $\omega_{1}, \hdots, \omega_{m}$. 
\end{prop}

\begin{proof}
That the characteristic polynomial of $\tau$ is equal to its minimal polynomial means that $V_{\textrm{dR}}^{*}$ admits only finitely many $\tau$-invariant subspaces, each of which is defined over the splitting field of $P$. The proof of \autoref{relprop} shows that $\hat{\gamma}^{*}_{B_{\textrm{dR}}}$ lies inside one of these subspaces, so one may take as a relation a product of a linear relation associated to this subspace with each of its Galois conjugates.
\end{proof}

\begin{lem}
\label{boundbyhodgenumber}
In the setting of \autoref{relprodsit}, the dimension of the $G_{K_{v}}$-invariant subspace of $\Hom(V_{\textrm{\'et}}, \mathbb{Q}_{p})$ is at most $\dim_{K} \textrm{gr}_{0} V^{*}_{\textrm{dR}}$.
\end{lem}

\begin{proof}
Letting $\mathbb{C}_{v}$ be the completion of the algebraic closure of $K_{v}$, the module $\Hom(V_{\textrm{\'et}}, \mathbb{Q}_{p})$ becomes isomorphic to $\Hom(V_{\textrm{\'et}} \otimes \mathbb{C}_{v}, \mathbb{C}_{v})$ after scalar extension, which the Hodge-Tate comparison shows is isomorphic to a sum $\bigoplus_{i} \textrm{gr}_{i} V_{\textrm{dR}}^{*} \otimes_{K} \mathbb{C}_{v}(-i)$. The subspace of $\Hom(V_{\textrm{\'et}}, \mathbb{Q}_{p})$ spanned by Galois invariants is then mapped into $\textrm{gr}_{0} V^{*}_{\textrm{dR}} \otimes_{K} \mathbb{C}_{v}$. 
\end{proof}

In the situation where the endomorphism $\tau$ appearing in \autoref{relprop} is defined over a finite extension $L$ of $K$, we will also want to control the degree of the extension $[L : K]$, for which the following fact, proven in \cite{papas2022height}, will be useful:

\begin{defn}
If $Y$ is an algebraic variety over $\overline{\mathbb{Q}}$, we say an element $\tau$ of $\textrm{End}(H^{a}(Y, \Omega^{\bullet}_{Y}))$ is an absolute Hodge endomorphism if the Betti-de Rham comparison $\rho_{\iota} : H^{a}(Y, \Omega^{\bullet}_{Y}) \otimes_{\iota} \mathbb{C} \xrightarrow{\sim} H^{a}(Y_{\iota}, \mathbb{C})$ identifies $\tau$ with an endomorphism of the rational Hodge structure on $H^{a}(Y_{\iota}, \mathbb{Q})$ under every choice of embedding $\iota : \overline{\mathbb{Q}} \hookrightarrow \mathbb{C}$. We say the absolute Hodge conjecture for endomorphisms holds for $Y$ in degree $a$ if there exists some fixed $\iota$ such that the set of absolute Hodge endomorphisms of $Y$ is identified with the image under $\rho^{-1}_{\iota} \circ (-) \circ \rho_{\iota}$ of the endomorphisms of the Hodge structure on $H^{a}(Y_{\iota}, \mathbb{Q})$. 
\end{defn}

\begin{prop}
\label{endoalgextprop}
Suppose that $Y$ is an algebraic variety defined over a number field $K$, and that the absolute Hodge conjecture for endomorphisms holds for $Y_{\overline{\mathbb{Q}}}$ in degree $a$. Then the algebra of absolute Hodge endomorphisms may be identified with a subalgebra of $\textrm{End}(H^{a}(Y, \Omega^{\bullet}_{Y}))_{L}$, where $L/K$ is a finite Galois extension with degree bounded only in terms of $m = \dim_{K} H^{a}(Y, \Omega^{\bullet}_{Y})$.
\end{prop}

\begin{proof}
This is \cite[Prop. 5.1]{papas2022height} and \cite[Prop. 5.2]{papas2022height}. Note that the absolute Hodge conjecture is only assumed for endomorphisms associated to $Y$ in the proof. 
\end{proof}

\section{Height Bounds for Families over Curves}
\label{heightboundsec}

In this section we prove \autoref{onlycontrolinfplaces} and \autoref{sndsinggivesheights}. In sections \S\ref{monocompframesec} and \S\ref{funcrelgencase} we work in an abstract setting starting with an arbitrary complex local system $\mathbb{L}$; in particular, we do not use the setup of \S\ref{notconvtpt1}. Starting in \S\ref{appoffuncreltoGfunc} we then specialize the discussion back to the setting of \S\ref{notconvtpt1}, which is then adjusted in \S\ref{chooseparamsec} and \S\ref{addsetupsec} to the setup required for the main arguments, which are carried out in \S\ref{fstproofsec6} and \S\ref{sndproofsec6}. 

\subsection{Period Torsors}
\label{monocompframesec}

To rule out functional relations on collections of $G$-functions, we will require the introduction of a period torsor $P'$ which we construct using the Riemann-Hilbert correspondence.

\begin{defn}
\label{algmonodef}
Let $T$ be a connected smooth complex algebraic variety, and $\mathbb{L}$ a complex local system on $T$. Then the algebraic monodromy group $\mathbf{H}_{\mathbb{L}}$ of $(T, \mathbb{L})$ is the abstract algebraic group defined as the identity component of the Zariski closure
\begin{equation}
\label{algmonoeq}
\overline{\textrm{im}[\pi_{1}(T,t) \to \GL(\mathbb{L}_{t})]}^{\textrm{Zar}} .
\end{equation}
\end{defn}

\noindent If $Y$ is an irreducible normal complex algebraic variety and $U \subset Y$ is a Zariski open subvariety, then the natural map $\pi_{1}(U) \to \pi_{1}(Y)$ is surjective (see \cite[0.7 (B)]{zbMATH03760312}). As a consequence of this, the formation of the algebraic monodromy group for a pair $(T, \mathbb{L})$ as in \autoref{algmonodef} is stable under replacing $T$ with an open subvariety $U$ (resp. a finite covering $C \to T$), and $\mathbb{L}$ with its restriction to $U$ (resp. $C$). 

Now recall the statement of the Riemann-Hilbert correspondence (c.f. the first paragraph of \cite{esnault2019survey}): 

\begin{thm}[Thm. 5.9 in \cite{Deligne1970}]
Let $T$ be a smooth irreducible complex algebraic variety, and consider the functor 
\[ \operatorname{Sol} : \operatorname{MIC}_{\textrm{reg}}(T) \to \mathbb{C}\operatorname{LocSys}(T^{\textrm{an}}) \]
which sends a pair $(\mathcal{V}, \nabla)$ consisting of a locally free $\mathcal{O}_{T}$-module $\mathcal{V}$ with regular singular integral connection $\nabla$ to its associated complex local system $\mathbb{L} := \operatorname{Sol}(\mathcal{V}, \nabla)$. Then $\operatorname{Sol}$ is an equivalence of neutral Tannakian categories.
\end{thm}

\begin{notn}
Given a local system $\mathbb{L}$ we write $\mathbb{L}^{a,b} = \mathbb{L}^{\otimes a} \otimes (\mathbb{L}^{*})^{\otimes b}$, where $(-)^{*}$ takes $\mathbb{L}$ to its dual. We have the analogous definition of $(\mathcal{V}^{a,b}, \nabla^{a,b})$ given an object $(\mathcal{V}, \nabla)$ of $\operatorname{MIC}_{\textrm{reg}}$, and for $L^{a,b}$ when $L$ is a fixed vector space.
\end{notn}

\begin{defn}
Given a locally free $\mathcal{O}_{T}$-module $\mathcal{V}$, we write $\mathbb{B}(\mathcal{V}) \to T$ for the associated algebraic vector bundle: here $\mathbb{B}(\mathcal{V})$ is an algebraic variety locally of the form $U \times \mathbb{A}^{m}$ over each Zariski open $U \subset T$ where $\restr{\mathcal{V}}{U}$ is free, satisfying the usual properties.
\end{defn}

We now construct the ``period torsor'' $P'$ which will be of interest. Start by replacing $T$ with a finite \'etale covering so that $\mathbf{H}_{\mathbb{L}}$ agrees with the expression in (\ref{algmonoeq}); that is, so that the algebraic group defined by (\ref{algmonoeq}) is connected. Recall that, by a theorem of Chevalley (e.g. \cite[I, Prop. 3.1]{zbMATH03728195}), for any algebraic subgroup $H \subset \GL(L)$, with $L$ a fixed $\mathbb{C}$-vector space, there exists $a, b \geq 0$ and a $1$-dimensional subspace $M \subset L^{a,b}$ such that $H$ is the stabilizer of $M$. Let $\mathbb{L}$ be the local system associated to $(\mathcal{V}, \nabla)$. Using the equivalence between local systems on $T$ and monodromy representations, one concludes that there exists $a, b \geq 0$ and a local subsystem $\mathbb{M} \subset \mathbb{L}^{a,b}$ such that, for each $t \in T$, the realization of $\mathbf{H}_{\mathbb{L}}$ inside $\GL(\mathbb{L}_{t})$ is the stabilizer of $\mathbb{M}_{t}$. Let $\mathcal{M} \subset \mathcal{V}^{a,b}$ be the locally free sheaf associated to $\mathbb{M}$ by Riemann-Hilbert. (From here on out, we will not use the assumption that $\mathbb{M}$ is $1$-dimensional.)

Fixing a basepoint $c \in T$, we now consider the Hom-sheaf $\sheafhom(\mathcal{V}, \mathcal{O}_{T} \otimes_{\mathbb{C}} \mathcal{V}_{c})$, which is a locally free $\mathcal{O}_{T}$-module constructed as in \cite[\href{https://stacks.math.columbia.edu/tag/01CM}{Section 01CM}]{stacks-project}, and whose associated algebraic bundle $\mathbb{B}(\sheafhom(\mathcal{V}, \mathcal{O}_{T} \otimes_{\mathbb{C}} \mathcal{V}_{c}))$ has as points pairs $(t, \varphi)$ where $t \in T$ and $\varphi : \mathcal{V}_{t} \to \mathcal{V}_{c}$ is a linear map. We then define
\begin{equation}
\label{Pdefeq}
P' := \{ (t, \varphi) \in \mathbb{B}(\sheafhom(\mathcal{V}, \mathcal{O}_{T} \otimes_{\mathbb{C}} \mathcal{V}_{c})) : \varphi(\mathcal{M}_{t}) = \mathcal{M}_{c}, \hspace{0.5em} \varphi \textrm{ invertible} \} .
\end{equation}
Because $\mathcal{M}$ is algebraic, $P'$ is an algebraic subvariety of $\mathbb{B}(\sheafhom(\mathcal{V}, \mathcal{O}_{T} \otimes_{\mathbb{C}} \mathcal{V}_{c}))$. It also admits a natural action by the algebraic group $\mathbf{H}_{\mathbb{L}}$, which we regard through its realization inside $\GL(\mathcal{V}_{c})$ and which acts on $(t, \varphi)$ via $(t, \varphi) \mapsto (t, \alpha \circ \varphi)$ for $\alpha \in \mathbf{H}_{\mathbb{L}}(\mathbb{C})$. Using the relationship between $\mathcal{M}$ and $\mathbf{H}_{\mathbb{L}}$ one checks this makes $P'$ into an $\mathbf{H}_{\mathbb{L}}$-torsor.

\subsection{Ruling out Functional Relations}
\label{funcrelsec}

\subsubsection{For General Local Systems}
\label{funcrelgencase}

We continue with the notation established in \S\ref{monocompframesec}, and let $B \subset T(\mathbb{C})$ be an analytic ball around $c \in T(\mathbb{C})$. 

\begin{prop}
\label{axschancor}
Let $\mathcal{P} \subset P'(\mathbb{C})$ be the image of the unique section $\textrm{id} \in \sheafhom(\mathbb{L}, \mathbb{L}_{c})(B)$ which extends the identity $\textrm{id}_{c} : \mathbb{L}_{c} \to \mathbb{L}_{c}$. Then $\mathcal{P}$ is Zariski dense in $P'$.
\end{prop}

\begin{proof}
We argue analogously to \cite[Lem. 7.2]{papas2022height} and \cite[Lem 2.7]{bakker2023functional}. Analyticially continue $\mathcal{P}$ to a closed irreducible analytic locus $\widetilde{\mathcal{P}} \subset P'$, and let $Z$ be the Zariski closure of $\mathcal{P}$ (and hence of $\widetilde{\mathcal{P}}$). Then because $\widetilde{\mathcal{P}}$ is obtained by parallel translation via the action of the monodromy group $\Gamma \subset \mathbf{H}_{\mathbb{L}}(\mathbb{C})$, the locus $Z$ contains an orbit of $\Gamma$ within each fibre of $P' \to T$. Since $P'$ is an $\mathbf{H}_{\mathbb{L}}(\mathbb{C})$-torsor and $\Gamma$ is Zariski dense in $\mathbf{H}_{\mathbb{L}}(\mathbb{C})$, this means that $Z$ contains each fibre of $P' \to T$, and hence contains $P'$.
\end{proof}

We now consider a fixed linear functional $\alpha^{*}_{1} : \mathbb{L}_{c} \to \mathbb{C}$ which we extend to a basis $\alpha^{*}_{1}, \hdots, \alpha^{*}_{m}$ for the dual space $(\mathbb{L}_{c})^{*}$, and write $\alpha_{1}, \hdots, \alpha_{m}$ for the corresponding basis of $\mathbb{L}_{c}$. We will also use the same notation for the corresponding elements of $\mathcal{V}^{*}_{c}$ and $\mathcal{V}_{c}$. Then we obtain an induced isomorphism $\mathcal{V}_{c} \cong \mathbb{C}^{m}$, and therefore a decomposition of locally free $\mathcal{O}_{T}$-modules
\[ \sheafhom(\mathcal{V}, \mathcal{O}_{T} \otimes_{\mathbb{C}} \mathcal{V}_{c}) \cong \sheafhom(\mathcal{V}, \mathcal{O}_{T} \otimes_{\mathbb{C}} \mathbb{C}^{m}) \cong (\mathcal{V}^{*})^{m} . \]
We get induced projection maps $p_{i} : \mathbb{B}(\sheafhom(\mathcal{V}, \mathcal{O}_{T} \otimes_{\mathbb{C}} \mathcal{V}_{c})) \to \mathbb{B}(\mathcal{V}^{*})$.

\begin{notn}
We write $\Gamma_{1}$ for the the image in $\mathbb{B}(\mathcal{V}^{*})$ of the unique extension of $\alpha^{*}_{1}$ to a section of $\mathbb{L}^{*}(B)$.
\end{notn}

\begin{lem}
\label{imagegivesGfunc}
We have $p_{1}(\mathcal{P}) = \Gamma_{1}$. 
\end{lem}

\begin{proof}
Immediate from the definitions.
\end{proof}

\begin{notn}
For a polynomial $Q$, we write $V(Q)$ for its vanishing locus.
\end{notn}

\begin{defn}
\label{Qreldef}
A relation on $\alpha^{*}_{1}$ is defined to be a homogeneous polynomial function $Q$ on $\mathcal{V}^{*}_{c}$ such that $V(Q)$ contains $\alpha^{*}_{1}$. It is said to be non-trivial if there does not exist a closed algebraic subvariety $\widetilde{R} \subset \mathbb{B}(\mathcal{V}^{*})$ containing $\Gamma_{1}$ whose intersection with $\mathcal{V}^{*}_{c}$ agrees with $V(Q)$. Otherwise it is said to be trivial.
\end{defn}

\begin{lem}
\label{trivialchar}
Suppose that $Q$ is a trivial relation on $\alpha^{*}_{1}$. Then the orbit $\mathbf{H}_{\mathbb{L}}(\mathbb{C}) \cdot \alpha^{*}_{1} \subset \mathbb{L}^{*}_{c}$ lies inside $V(Q)$. 
\end{lem}

\begin{proof}
Take $\widetilde{R}$ as in (\ref{Qreldef}). By combining \autoref{axschancor} and \autoref{imagegivesGfunc}, the inverse image $p^{-1}_{1}(\widetilde{R})$ contains $P'$. Specializing at $c$, this means that $p^{-1}_{1}(V(Q))$ contains the fibre $P'_{c}$ of $P' \to T$ over $c$. The fibre $P'_{c}$ is, by construction, the $\mathbf{H}_{\mathbb{L}}(\mathbb{C})$ orbit of $\textrm{id}_{c} \in \Hom(\mathcal{V}_{c}, \mathcal{V}_{c})$ (it evidently makes no difference whether we act on the left or the right). Thus, one learns that for any element $g \in \mathbf{H}_{\mathbb{L}}(\mathbb{C})$ the projection $p_{1}(\textrm{id}_{c} \circ g^{-1})$ lies in the vanishing locus of $Q$. From the definitions one has 
\begin{align*}
p_{1}(\textrm{id}_{c} \circ g^{-1}) &= p_{1}\left( (\alpha^{*}_{1} \circ g^{-1}) \otimes \alpha_{1} + \cdots + (\alpha^{*}_{m} \circ g^{-1}) \otimes \alpha_{m} \right) = g \cdot \alpha^{*}_{1} .
\end{align*}
\end{proof}
    
\begin{prop}
\label{relisnontrivial}
Let $Q$ be a relation on $\alpha^{*}_{1}$, and suppose that $Q$ is a product of linear relations. Let $E \subset \mathbb{C}$ be a field, $\mathbb{L}_{E}$ an $E$-structure on $\mathbb{L}$ for which $\mathbb{L}^{*}_{E}$ is $E$-simple, and suppose that $\alpha^{*}_{1}$ is defined over $E$. Then $Q$ is non-trivial.
\end{prop}

\begin{proof}
Assume that $Q$ is trivial. Then by \autoref{trivialchar} $Q$ vanishes on the orbit $\mathbf{H}_{\mathbb{L}}(\mathbb{C}) \cdot \alpha^{*}_{1}$. Since $\mathbf{H}_{\mathbb{L}}(\mathbb{C})$ is irreducible as an algebraic set, $\mathbf{H}_{\mathbb{L}}(\mathbb{C}) \cdot \alpha^{*}_{1}$ gives an irreducible algebraic subset of $\mathcal{V}^{*}_{c}$, and therefore some linear factor $B_{1}$ of $Q$ vanishes on $\mathbf{H}_{\mathbb{L}}(\mathbb{C}) \cdot \alpha^{*}_{1}$. If one has two elements $\beta^{*}, \beta'^{*} \in \mathcal{V}^{*}_{c}$ such that $B_{1}$ vanishes on $\mathbf{H}_{\mathbb{L}}(\mathbb{C}) \cdot \beta^{*}$ and $\mathbf{H}_{\mathbb{L}}(\mathbb{C}) \cdot \beta'^{*}$, then, because $B_{1}$ is linear, the same is true for $\beta^{*} + \beta'^{*}$. It follows that $W := \textrm{span}_{\mathbb{C}} \{ \mathbf{H}_{\mathbb{L}}(\mathbb{C}) \cdot \alpha^{*}_{1} \}$ is a strict subspace of $\mathbb{L}^{*}_{c}$ which is invariant under $\mathbf{H}_{\mathbb{L}}$. But because $\alpha^{*}_{1}$ is defined over $E$, this contradicts the $E$-simplicity of $\mathbb{L}^{*}_{E}$. It follows that $Q$ is non-trivial.
\end{proof}

\begin{prop}
\label{relisnontrivial2}
Assume the same setup as \autoref{relisnontrivial}. Let $\mathbb{W}$ be a second complex local system on $T$ with associated module with regular singular integral connection $(\mathcal{W}, \delta)$. Set $\mathbb{M} = \mathbb{L} \oplus \mathbb{W}$ and $\mathcal{M} = \mathcal{V} \oplus \mathcal{W}$. Then we may consider an element $(\alpha^{*}_{1}, \zeta^{*}_{1}) \in \mathcal{M}_{c}^{*} = \mathcal{V}^{*}_{c} \oplus \mathcal{W}^{*}_{c}$, and let $\overline{Q}$ be the pullback of $Q$ to $\mathcal{M}^{*}_{c}$. Then $\overline{Q}$ is a non-trivial relation on $(\alpha^{*}_{1}, \zeta^{*}_{1})$.
\end{prop}

\begin{proof}
Supposing $\overline{Q}$ is trivial, we may apply \autoref{trivialchar} to $(\mathbb{M}, \overline{Q}, (\alpha^{*}_{1}, \zeta^{*}_{1}))$ and learn that $\mathbf{H}_{\mathbb{M}}(\mathbb{C}) \cdot (\alpha^{*}_{1}, \zeta^{*}_{1}) \subset V(\overline{Q})$. But $\mathbf{H}_{\mathbb{M}}$ acts on each factor separately, and on $\mathbb{L}^{*}_{c} \subset \mathbb{M}^{*}_{c}$ through a surjective map $\rho : \mathbf{H}_{\mathbb{M}} \to \mathbf{H}_{\mathbb{L}}$ (the image of $\rho$ is an algebraic group containing the image of $\pi_{1}(T)$, hence equal to $\mathbf{H}_{\mathbb{L}}$), so in fact we obtain $\mathbf{H}_{\mathbb{L}}(\mathbb{C}) \cdot \alpha^{*}_{1} \subset \mathcal{V}^{*}_{c} \cap V(\overline{Q}) = V(Q)$. The proof then proceeds as in \autoref{relisnontrivial}.
\end{proof}

\subsubsection{Application to $G$-functions}
\label{appoffuncreltoGfunc}

We now apply the above discussion to our situation. For this we take $\mathbb{L} = \mathbb{V}'_{\mathbb{C}}$, $\mathcal{V} = \mathcal{H}'$ with its Gauss-Manin connection, and $T = S'$. Note that the Gauss-Manin connection is regular singular \cite[II, \S7]{Deligne1970}, so $(\mathcal{H}', \nabla)$ indeed corresponds to $\mathbb{L}$ under Riemann-Hilbert. 

We first comment on the relationship between \autoref{Qreldef} and \autoref{hasseprincip}. Coming from \autoref{hasseprincip}, one has the following notion of non-trivial relations.

\begin{defn}
\label{Gfuncnontriv}
Let $\mathcal{G} = \{ G_{1}, \hdots, G_{m} \}$ be a set of $G$-functions in a parameter $t$, and let $Q \in K[y_{1}, \hdots, y_{m}]$ be a homogeneous polynomial. Then $Q$ is non-trivial if there does not exist a polynomial $\widetilde{Q} \in K[t][y_{1}, \hdots, y_{m}]$, homogeneous in $y_{1}, \hdots, y_{m}$, specializing to $Q$ at some point $t = \xi \in \overline{\mathbb{Q}}$, such that
\begin{equation}
\label{trivQvanishing}
\widetilde{Q}(t, G_{1}, \hdots, G_{m}) = 0 .
\end{equation}
\end{defn}

\noindent That $Q$ is non-trivial can be checked at any place $v$ of $K$, since the equality (\ref{trivQvanishing}) holds formally if it holds $v$-analytically for some $v$. 

\vspace{0.5em}

To relate this to our notion of non-triviality in \autoref{Qreldef} consider the trivialization $\mathcal{H}'^{*} \cong S' \times \mathbb{A}^{m}$ obtained using the dual frame associated to $\omega_{1}, \hdots, \omega_{m}$, and let $(G_{1}, \hdots, G_{m}) = (\Gamma(\omega_{1}), \hdots, \Gamma(\omega_{m}))$ be constructed as in \autoref{bigsupergfuncthm}. Then the graph of $(G_{1}, \hdots, G_{m})$ is identified, using the parameter $s$, with a section $\widetilde{\alpha}^{*}_{1}$ of the analytification of $\mathcal{H}^{*}$ over the neighbourhood $\mathcal{D}_{R}$ of \autoref{bigsupergfuncthm}(i), and this section agrees over $\mathcal{D}_{R} \setminus \{ 0 \}$ with a section of $\mathbb{V}'^{*}_{\mathbb{Q}}$ (up to a $2 \pi i$ multiple). One can then analytically continue $\widetilde{\alpha}^{*}_{1}$ to a neighbourhood of a point $c \in s^{-1}(\xi)$, let $\Gamma_{1}$ be the restriction of this analytic continuation to ball $B$ around $c$, and $\alpha^{*}_{1}$ the restriction of $\widetilde{\alpha}^{*}_{1}$ to the fibre above $c$. Then if $\widetilde{Q}$ is a polynomial in $s, y_{1}, \hdots, y_{m}$ which restricts at $\xi$ to a polynomial $Q$ of $y_{1}, \hdots, y_{m}$, we may take $\widetilde{R} = V(\widetilde{Q})$ and be in the setup of \autoref{Qreldef}, where we regard $\widetilde{Q}$ (resp. $Q$) as a function on $S' \times \mathbb{A}^{m}$ (resp. $\{ \xi \} \times \mathbb{A}^{m}$). The notion of non-triviality in \autoref{Qreldef} then implies the notion in \autoref{Gfuncnontriv}.

Finally, let us note that when $\mathcal{G} = \mathcal{G}_{1} \sqcup \mathcal{G}_{2}$ with $\mathcal{G}_{1}$ associated to data $(S', \mathbb{V}'_{1}, \mathcal{H}'_{1}, \hdots)$ and $\mathcal{G}_{2}$ associated to data $(S', \mathbb{V}'_{2}, \mathcal{H}'_{2}, \hdots)$, \autoref{relisnontrivial2} allows us to construct non-trivial linear relations on $\mathcal{G}$ from relations only involving the values of functions in $\mathcal{G}_{1}$ (resp. $\mathcal{G}_{2}$). Moreover products of non-trivial relations, which are individually verified to be non-trivial by applying \autoref{relisnontrivial2}, are also non-trivial: this is ultimately a consequence of the irreducibility of the algebraic monodromy group (or the primality of its defining ideal). See \cite[pg. 140, Rem. 5.3]{zbMATH00041964} and \cite[pg. 138, Rem. 2]{zbMATH00041964} for further details.

\subsection{Choosing Parameters} 
\label{chooseparamsec}

We will need to replace $S$ with a covering due to an important subtlety that occurs in applications of the $G$-function method to bound heights on curves. To understand why, suppose that $s$ is some uniformizing parameter on our curve $S$ at $s_{0}$ and our $G$-functions are, as in \autoref{bigsupergfuncthm}, obtained from expanding periods near $s_{0}$ in terms of $s$. Then if $\xi \in \mathcal{S} \subset S(\mathbb{C})$ is a special point and $v$ is a place of $K(\xi)$ which is relevant for $s(\xi)$, it could be that $\xi$ does not lie sufficiently close to $s_{0}$; more precisely, if $\mathcal{D} \subset \mathbb{A}^{1}$ is a $v$-adic disk on which the $G$-functions are defined, $\xi$ and $s_{0}$ could lie in different components of $s^{-1}(\mathcal{D})$. What we need instead is to be in the situation where every component of $s^{-1}(\mathcal{D})$ contains an appropriate degeneration point. 

To achieve this we use the following, proven in \cite[Lemma 5.1]{daw2022zilber}:

\begin{lem}
\label{daworrcovlem}
Let $C'$ be a smooth irreducible projective curve over a characteristic zero field $K$ with $s_{0} \in C'(K)$ a point. Then there exists a finite extension $L/K$, a smooth irreducible projective curve $C$ over $L$, a finite map $c : C \to C'_{L}$, and a rational function $s$ on $C$ such that
\begin{itemize}
\item[-] every zero $d \in C(\overline{L})$ of $s$ is simple;
\item[-] every zero $d \in C(\overline{L})$ of $s$ maps under $c$ to $s_{0}$;
\item[-] $s : C \to \mathbb{P}^{1}_{L}$ is a finite Galois covering (not necessarily \'etale).
\end{itemize}
\end{lem}

\begin{lem}
\label{eachcomphaszero}
In the setting of \autoref{daworrcovlem}, for any finite extension $F$ of $L$, any place $v$ of $L$, and any $R > 0$, each connected component $W$ of the analytic neighbourhood $|s|_{v} < R$ in the analytic space associated to $C_{F_{v}}$ (complex or rigid) contains a zero of $s$.
\end{lem}

\begin{proof}
Let $\mathcal{D}_{R} \subset \mathbb{A}^1 = \Spec F[x]$ be the analytic neighbourhood defined by $|x| < R$, which we take to be either complex analytic or rigid, as appropriate. Because the map $s$ is finite, so is the base-change map $s^{-1}(\mathcal{D}_{R}) \to \mathcal{D}_{R}$ (see the discussion following \cite[\S6.3 Def. 8]{zbMATH06255263} for the rigid case). Then $W \hookrightarrow s^{-1}(\mathcal{D}_{R})$ is a closed embedding, so the composition $W \to \mathcal{D}_{R}$ is again finite (in the rigid case, by loc. cit.). The image of $W$ is then a closed analytic subspace of $\mathcal{D}_{R}$ of positive dimension (see the discussion following \cite[\S6.3 Lem. 10]{zbMATH06255263} for the rigid case), hence necessarily equal to $\mathcal{D}_{R}$. It follows that the image of $W$ contains $0$, i.e., $W$ contains a zero of $s$.
\end{proof}

\subsection{Additional Setup}
\label{addsetupsec}

We now adjust the setup of \S\ref{notconvtpt1} to take into account the revised choice of parameters, and explain how we construct the set of $G$-functions we work with.

\subsubsection{Choosing the $G$-function parameter}
\label{choosingtheparamsubsubsec}

We make the following adjustments to our setup in \S\ref{notconvtpt1}:

\begin{itemize}
\item[(1)] we have a finite Galois covering $s : S \to \mathbb{P}^{1}$ with simple zeros, which we label $\{ d_{1}, \hdots, d_{k} \}$;
\item[(2)] the family $f : X \to S$ has been obtained as the pullback of a family $X_{0} \to S_{0}$ along a finite map $S \to S_{0}$ which sends $\{ d_{1}, \hdots, d_{k} \}$ to a single point $s_{0}$;
\item[(3)] for any extension $L$ of $K$, any place $v$ of $L$, and any $R > 0$, each component of the analytic neighbourhood $|s|_{v} < R$ in the associated analytic space $S_{L, v}$ (complex or rigid) contains a zero of $s$; and
\item[(4)] we no longer assume that $S \setminus S'$ is an isolated point.
\end{itemize}
To explain the reduction, complete $X \to S$ to a projective family $\overline{X} \to \overline{S}$, and apply Mumford's semistable reduction theorem \cite[Ch II]{zbMATH03425769} to reduce, after a finite base-change $\overline{S}' \to \overline{S}$, to the case where $\overline{X} \to \overline{S}$ is semistable: the divisor $E \subset \overline{X}$ consisting of the singular fibres is a simple normal crossing divisor. Replace $S$ with $\overline{S}$ and $X$ with $\overline{X}$. Then use \autoref{daworrcovlem} to obtain a covering $c : C \to S_{L}$. Then set $S_{0} = S$, $X_{0} = X$, and pull back all the data to $C$; afterwards replace $K$ with $L$, $S$ with $C$, and $S'$ with $c^{-1}(S')$. Lastly, we apply \autoref{eachcomphaszero} to obtain (3).

\subsubsection{$G$-function construction}

We construct a new, larger set of $G$-functions $\mathcal{G}$, replacing the functions $h_{1}, \hdots, h_{m}$ in \S\ref{notconvtpt1}. With an eye toward future applications, and so that our setup here is easily reused, we will do it in the general setting where we start with $\ell$ different order-$w$ normal crossing singularities in the fibre $X_{0, s_{0}}$, where $\ell$ is arbitrary. We then construct $\mathcal{G} = \{ G_{1}, \hdots, G_{M} \}$, which we regard as elements of the formal power series ring $K[[t]]$, as a union $\mathcal{G} = \mathcal{G}_{d_{1}} \sqcup \cdots \sqcup \mathcal{G}_{d_{k}}$ of sets of $G$-functions associated to the elements of the fibre $s^{-1}(0) = \{ d_{1}, \hdots, d_{k} \}$. Fixing a point $d = d_{i} \in s^{-1}(0)$, we will further sub-divide the set $\mathcal{G}_{d} = \mathcal{G}_{d_{i}}$ as a union $\mathcal{G}_{d} = \mathcal{G}_{q_{1}} \sqcup \cdots \sqcup \mathcal{G}_{q_{\ell}}$, where $q_{1}, \hdots, q_{\ell}$ are points of the fibre $X_{d}$ where the normal crossings occur. All points involved we can assume are defined over $K$ after passing to a finite extension.

Choosing a sufficiently small affine neighbourhood $U^{q_{i}} \subset X$ containing $q_{i}$ we may find functions $z_{i1}, \hdots, z_{i\nu}$ on $U^{q_{i}}$ such that:
\begin{itemize}
\itemsep0em
\item[-] the equation $z_{i1} \cdots z_{i\mu} = 0$ defines $X_{d} \cap U^{q_{i}}$ where $\mu = w +1$;
\item[-] the point $q_{i}$ lies in the locus $z_{i1} = \cdots = z_{i\nu} = 0$;
\item[-] the function $s$ maps to $z_{i1} \cdots z_{i\mu}$; and
\item[-] the map $U^{q_{i}} \to \mathbb{A}^{\nu}_{K}$ defined by $(z_{i1}, \hdots, z_{i\nu})$ is \'etale.
\end{itemize}
Indeed, the pullback of $s$ to $U^{q_{i}}$ necessarily vanishes on $X_{d}$, hence lies in the ideal defining $X_{d} \cap U^{q_{i}}$, which is necessarily locally principal, generated near $q_{i}$ by a function $z_{i1} \cdots z_{i\mu}$ such that the locus defined by $z_{ij} = 0$ is smooth at $q_{i}$ for $1 \leq j \leq \mu$. Moreover, the differentials $dz_{i1}, \hdots, dz_{i\mu}$ are independent in a neighbourhood of $q_{i}$. After shrinking $U^{q_{i}}$ we may then extend to a local \'etale coordinate system $z_{i1}, \hdots, z_{i\nu}$ at $q_{i}$. With this setup we are now, up to removing finitely many points from $S$ so that $ds$ trivializes $\Omega^{1}_{S}$, in the situation of the setup in \S\ref{realgfuncsec}; we define $\mathcal{G}_{q_{i}} = \{ \Gamma(\omega_{1}), \hdots, \Gamma(\omega_{m}) \}$, where $\omega_{1}, \hdots, \omega_{m}$ is a fixed frame for $\mathcal{H}$ near $d$.

In general an application of \autoref{bigsupergfuncthm} results in functions in a scaled parameter $\lambda s$, where $\lambda = N^{-1}$ for $N = N(d_{i}, q_{j}) \in \mathbb{Z}$ depending on $d_{i}$ and $q_{j}$. However because, by \autoref{scalesoball}, any $N'$ with $N'$ a multiple of $N$ will suffice, we may arrange for there to be some common $N$ for all elements of $\mathcal{G}$ by taking a product (or greatest common multiple) of the individual $N$'s that arise. We then replace $s$ with $N^{-1} s$.

\vspace{0.5em}

When it comes time to verify that relations we construct on the evaluations of $\mathcal{G}$ are non-trivial, which following \S\ref{appoffuncreltoGfunc} we may do at a fixed complex place $\iota : K \hookrightarrow \mathbb{C}$, we will consider the graph of the (complex realization of the) tuple $(G_{1}, \hdots, G_{M})$ as a section of $(\tau_{1}^{*} \mathcal{H}^{*})^{\ell} \oplus \cdots \oplus (\tau^{*}_{k} \mathcal{H}^{*})^{\ell}$, where $\tau_{i} : S \to S$ is an automorphism sending $d_{1}$ to $d_{i}$, which is guaranteed to exist by \S\ref{choosingtheparamsubsubsec}(1). (The automorphisms are added so that the functions in $\mathcal{G}$ may be lifted to a common neighbourhood in $S^{\textrm{an}}$.)

\subsection{Proof of \ref{onlycontrolinfplaces}}
\label{fstproofsec6}

Choose $\xi \in \mathcal{S}$ and a finite place $v$ of $K(\xi)$ which is relevant for $s(\xi)$. Let $R$ be the minimum $v$-adic convergence radius of the functions in $\mathcal{G}$. Then in particular $\xi$ lies in some component $\mathcal{D}_{R}$ of the neighbourhood of $\ad{S}_{K(\xi)_{v}}$ defined by $|s|_{v} < R$. By \S\ref{choosingtheparamsubsubsec}(3) above, necessarily such a component must contain a point $d = d_{i}$ in the fibre $s^{-1}(0)$. We may then identify our neighbourhood $\mathcal{D}_{R}$ with the one in the statement of \autoref{bigsupergfuncthm}, and construct linear relations on the values of the $G$-functions in the set $\mathcal{G}_{d} = \mathcal{G}_{q}$ corresponding to $d$. We label them $h_{1}, \hdots, h_{m}$, as in the statement of \autoref{onlycontrolinfplaces}. 

\subsubsection{Constructing Linear Relations}

Write $\pi : \mathbb{V}'_{\xi} \to W$ for the Hodge-theoretic projection. By assumption, this projection is induced by an algebraic cycle class, and so has cohomological realizations in both \'etale and algebraic de Rham cohomology, compatibly with the comparison isomorphisms. In particular, the image of $\pi$ corresponds to a summand $W_{\textrm{\'et}} \subset H^{w}(X_{\xi,\textrm{\'et}}, \mathbb{Q}_{p}(w))$ and a summand $W_{\textrm{dR}} \subset H^{w}(X_{\xi}, \Omega^{\bullet}_{X})$ which correspond under the $p$-adic Hodge comparison maps. 

\vspace{0.5em}

Set $\hat{\gamma}^{*}_{\xi} := \hat{\gamma}^{*}_{1}$ to be the functional from \autoref{bigsupergfuncthm}(ii). We first consider the situation where the restriction $\restr{\hat{\gamma}^{*}_{\xi}}{W_{\textrm{\'et}}}$ vanishes. From the compatibility with the $p$-adic Hodge comparison, this means that $\hat{\gamma}^{*}_{\xi, B_{\textrm{dR}}} \circ \rho^{-1}$ vanishes when restricted to $W_{\textrm{dR}}$. The subspace $W_{\textrm{dR}} \subset H^{w}(X_{\xi}, \Omega^{\bullet}_{X_{\xi}})$ is defined over an extension $L$ of $K(\xi)$ which, by \autoref{endoalgextprop}, has degree bounded independently of $\xi$. Let $\omega \in W_{\textrm{dR}}$ be a vector expressed as an $L$-linear combination $\omega = \sum_{i} a_{i} \omega_{i,\xi}$. Applying \autoref{bigsupergfuncthm}(ii) we find that $\sum_{i} a_{i} h_{i}(s(\xi)) = 0$, so we obtain a linear relation on the coordinates of $\hat{\gamma}^{*}_{\xi}$ in this case. 

\vspace{0.5em}

We now consider the case where $\restr{\hat{\gamma}^{*}_{\xi}}{W_{\textrm{\'et}}}$ is non-vanishing. As $W$ is a CM Hodge structure, its endomorphism algebra is a number field of degree $\dim W$, and has a primitive element $\varphi$. The characteristic polynomial $P$ of $\varphi$ has degree $\dim W$, as does its minimal polynomial (since $\varphi$ is a primitive element for a degree $\dim W$ field extension), and this is true regardless of the cohomological realization of $\varphi$ chosen. Applying \autoref{charpolyrelversion} and \autoref{boundbyhodgenumber}, we obtain, as a factor of the relation produced by \autoref{charpolyrelversion}, a linear algebraic relation $R$ on the dual coordinates of the restriction $\restr{\hat{\gamma}_{\xi}^{*}}{W_{\textrm{\'et}}}$ with respect to a basis of $W_{\textrm{dR}}$; note here that $W$ being a CM Hodge structure implies that it has more than one non-zero Hodge number, which combined with the characteristic polynomial of $\varphi$ having degree $\dim W$ puts us in \autoref{relprodsit}. If one expresses this basis as a linear combination of the basis vectors $\omega_{1,\xi}, \hdots, \omega_{m,\xi}$, one then obtains a linear relation on the dual coordinates of $\hat{\gamma}^{*}_{\xi}$ in the basis $\omega_{1}, \hdots, \omega_{m}$. This relation is defined over $K F L$, where $F$ is the splitting field of $P$ and $L$ is the field of definition of the de Rham realization of $\pi$. If one includes all factors of the relation given to us by the statement of \autoref{charpolyrelversion}, the resulting relation is defined over $K L$ and holds at all finite places as it depends only on the de Rham coordinates of $\varphi$.

\subsubsection{Ruling out functional relations:}

As explained in \S\ref{appoffuncreltoGfunc} we may reduce the task of verifying the non-triviality of linear relations on the values of the $G$-functions in $\mathcal{G}$ to the same question for the complex analytic realizations of such functions under some fixed embedding $\iota : K \hookrightarrow \mathbb{C}$. By symmetry, it suffices to consider the relations that involve just the evaluations of elements of $\mathcal{G}_{d_{1}}$. One then reduces as in \S\ref{appoffuncreltoGfunc} to the same problem for the graph of a complex analytic section $\widetilde{\gamma}^{*}_{1} \oplus \cdots \oplus \widetilde{\gamma}^{*}_{k}$ of $\tau_{1}^{*} \mathbb{V}'^{*} \oplus \cdots \oplus \tau^{*}_{k} \mathbb{V}'^{*}$ on some ball $B \subset S'(\mathbb{C})$, and where non-triviality is understood in the sense of \autoref{Qreldef}; here $\widetilde{\gamma}^{*}_{i}$ is obtained by analytically continuing the graph of the functions in $\mathcal{G}_{d_{i}}$ regarded as functions on $S'$ near $d_{1}$ using the parameter $s \circ \tau_{i}$. Then \autoref{relisnontrivial2} applied with $\mathbb{L} = \tau^{*}_{1} \mathbb{V}'_{\mathbb{C}}$ and $\mathbb{W} = \tau_{2}^{*} \mathbb{V}'^{*}_{\mathbb{C}} \oplus \cdots \oplus \tau^{*}_{k} \mathbb{V}'^{*}_{\mathbb{C}}$, and taking $E = \mathbb{Q}$ for the field in the statement of \autoref{relisnontrivial}, verifies that the relations are non-trivial (note the $\mathbb{Q}$-simplicity assumption on $\mathbb{V}'$ verifies the corresponding condition in \autoref{relisnontrivial}).

\subsubsection{Conclusion:} We have produced, by the arguments in the vanishing and non-vanishing cases, non-trivial linear relations $Q_{\textrm{nv},1}, \hdots Q_{\textrm{nv},j}$ and $Q_{\textrm{van}}$ on the value at $\xi$ of the elements of $\mathcal{G}$ whenever $\xi$ is sufficiently close to $d$. These relations were defined over a field $K L F$. Repeating this for all $d \in \{ d_{1}, \hdots, d_{k} \}$, we then obtain non-trivial relations on all evaluations at $s(\xi)$ of the functions in $\mathcal{G}$ at all relevant finite places. To obtain the conclusion we just have to argue that the relation obtained as the product of all Galois conjugates of $Q_{\textrm{nv},1}, \hdots Q_{\textrm{nv},j}$ and $Q_{\textrm{van}}$ over $K(\xi)$ has degree bounded independently of $\xi$. This is a matter of bounding $[L : K(\xi)]$ and $[F : \mathbb{Q}]$ independently of $\xi$. But the field $F$ is a splitting field of a polynomial of bounded degree, and the degree of $[L : K(\xi)]$ is bounded by \autoref{endoalgextprop}.

\subsection{Proof of \ref{sndsinggivesheights}}
\label{sndproofsec6}

The start of the argument proceeds exactly as in the proof of \autoref{onlycontrolinfplaces} to establish non-trivial relations exist at relevant finite places.

We consider a point $\xi \in \mathcal{S}_{\textrm{split}}$ and a relevant infinite place $v$ for $s(\xi)$. Let $R$ be the minimum of the convergence radii of the elements in $\mathcal{G}$. As before, using property in \S\ref{choosingtheparamsubsubsec}(3), we can find a component $\mathcal{D}_{R}$ of the complex analytic neighbourhood defined by $|s|_{v} < R$ which contains both $\xi$ and a degeneration point $d$. We write $(h_{1}, \hdots, h_{m})$ and $(h'_{1}, \hdots, h'_{m})$ for the $G$-functions associated to the two normal crossing singularities $q$ and $q'$ in the fibre $X_{d}$, which we regard as complex analytic functions.

\subsubsection{Constructing Linear Relations} We now construct a $K(\xi)$-algebraic relation on the values 
\begin{equation}
\label{hhprimevalues}
h_{1}(s(\xi)), \hdots, h_{m}(s(\xi)), h'_{1}(s(\xi)), \hdots, h'_{m}(s(\xi))  
\end{equation}
Our approach is the same as the one given in \cite[Lem. 8.3(1)]{papas2022height}. We denote by $\gamma^{*}_{1}$ and $\gamma'^{*}_{1}$ the functionals given to us by \autoref{bigsupergfuncthm}(i) corresponding to the point $\xi$. We write $L$ for the field of definition of the de Rham realization $W_{\textrm{dR}}$ of $W$, and we write $U$ for a complementary Hodge summand such that $\mathbb{V}'_{\xi} = W \oplus U$. Let $E$ be the endomorphism algebra of $W$, which is a field of dimension $\dim_{\mathbb{Q}} W$, and let $F$ be its Galois closure. 

Now $W \otimes_{\mathbb{Q}} F$ decomposes as a direct sum $\bigoplus_{\sigma : E \hookrightarrow \mathbb{C}} W_{\sigma}$ of one-dimensional weight spaces for the action of $E$. Fix an embedding $\sigma_{1} : E \hookrightarrow \mathbb{C}$, and write $E_{1}$ for its image. Then the subspace $V_{F} := U_{F} \oplus \bigoplus_{\sigma \neq \sigma_{1}} W_{\sigma}$ is stable under the action of $\textrm{Aut}(F/E_{1})$, and therefore admits a model $V \subset \mathbb{V}'_{E_{1}, \xi}$. 

Taking duals we have $\mathbb{V}'^{*}_{E_{1}, \xi} = V^{*} \oplus W^{*}_{\sigma_{1}}$, and we may choose a non-zero element $\beta^{*}_{1} \in V^{*} \cap [ \textrm{span} \{ \gamma^{*}_{1}, \gamma'^{*}_{1} \} ]_{E_{1}}$, which exists because $\gamma^{*}_{1}$ and $\gamma'^{*}_{1}$ are assumed linearly independent. Choose a second element $\beta'^{*}_{1} \in [ \textrm{span} \{ \gamma^{*}_{1}, \gamma'^{*}_{1} \} ]_{E_{1}}$ linearly independent from $\beta^{*}_{1}$.

To construct a linear relation on $\beta^{*}_{1}$ in de Rham coordinates, take an element $\omega \in W_{\sigma_{1}}$ which is defined over $LF$. Then $\beta^{*}_{1}(\omega) = 0$, so writing $\omega = \sum_{i} a_{i} \omega_{i,\xi}$ for $a_{i} \in L F$ we obtain a relation $\sum_{i} a_{i} \beta^{*}_{1}(\omega_{i,\xi}) = 0$ defined over $LF$. Call this relation $Q$. Expressing $\beta^{*}_{1}$ as a $E_{1}$-linear combination of $\gamma^{*}_{1}$ and $\gamma'^{*}_{1}$, this can be interpreted as a relation on the de Rham coordinates of $\gamma^{*}_{1}$ and $\gamma'^{*}_{1}$ which is defined over $L F$.

\subsubsection{Ruling out functional relations:} Once again we follow \S\ref{appoffuncreltoGfunc} in reducing the verification of non-triviality of linear relations on the values of the $G$-functions in $\mathcal{G}$ to the same question for the complex analytic realizations of such functions under some fixed embedding $\iota : K \hookrightarrow \mathbb{C}$. This time, one is reduced to the same problem for the graph of a complex analytic section 
\[ (\widetilde{\gamma}_{1}^{*} \oplus \widetilde{\gamma}'^{*}_{1}) \oplus \cdots \oplus (\widetilde{\gamma}_{k}^{*} \oplus \widetilde{\gamma}'^{*}_{k}) \hspace{2em} \textrm{ of } \hspace{2em} (\tau^{*}_{1} \mathbb{V}'^{*})^{2} \oplus \cdots \oplus (\tau^{*}_{k} \mathbb{V}'^{*})^{2} , \]
where $\widetilde{\gamma}^{*}_{i}$ (resp. $\widetilde{\gamma}'^{*}_{i}$) is the analytic continuation of $\gamma^{*}_{i}$ (resp. $\gamma'^{*}_{i}$) to some analytic ball $B \subset S'(\mathbb{C})$ containing $\xi$.

Using the functionals $\beta^{*}_{1}$ and $\beta'^{*}_{1}$ constructed above, we may then reduce to the same problem for the graph of a complex analytic section
\[ (\widetilde{\beta}_{1}^{*} \oplus \widetilde{\beta}'^{*}_{1}) \oplus (\widetilde{\gamma}_{2}^{*} \oplus \widetilde{\gamma}'^{*}_{2}) \oplus \cdots \oplus (\widetilde{\gamma}_{k}^{*} \oplus \widetilde{\gamma}'^{*}_{k}) \hspace{1em} \textrm{ of } \hspace{1em} (\tau^{*}_{1} \mathbb{V}'^{*}_{E_{1}})^{2} \oplus \cdots \oplus (\tau^{*}_{k} \mathbb{V}'^{*}_{E_{1}})^{2} , \]
where $\widetilde{\beta}^{*}_{1}$ (resp. $\widetilde{\beta}'^{*}_{1}$) extends $\beta^{*}_{1}$ (resp. $\beta'^{*}_{1}$). Recall that the relation $Q$ only involves $\beta^{*}_{1}$. We then apply \autoref{relisnontrivial2} with $T = S'$,
\begin{align*}
\mathbb{L} &= \tau^{*}_{1} \mathbb{V}'_{\mathbb{C}} & \alpha^{*}_{1} &= \beta^{*}_{1} \\ 
\mathbb{W} &= \tau^{*}_{1} \mathbb{V}'_{\mathbb{C}} \oplus (\tau^{*}_{2} \mathbb{V}'^{*}_{\mathbb{C}})^{2} \oplus \cdots \oplus (\tau^{*}_{k} \mathbb{V}'^{*}_{\mathbb{C}})^{2} & \zeta^{*}_{1} &= \widetilde{\gamma}'^{*}_{1} \oplus \cdots \oplus (\widetilde{\gamma}_{k}^{*} \oplus \widetilde{\gamma}'^{*}_{k}) ,
\end{align*}
and $E = E_{1}$. As $\mathbb{V}'$ was assumed $E$-simple, we conclude that $Q$ is non-trivial.

\subsubsection{Conclusion} The above reasoning shows that at each relevant place $v$ for $s(\xi)$ we may find a linear $K L F$-algebraic relation on the values of $\mathcal{G}$ at $s(\xi)$ which is non-trivial, and such that the degrees $[L : K(\xi)]$ and $[F : \mathbb{Q}]$ are bounded independently of $\xi$; note that to obtain the bound on $[L : K(\xi)]$ we are using \autoref{endoalgextprop}. Taking the product of conjugates of these relations we obtain $K(\xi)$-algebraic relations $Q_{v}$ of degree bounded by a constant $\kappa$ independent of $\xi$ and the place $v$; as we saw in the proof of \autoref{onlycontrolinfplaces} the relations $Q_{v}$ may be assumed to be the same for all finite places. The product $Q = \prod_{v} Q_{v}$ taken over the relevant infinite places and a relevant finite place (if it exists) is then a $K(\xi)$-algebraic relation of degree at most $\kappa \,  ([K(\xi) : \mathbb{Q}] + 1)$ which holds at all relevant places for $s(\xi)$. Using the Weil height $\theta$ induced by the parameter $s$ one concludes from \autoref{hasseprincip} that for all $\xi \in \mathcal{S}_{\textrm{split}}$
\[ \theta(\xi) \leq \kappa' \, [K(\xi) : \mathbb{Q}]^{a} , \]
for some constants $\kappa', a \in \mathbb{R}_{> 0}$.

\begin{rem}
The reader may wonder how Bombieri's Theorem \autoref{hasseprincip} applies for those $\xi$ for which there are no relevant places. In this case we are always in branch (2)(i) of \autoref{hasseprincip}, and so the condition in (2)(ii) in \autoref{hasseprincip} is never relevant. In particular, one can choose for $Q$ any linear relation for which \autoref{hasseprincip}(1) holds, and obtain an absolute bound on such points. (See also the treatment in \cite[\S12]{papas2023unlikely}.)
\end{rem}

\section{Pila-Zannier for General Intersections}
\label{pilazansec7}

In this section we drop all ``primed'' superscripts; i.e., write $S$ instead of $S'$, etc., as compactifications and degeneration points will not be relevant. A general overview of the ideas in this section is given in \S\ref{pilazannatypintro}.

We make an additional comment about our presentation. The Zilber-Pink conjecture has a ``geometric'' part and an ``arithmetic'' part. With reference to the formulation appearing in \cite{BKU}, which works in the setting of a general polarizable integral variation of Hodge structure $\mathbb{V}$ on a complex algebraic variety $S$, the difference between the two is whether the ``atypical'' Hodge loci one considers in the base of the variation of Hodge structures are positive or zero dimensional (in the sense of period dimension, see \cite[Def. 1.2]{BKU}). It is known as a consequence of \cite[Theorem 3.1]{BKU} that geometric Zilber-Pink statements --- statements about non-Zarski density of positive-dimensional atypical Hodge loci --- can be be proven without any arithmetic input, using only general Hodge-theoretic facts. (There is an exception in that one cannot resolve Zilber-Pink-type statements for positive-dimensional Hodge loci coming from splittings of the generic Hodge datum associated to $S$, but this can be reinterpreted as a failure to resolve the conjecture for zero-dimensional Hodge loci for an auxiliary variation constructed from this datum.) We therefore focus exclusively on the zero-dimensional, i.e. arithmetic, part of Zilber-Pink here.

Throughout, ``definable'' is taken to mean with respect to the o-minimal structure $\mathbb{R}_{\textrm{an,exp}}$.

\subsection{Pila-Zannier Preliminaries}

\subsubsection{Hodge-theoretic Setup}
\label{Hodgetheosetupsec}

Throughout \S\ref{pilazansec7} we will assume we have a smooth projective $K$-algebraic family $f : X \to S$, with $K \subset \mathbb{C}$ a number field, over a smooth base $S$ with geometrically connected fibres; we note in particular we do not assume $S$ is a curve. Subvarieties of $S$ will be understood to mean subvarieties of $S_{\mathbb{C}}$. We set $\mathbb{V} = R^{w} f^{\textrm{an}}_{*} \mathbb{Z} / \textrm{tor}$.

Let $\pi : \widetilde{S} \to S^{\textrm{an}}$ be the universal covering, and let us fix a point $\widetilde{c} \in \widetilde{S}$ with image $c \in S(\mathbb{C})$. We choose a frame $\alpha_{1}, \hdots, \alpha_{m}$ of $\widetilde{\mathbb{V}}_{\widetilde{c}}$, where $\widetilde{\mathbb{V}} = \pi^{*}\mathbb{V}$. We then have an induced identification $\widetilde{\mathbb{V}}_{\widetilde{c}} \cong \mathbb{Z}^{m}$, and using parallel transport we obtain a natural height function $\widetilde{\theta}$ on all tensor spaces associated to any fibre of $\widetilde{\mathbb{V}}$ induced by the standard Weil height on the affine space associated to $\mathbb{Z}^{m}$. Fixing a polarization $Q : \mathbb{Z}^{m} \otimes \mathbb{Z}^{m} \to \mathbb{Z}$, we then consider the space $D$ of all polarized Hodge structures on $\mathbb{Z}^{m}$ with the same Hodge numbers $\{ h^{p,q} \}_{p+q=n}$ as the variation $\mathbb{V}$. Applying \cite[I.A Def.(iii)]{GGK} and \cite[II.A]{GGK}, this set has the following description: $D$ is the set of morphisms $h : \mathbb{S} \to \textrm{GAut}(\mathbb{Z}^{m}, Q)_{\mathbb{R}}$, where $\mathbb{S} = \textrm{Res}_{\mathbb{C}/\mathbb{R}} \mathbb{G}_{m,\mathbb{C}}$ is the Deligne torus, such that the induced weight space decomposition on $\mathbb{C}^{m}$ as in \cite[I.A Def.(iii)]{GGK} defines a pure polarized Hodge structure with Hodge numbers $\{ h^{p,q} \}_{p+q=n}$. We write $\ch{D}$ for the compact dual of $D$, defined as in \cite[Ch. II]{GGK}. The Mumford-Tate group $M_{h}$ of a Hodge structure $h \in D$ is the $\mathbb{Q}$-Zariski closure of $h(\mathbb{S})$. 

For $h \in D$, we will also have use for the \emph{special Mumford-Tate} subgroup $_{\textrm{sp}} M_{h} \subset M_{h}$. To define this, one considers the natural $\mathbb{R}$-algebraic subgroup $\mathbb{U} \subset \mathbb{S}$ defined on \cite[pg. 29]{GGK}, and takes $_{\textrm{sp}} M_{h}$ to be the $\mathbb{Q}$-Zariski closure of $h(\mathbb{U})$. As explained on \cite[pg. 33]{GGK}, assuming $h$ is not pure of weight zero, the groups $_{\textrm{sp}} M_{h}$ and $M_{h}$ differ by a multiplicative almost-direct $\mathbb{Q}$-factor $\mathbb{G}_{m}$ obtained as the image of the restriction of $h$ to $\mathbb{G}_{m} \subset \mathbb{S}$ (c.f. \cite[pg. 31, Def.]{GGK}). We will assume that our Hodge structures on $\mathbb{Q}^{m}$ are pure of fixed weight $w$, which implies that the $\mathbb{G}_{m}$ factors are all identified with the central scalar subgroup $Z \subset \GL_{m}$ in the case $w \neq 0$ (c.f. \cite[pg. 31, Def.]{GGK}) and that $M_{h} ={_{\textrm{sp}}}M_{h}$ in the case $w = 0$. 

\begin{lem}
\label{spindepofh}
With the above setup, we have ${_{\textrm{sp}}}M_{h} = (M_{h} \cap \textrm{SL}_{m})^{\circ}$ for each $h \in D$, where $(-)^{\circ}$ denotes passage to the identity component. In particular, ${_{\textrm{sp}}}M_{h}$ depends only on $M_{h}$ and not $h$. 
\end{lem}

\begin{proof}
From \cite[Rem. 2.1.12]{zbMATH06829467} one has that ${_{\textrm{sp}}}M_{h} \subset \textrm{SL}_{m}$, that ${_{\textrm{sp}}}M_{h} = M_{h}$ when the weight $w = 0$, and that $M_{h} = Z \cdot {_{\textrm{sp}}}M_{h}$ when $w \neq 0$. The algebraic group ${_{\textrm{sp}}} M_{h}$ is (geometrically) connected, since it is the $\mathbb{Q}$-Zariski closure of the image of the connected $\mathbb{R}$-algebraic group $\mathbb{U}$, and the identity component of a $\mathbb{Q}$-algebraic group is always defined over $\mathbb{Q}$. When $w \neq 0$ intersecting $M_{h}$ with $\textrm{SL}_{m}$ drops the dimension by $1$, and so the identity component of the resulting intersection, which is a (geometrically) connected algebraic group containing the (geometrically) connected algebraic group ${_{\textrm{sp}}} M_{h}$, must be equal to ${_{\textrm{sp}}} M_{h}$.
\end{proof}

We consider the map $\widetilde{\varphi} : \widetilde{S} \to D$ given by sending $\widetilde{s}$ to the parallel translate of its Hodge flag $F^{\bullet}_{\widetilde{s}}$ to the fibre $\widetilde{\mathbb{V}}_{\widetilde{c}} \cong \mathbb{Z}^{m}$. We assume that $\widetilde{\varphi}$ has discrete fibres. This is the same as saying that, if we set $\Gamma = \textrm{Aut}(\mathbb{Z}^{m}, Q)(\mathbb{Z})$, the associated period map $\varphi : S^{\textrm{an}} \to \Gamma \backslash D$ is quasi-finite. 

In the next section, we review the notion of a \emph{Siegel (sub)set} of the space $D$. Applying \cite[Thm. 1.5]{defpermap} and \cite[\S1.5]{zbMATH07720398}, we may find a finite union $\mathcal{F} = \bigsqcup_{i=1}^{n} \mathfrak{G}_{i}$ of Siegel subsets of $D$ such that $\mathcal{I} := \widetilde{\varphi}(\widetilde{S}) \cap \mathcal{F}$ and $\widetilde{\varphi}(\widetilde{S})$ have the same image in $\Gamma \backslash D$, and the sets $\{ \mathfrak{G}_{i} \}_{i=1}^{n}$ and $\mathcal{I}$ are definable subsets of $D$. (See also the discussion in \cite[\S4.1, \S4.2]{defpermap} for an explanation of how \cite[Thm. 1.5]{defpermap} implies a global statement about the lift $\widetilde{\varphi}$ of $\varphi$.)

\subsubsection{Siegel Set Background}

Following \cite[\S2]{zbMATH06880895}, a Siegel set for a reductive $\mathbb{Q}$-group $G$ is a product $\Omega A_{t} K$ associated to a pair $(P, K)$ consisting of a minimal $\mathbb{Q}$-parabolic subgroup $P \subset G$ and a maximal compact subgroup $K \subset G(\mathbb{R})$. Associated to $P$ and $K$, there is a distinguished real torus $S \subset P_{\mathbb{R}}$ characterized by the conditions in \cite[Lem 2.1]{zbMATH06880895}. Then to construct the product $\Omega A_{t} K$ one takes:
\begin{itemize}
\item[(i)] $\Delta$ to be the set of simple roots for $G$ with respect to $S$, with the order from $P$;
\item[(ii)] for some real $t > 0$ the set
\[ A_{t} := \{ \alpha \in S(\mathbb{R})^{+} : \chi(\alpha) \geq t\textrm{ for all } \chi \in \Delta \} ; \]
\item[(iii)] and $\Omega \subset P(\mathbb{R})$ to be a compact set (c.f. \cite[\S2.C]{zbMATH06880895}).
\end{itemize}

\noindent Set $G = \textrm{GAut}(\mathbb{Z}^{m}, Q)$, let $Z \subset G$ be its center, and let $h : \mathbb{S} \to G_{\mathbb{R}}$ be a map which defines a point of $D$. 

\begin{notn}
For a polarized vector space $(V, Q)$ we will write $\textrm{SAut}(V,Q)$ for the subgroup of $\textrm{GAut}(V,Q)$ consisting of elements with determinant $1$.
\end{notn}

\begin{lem}
\label{maxcomplem}
Let $H_{h} \subset G(\mathbb{R})$ be the centralizer of $h$ (c.f. \cite[Ch. II]{GGK}). Then there is a unique maximal compact subgroup $K_{h} \subset \textrm{SAut}(\mathbb{Q}^{m},Q)(\mathbb{R})$ such that $H_{h} \subset Z(\mathbb{R}) K_{h}$. 
\end{lem}

\begin{proof}
The group in question is constructed by Schmid in \cite[(8.4)]{Schmid1973}. Its uniqueness is claimed on \cite[pg. 75]{GGK} and can be verified on the level of Lie algebras using the Hodge-Riemann bilinear relations for the adjoint Hodge structure on the Lie algebra of $G_{\mathbb{R}}$ induced by $h$.
\end{proof}

\begin{defn} A Siegel set $\mathfrak{O}$ for $D$ is an orbit $\mathfrak{G} \cdot h$ where $\mathfrak{G} = \Omega A_{t} K_{h}$ is a Siegel set for $G$ associated to a pair $(P, K_{h})$.
\end{defn} 

\noindent By \cite[II.A]{GGK} the action of $G(\mathbb{R})$ induces an identification $D = G(\mathbb{R})/H_{h}$, which implies:

\begin{lem}
\label{invimageSiegel}
For a Siegel set $\mathfrak{O} \subset D$, the preimage of $\mathfrak{O}$ in $G(\mathbb{R})$ under the orbit map $g \mapsto g \cdot h$ is exactly $\mathfrak{G}$. \qed
\end{lem}

We will want to use the natural embedding $G \subset \GL_{m}$ to embed a Siegel set $\mathfrak{G} \subset G(\mathbb{R})$ into a Siegel set $\mathfrak{S}$ for $\GL_{m}(\mathbb{R})$, and we will want to choose the compact group $K$ associated to $\mathfrak{S}$ carefully. We let $C_{h}$ denote the Weil operator \cite[Ch. I]{GGK} associated to the Hodge structure $h$, and $\psi_{h}(-,-) = Q(C_{h}(-),-)$ the associated positive-definite form.

\begin{prop}
\label{bettersiegelembedding}
For any Siegel set $\mathfrak{G} = \Omega A_{t} K_{h}$ associated to the Hodge structure $h$, there exists:
\begin{itemize}
\item[-] a Siegel set $\mathfrak{S} = \Lambda B_{t} K$ for $\GL_{m}$ with $K$ the orthogonal group of the form $\psi_{h}$; and
\item[-] a finite set $C \subset \GL_{m}(\mathbb{Q})$ such that $\mathfrak{G} \subset C \cdot \mathfrak{S}$.
\end{itemize}
\end{prop}

\begin{proof}
One checks directly that $K_{h} = K \cap G(\mathbb{R})$ (c.f. \cite[\S3.8]{bakker2023finiteness}). Applying \cite[Thm. 4.1]{orrarXivv6} and the discussion on pg. 16, it suffices to verify condition (2) in \cite{orrarXivv6} appearing below \cite[Lem 4.3]{orrarXivv6}. The Cartan involution associated to $\psi_{h}$ is given by $X \mapsto C_{h} (X^{\dagger})^{-1} C^{-1}_{h}$ (c.f. \cite[\S3.8]{bakker2023finiteness} in the even weight case), where $\dagger$ is the adjunction coming from $Q$. On the other hand by \cite[\S2]{milne2023classification} (c.f \cite[\S3.4]{bakker2023finiteness}) the restricted map $g \mapsto C_{h} g C_{h}^{-1}$ on $G(\mathbb{R})$ is the Cartan involution associated to $K_{h}$. The reasoning in \cite[Lem 4.A]{orrarXivv6} then verifies (2).
\end{proof}

\subsection{Reduction to Height Bounds on Tensors}

In what follows we will freely use the notions of \emph{special} and \emph{weakly special} subvariety of $S$ for the variation $\mathbb{V}$, defined for instance as in \cite[Def. 4.4]{BKU}. A special point is a zero dimensional special subvariety. We will also use:

\begin{defn}
\label{sppointIdef}
A special point of $\mathcal{I}$ is a point $h \in \mathcal{I}$ whose image in $\Gamma \backslash D$ agrees with the image of a special point of $S$.
\end{defn}

Let $Y \subset S$ be a geometrically irreducible complex subvariety.

\begin{defn}
The Mumford-Tate group $\mathbf{G}_{Y}$ is the Mumford-Tate group of $\mathbb{V}_{s}$ above a generic point $s \in Y$ (i.e., a point outside the Hodge locus of $\restr{\mathbb{V}}{Y}$). 
\end{defn}

For each point $y \in Y(\mathbb{C})$, the group $\mathbf{G}_{Y}$ has a unique realization $\mathbf{G}_{Y,y} \subset \GL(\mathbb{V}_{y})$, and the corresponding map $h_{y} : \mathbb{S} \to \GL(\mathbb{V}_{y})_{\mathbb{R}}$ factors through $\mathbf{G}_{Y,y,\mathbb{R}}$. If one uses the flat structure of $\restr{\mathbb{V}}{Y}$ to identify the groups $\mathbf{G}_{Y,y}$, the resulting family of maps $h_{y} : \mathbb{S} \to \mathbf{G}_{Y,\mathbb{R}}$ lies in a single $\mathbf{G}_{Y}(\mathbb{R})$-orbit $D_{Y}$ of maps $\mathbb{S} \to \mathbf{G}_{Y,\mathbb{R}}$, where the action is by conjugacy (see \cite[III.A Step 3]{GGK}). 

\begin{defn}
The data $(\mathbf{G}_{Y}, D_{Y})$ is the \emph{Hodge datum} associated to $(Y,  \restr{\mathbb{V}}{Y})$. 
\end{defn}

Fixing a point $y \in Y$, one obtains realizations $(\mathbf{G}_{Y,y}, D_{Y,y})$ of the pair $(\mathbf{G}_{Y}, D_{Y})$ associated to the fibre $\mathbb{V}_{y}$, and the flat structure of $\restr{\mathbb{V}}{Y}$ preserves these realizations. In the case of $Y = S$, our trivialization $\widetilde{\mathbb{V}} \cong \mathbb{Z}^{m}$ fixed in \S\ref{Hodgetheosetupsec} above gives a distinguished embedding $(\mathbf{G}_{S}, D_{S}) \hookrightarrow (G\textrm{Aut}(\mathbb{Q}^{m}, Q), D)$ and we will identify $(\mathbf{G}_{S}, D_{S})$ with its image. 

\begin{defn}
\label{MTsubgpdef}
By a \emph{Mumford-Tate subgroup} of $\mathbf{G}_{S}$, we will always mean the Mumford-Tate group $M = M_{h}$ of some Hodge structure $h \in D_{S} \subset D$. 
\end{defn}

\begin{notn}
In this section we write $\mathbf{H}_{S} = \mathbf{H}_{\mathbb{V}}$ for the algebraic monodromy group of $(S, \mathbb{V})$.
\end{notn}

\begin{defn}
\label{defaspvar}
Given a special point $\{ y \} \subset S$, we say that a pair $(M, h)$, with $M$ a Mumford-Tate subgroup of $\mathbf{G}_{S}$ and $h \in D_{S}$, \emph{defines} $\{ y \}$ if
\begin{itemize}
\item[-] there exists a point $h' \in D_{S}$ such that $M = M_{h'}$, and $h \in D_{M} := M(\mathbb{R}) \cdot h'$;
\item[-] there exists a lift $\widetilde{y} \in \widetilde{S}$ of $y$ with image $h = \widetilde{\varphi}(\widetilde{y})$; and
\item[-] $\{ y \}$ is obtained as an irreducible component of $\pi(\widetilde{\varphi}^{-1}(D_{M}))$.
\end{itemize}
In this situation we also say that $M$ defines the special point $\{ h \} \in \mathcal{I}$. 
\end{defn}

\begin{defn}
\label{defaspvar2}
Given a special point $\{ y \} \subset S$ and a subset 
\[ T \subset \bigoplus_{a, b \geq 0} (\mathbb{V}_{\mathbb{Q},y})^{\otimes a} \otimes (\mathbb{V}_{\mathbb{Q},y}^{*})^{\otimes b} \]
consisting of Hodge tensors for the Hodge structure on $\mathbb{V}_{y}$, we say that $y$ is \emph{defined by }$T$ if there exists a pair $(M,h)$ as in \autoref{defaspvar} which defines $\{ y \}$, and such that $M$ is identified with the stabilizer in $\mathbf{G}_{S}$ of $T$ after applying the isomorphism $\mathbb{V}_{y} \cong \mathbb{V}_{\widetilde{y}} \cong \mathbb{Z}^{m}$ induced by the $\widetilde{y}$ in \autoref{defaspvar}. 
\end{defn}

\begin{lem}
\label{finmanyconjclasses}
Let $G \in \{ \mathbf{G}_{S}, G\textrm{Aut}(\mathbb{Q}^{m}, Q) \}$. Then there are finitely many Mumford-Tate subgroups $M \subset G$ up to $G(\mathbb{R})$-conjugacy (and hence $G(\mathbb{C})$-conjugacy).
\end{lem}

\begin{proof}
For $\GL_{m}(\mathbb{C})$-conjugacy this is \cite[Theorem 4.14]{voisin2010hodge}; the same proof works for $G(\mathbb{C})$-conjugacy.

 To get $G(\mathbb{R})$-conjugacy, we note that there are finitely many real semisimple subgroups of real linear algebraic groups up to real conjugacy as a consequence of \cite{richardson1967} (c.f. \cite[Cor. 0.2]{zbMATH07411859}), and also finitely many $G(\mathbb{R})$-conjugacy classes of maximal real tori in $G_{\mathbb{R}}$ \cite[Ch. III, \S4.4]{zbMATH00707215}. Arguing as in \cite[Theorem 4.14]{voisin2010hodge}, we may therefore reduce to showing there are finitely many conjugacy classes of real reductive subgroups of $G_{\mathbb{R}}$ generated by images of maps $G(\mathbb{R})$-conjugate to those in a fixed finite set $F = \{ h_{i} : \mathbb{S} \to G_{\mathbb{R}} \}_{i=1}^{m}$. We then run the argument of \cite[Prop. 4.15]{voisin2010hodge} with $G$ a real-algebraic group, complex conjugacy replaced by real conjugacy, and where the set of $1$-parameter subgroups in \cite[Prop. 4.15]{voisin2010hodge} is replaced by $F$.
\end{proof}

In what follows we write $\ch{\textrm{NL}}_{M} \subset \ch{D}$ for the locus of Hodge flags $F^{\bullet}$ with Mumford-Tate group $M_{F^{\bullet}} \subset M$ (see \cite[Ch 2., pg. 49]{GGK}). 

\begin{lem}
\label{finMTgroups}
Let $M \subset \mathbf{G}_{S}$ be a Mumford-Tate subgroup. Then the number of special points (recall \autoref{sppointIdef}) in $\mathcal{I}$ defined by $M$ is bounded by a constant $\kappa$ independent of $M$. More generally, the number of isolated points in $\mathcal{I} \cap g \ch{\textrm{NL}}_{M}$ is bounded by a constant $\kappa$ independent of $g \in \mathbf{G}_{S}(\mathbb{C})$ and $M$. 
\end{lem}

\begin{proof}
The orbit $M(\mathbb{R}) \cdot h$ is open in $\ch{\textrm{NL}}_{M}$ as a consequence of \cite[II.B]{GGK} and \cite[VI.B.1]{GGK}, so the first problem reduces to the second. Because Mumford-Tate subgroups of $\mathbf{G}_{S}$ lie in finitely many $\mathbf{G}_{S}(\mathbb{C})$-conjugacy classes by \autoref{finmanyconjclasses}, we may reduce to the case where the $\mathbf{G}_{S}(\mathbb{C})$-conjugacy class of $M$ is fixed. Considering the intersection $\mathcal{I}_{M'} := \mathcal{I} \cap \ch{\textrm{NL}}_{M'}$ for some Mumford-Tate group $M' = g M g^{-1}$, where $g \in \mathbf{G}_{S}(\mathbb{C})$, we have $\ch{\textrm{NL}}_{M'} = g \cdot \ch{\textrm{NL}}_{M}$ and $\mathcal{I}_{M'} = \mathcal{I} \cap (g \cdot \ch{\NL}_{M})$. But the number of isolated points in any fibre of the definable family $\{ \mathcal{I} \cap (g \cdot \ch{\NL}_{M}) \}_{g \in \mathbf{G}_{S}(\mathbb{C})}$ is bounded by a uniform constant $\kappa$ as a consequence of definability.
\end{proof}

For an algebraic group $H \subset \textrm{SAut}(\mathbb{Q}^{m}, Q)$, write $\textrm{Stab}(H) \subset  \bigoplus_{a \geq 1} (\mathbb{Q}^{m})^{\otimes a}$ for the set pointwise stabilized by $H$. Consider a Mumford-Tate subgroup $M \subset \mathbf{G}_{S}$ associated to a point of $D_{S}$. 

\begin{defn}
We say that a set $T \subset \bigoplus_{a \geq 1} (\mathbb{Q}^{m})^{\otimes a}$ of tensors is a set of \emph{minimal type} defining $M$, if:

\begin{itemize}
\item[(i)] $T$ admits a partition $T = T_{1} \sqcup T_{2} \sqcup \cdots \sqcup T_{\ell}$ where $T_{i} \subset (\mathbb{Q}^{m})^{\otimes i}$;
\item[(ii)] $T$ is a basis for $\textrm{Stab}(_\textrm{sp} M) \cap \bigoplus_{\ell \geq a \geq 1} (\mathbb{Q}^{m})^{\otimes a}$; and
\item[(iiii)] $\ell$ is chosen as small as possible subject to (i), (ii) and the condition that $_\textrm{sp} M$ is exactly the pointwise stabilizer of $T$ in $_\textrm{sp} \mathbf{G}_{S}$.
\end{itemize}

\noindent We additionally say that $T$ is a \emph{minimal} set defining $M$ if in addition:

\begin{itemize}
\item[(iv)] Subject to (i), (ii) and (iii), the set $T$ is chosen to have minimal total height, where the height $\widetilde{\theta}(T)$ is the maximum of the heights of all its (necessarily finitely many) elements.
\end{itemize}
\end{defn}

\begin{rem}
We note that, as a consequence of \autoref{spindepofh}, the group $_\textrm{sp} M$ depends only on $M$ and not the point of $D_{S}$ which $M$ is associated to, so it makes sense to think of $T$ as being associated to $M$ itself.
\end{rem}

Note that given a fixed algebraic subgroup $M \subset \mathbf{G}_{S}$, the coset space $\mathbf{G}_{S}/M$ has, by \cite[Thm. 4.4.1]{zbMATH01746054}, the structure of a quasi-projective algebraic variety in such a way that the natural quotient map $q : \mathbf{G}_{S} \to \mathbf{G}_{S}/M$ is algebraic.

\begin{prop}
\label{heightconstrprop}
Fix the data of 
\begin{itemize}
\item[-] a Mumford-Tate group $M \subset \mathbf{G}_{S}$, a quasi-projective embedding $\iota : (\mathbf{G}_{S}/M) \hookrightarrow \mathbb{P}^{r}$, and the induced Weil height $\widetilde{\theta}_{M} : (\mathbf{G}_{S}/M)(\overline{\mathbb{Q}}) \to \mathbb{R}$; 
\item[-] a set $\mathcal{S} \subset S(\mathbb{C})$ of special points defined by pairs $(M_{\widetilde{\xi}}, \widetilde{\varphi}(\widetilde{\xi}))$, with $M_{\widetilde{\xi}}$ conjugate by $\mathbf{G}_{S}(\mathbb{C})$ to $M$ and $\widetilde{\xi} \in \widetilde{S}$ a lift of $\xi \in \mathcal{S}$.
\end{itemize}
Then if $P_{1}, \hdots, P_{r} \in \ch{\NL}_{M}(\mathbb{C})$ are representatives of the components of $\ch{\NL}_{M}$, there exists a map $g_{(-)} M : \mathcal{S} \to (\mathbf{G}_{S}/M)(\overline{\mathbb{Q}})$, denoted $\xi \mapsto g_{\xi} M$, with the following properties:
\begin{itemize}
\item[(1)] there exists a constant $d$ such that each point in the image of $g_{(-)} M$ lies inside a number field of degree at most $d$;
\item[(2)] there exists a constant $c$ such that each fibre of $g_{(-)} M$ has cardinality at most $c$;
\item[(3)] the height $\widetilde{\theta}_{M}(g_{\xi} M)$ is bounded by a polynomial in $\widetilde{\theta}(T_{\widetilde{\xi}})$, where $T_{\widetilde{\xi}}$ is a minimal set of definition for $M_{\widetilde{\xi}}$; and 
\item[(4)] for each $\xi \in \mathcal{S}$ we have $g^{-1}_{\xi} \cdot M_{\widetilde{\xi}} \cdot g_{\xi} = M$ and that $\widetilde{\varphi}(\widetilde{\xi}) \in (g_{\xi} \cdot M(\mathbb{C}) \cdot P_{j})$ for some $j$.
\end{itemize}
\end{prop}

\begin{rem}
We note that $\ch{\NL}_{M}$ is a finite union of $M(\mathbb{C})$-orbits; see \cite[VI.B.1]{GGK}.
\end{rem}

\begin{rem}
We emphasize that $M_{\widetilde{\xi}}$ need not be the Mumford-Tate group of $\widetilde{\varphi}(\widetilde{\xi})$, and is the Mumford-Tate group of a point $h' \in D_{S}$ for which $\widetilde{\varphi}(\widetilde{\xi}) \in M_{h'}(\mathbb{R}) \cdot h'$, where $(M_{h'}, \widetilde{\varphi}(\widetilde{\xi}))$ satisfies the hypotheses of \autoref{defaspvar}.
\end{rem}

\begin{proof}
We note that, because all (special) Mumford-Tate groups are defined by the tensor invariants they stabilize and the spaces of these tensor invariants are defined over $\mathbb{Q}$, any two minimal sets defining Mumford-Tate groups conjugate to $M$ have the same number of tensors in each degree. Denoting by $d_{1}, d_{2}, \hdots, d_{\ell}$ the dimensions of the associated subspaces of tensors, we may then consider an affine parameter space 
\[ \mathcal{T} = \textrm{Aff}\left[(\mathbb{Q}^{m})\right]^{d_{1}} \times \textrm{Aff}\left[(\mathbb{Q}^{m})^{\otimes 2}\right]^{d_{2}} \times \cdots \times \textrm{Aff}\left[(\mathbb{Q}^{m})^{\otimes \ell}\right]^{d_{\ell}} , \]
where for a vector space $V$ we write $\textrm{Aff}[V]$ for the associated affine space. To any point $T \in \mathcal{T}$ we may associate a group $M_{T}$ defined as the point-wise stabilizer of $T$. We consider the $\mathbb{Q}$-algebraic family $v : \mathcal{C} \to \mathcal{T} \times \GL_{m}$ whose fibre above $(T, g)$ is $g M_{T} g^{-1}$. Applying \cite[\href{https://stacks.math.columbia.edu/tag/02FZ}{Lemma 02FZ}]{stacks-project}, \cite[\href{https://stacks.math.columbia.edu/tag/05F9}{Lemma 05F9}]{stacks-project} and \cite[\href{https://stacks.math.columbia.edu/tag/055I}{Lemma 055I}]{stacks-project} there is a constructible locus $\mathcal{D}_{1} \subset \mathcal{T} \times \GL_{m}$ above which the fibre of $v$ is a connected algebraic group of dimension $\dim (_\textrm{sp}M)$. By considering a family $v' : \mathcal{C} \to \mathcal{T} \times \GL_{m}$ whose fibre above $(T, g)$ is $[g M_{T} g^{-1} \cap\ _\textrm{sp}M]$, we also obtain a constructible algebraic locus $\mathcal{D}_{2} \subset \mathcal{T} \times \GL_{m}$ above which the fibre of $v'$ has dimension $\dim _\textrm{sp}M$. Then over the locus $\mathcal{D}_1 \cap \mathcal{D}_2$ one has $\dim [g M_{T} g^{-1} \cap\ _\textrm{sp}M] = \dim g M_{T} g^{-1} = \dim _\textrm{sp}M$ with all varieties connected, hence $g M_{T} g^{-1} =\ _\textrm{sp}M$. Let $\mathcal{T}_{0} \subset \mathcal{T}$ be the image in $\mathcal{T}$ of $\mathcal{D}_1 \cap \mathcal{D}_2$.

The locus $\mathcal{T}_{0}$ is ($\mathbb{Q}$-algebraically) constructible in $\mathcal{T}$, and on the level of complex points is the locus of $T \in \mathcal{T}(\mathbb{C})$ for which $M_{T}$ is conjugate to $_\textrm{sp} M$. Being a constructible algebraic subset of $\mathcal{T}$, it is (by definition, see \cite[\href{https://stacks.math.columbia.edu/tag/04ZC}{Section 04ZC}]{stacks-project}) a finite union of locally closed $\mathbb{Q}$-subvarieties $\mathcal{T}'_{0} \subset \mathcal{T}$. For such a subvariety, we write $y : \mathcal{G} \to \mathcal{T}'_{0}$ for the $\mathbb{Q}$-algebraic family whose fibre above $T$ consists of those $g \in \mathbf{G}_{S}(\mathbb{C})$ satisfying the property that $g^{-1} \cdot M_{T} \cdot g ={_\textrm{sp}M}$. 

Replacing $\mathcal{T}_{0}$ by $\mathcal{T}'_{0}$ we claim it suffices to show that, for any $T \in \mathcal{T}_{0}(\mathbb{Q})$, we can construct $g(T) \in y^{-1}(T)(\overline{\mathbb{Q}})$ defined over a number field of uniformly bounded degree and whose height is bounded by a polynomial in $\widetilde{\theta}(T)$. Indeed, suppose we can do this, and define $g_{\xi} M = g(T_{\widetilde{\xi}}) M$ where $T_{\widetilde{\xi}}$ is a minimal set of definition for $M_{\widetilde{\xi}}$. We may then see that (1) is true by assumption. For (2) we note that if $g(T_{\widetilde{\xi}_{1}}) M = g(T_{\widetilde{\xi}_{2}}) M$ then $M_{\widetilde{\xi}_{1}} = M_{\widetilde{\xi}_{2}}$, and the number of lifts $\widetilde{\xi}$ mapping to $\mathcal{I}$ of points in $\mathcal{S}$ for which this can occur is bounded by \autoref{finMTgroups}. Then (3) is true by the properties of heights under polynomial maps, and as $\ch{\NL}_{M_{\widetilde{\xi}}} = g_{\xi} \cdot \ch{\NL}_{M}$ statement (4) reduces to the claim that $\widetilde{\varphi}(\widetilde{\xi})$ has Mumford-Tate group contained in $M_{\widetilde{\xi}}$, which is true because $(M_{\widetilde{\xi}}, \widetilde{\varphi}(\widetilde{\xi}))$ defines $\xi$. As only finitely many possible $\mathcal{T}_{0}$ occur, we obtain the result after possibly increasing the constants finitely many times.

Continuing now with our claim about the family $y$, we may consider a stratification $\mathcal{T}_{0} = \mathcal{T}_{1} \sqcup \cdots \sqcup \mathcal{T}_{k}$ such that for each $i$ the base-change $y_{i} : \mathcal{G}_{i} \to \mathcal{T}_{i}$ along $\mathcal{T}_{i} \hookrightarrow \mathcal{T}_0$ of $y$ is flat and each $\mathcal{T}_{i}$ is smooth; base-changing to one of these strata, we may assume that $y : \mathcal{G} \to \mathcal{T}_{0}$ is flat over a smooth base. Because $y$ is flat with smooth fibres, it is then necessarily a smooth morphism, and then this necessarily implies (because $\mathcal{T}_{0}$ is smooth) that $\mathcal{G}$ is a smooth variety. Each fibre of $y$ is naturally a torsor for the normalizer $N$ of $_\textrm{sp} M$ in $\mathbf{G}_{S}$, and one has a natural map $N \times \mathcal{G} \to \mathcal{G} \times_{\mathcal{T}_{0}} \mathcal{G}$ bijective on complex points; since this is a map of two smooth varieties, this means it is necessarily an isomorphism. We conclude that $y$ is an fppf $N$-torsor, and then because $N$ is a smooth algebraic group, an \'etale $N$-torsor by \cite{localtorsorMO}.

Choose an \'etale covering $e : \mathcal{U} \to \mathcal{T}_{0}$ and a trivialization $\sigma : N \times \mathcal{U} \xrightarrow{\sim} \mathcal{G} \times_{\mathcal{T}_{0}} \mathcal{U}$. For $T \in \mathcal{T}_{0}$, construct $g(T)$ by choosing any element $\zeta$ of the fibre $e^{-1}(T)$ and then defining $g(T) = \sigma(1,\zeta)$. That the degree of the resulting field of definition is bounded is a consequence of the fact that the degree of $e$ is bounded, and that the resulting height is polynomial in the height of $T$ is an easy consequence of how heights behave under polynomial maps.
\end{proof}

To give a finiteness criterion for special points, we now need to introduce some language to talk about atypical intersections. Fix a map $h : \mathbb{S} \to \textrm{GAut}(\mathbb{Q}^{m}, Q)_{\mathbb{R}}$ corresponding to a point of $D$, and which factors through a Mumford-Tate subgroup $M \subset \textrm{GAut}(\mathbb{Q}^{m}, Q)$. Then one obtains an induced Hodge structure on the Lie algebra $\mathfrak{m}$ of $M$ through the adjoint action; in particular this is true for the Mumford-Tate group $M = M_{h}$ of the Hodge structure.

\begin{notn}
Given a Hodge structure $h \in D$ with map $\mathbb{S} \to \textrm{GAut}(\mathbb{Q}^{m}, Q)_{\mathbb{R}}$ factoring through a Mumford-Tate group $M$ with Lie algebra $\mathfrak{m}$, we write $\delta(M,h)$ for the sum of the positive Hodge numbers associated to the induced Hodge structure on $\mathfrak{m}$. 
\end{notn}

As discussed in \cite[Ch. II, pg. 56]{GGK}, the quantity $\delta(M,h)$ is the same as the dimension of $\ch{D}_{M} = M(\mathbb{C}) \cdot h$, where we regard $h$ as a point of $\ch{D}$.

\begin{defn}
\label{atyppointdef}
Given a special point $\{ y \} \subset S$ and a pair $(M, h = \widetilde{\varphi}(\widetilde{y}))$ defining $\{ y \}$, we say that $\{ y \}$ is atypical for $(M, h)$ if $\dim S < \dim D_{S} - \delta(M,h)$. We say that $\{ y \}$ is atypical if this inequality holds with $M = M_{h}$.
\end{defn}

\begin{rem}
The stated definition is equivalent to the definition of atypicality for a zero-dimensional Hodge locus given in \cite[Def. 2.2]{BKU} (here we recall that $\varphi$ is assumed to have discrete fibres). This refines the notion of atypicality given in Klingler's paper \cite{atyp}; see \cite[\S5.3]{BKU} for a discussion concerning the relationship between the two notions. 
\end{rem}

\begin{cor}
\label{finspptcor}
Suppose $\mathbf{H}_{S}$ agrees with the derived subgroup of $\mathbf{G}_{S}$. Fix the same data as in \autoref{heightconstrprop}, and suppose additionally that:
\begin{itemize}
\item[(i)] each point $\xi \in \mathcal{S}$ is atypical for the pair $(M_{\widetilde{\xi}}, \widetilde{\varphi}(\widetilde{\xi}))$ which defines it;
\item[(ii)] the $\textrm{Gal}(\overline{\mathbb{Q}}/K)$-action preserves $\mathcal{S}$; and
\item[(iii)] there exists constants $a, b \in \mathbb{R}_{> 0}$ depending only on the fixed data such that 
\[ \widetilde{\theta}(T_{\widetilde{\xi}}) \leq a \, [K(\xi) : K]^{b} , \]
for each $\widetilde{\xi}$ as in (i), where $T_{\widetilde{\xi}}$ is a set of definition of minimal type for $M_{\widetilde{\xi}}$.
\end{itemize}
Then if $\mathcal{S}$ is infinite, there exists $\xi_{1} \in \mathcal{S}$ which lies in a strict positive-dimensional weakly special subvariety of $S$.
\end{cor}

\begin{proof}
We apply \autoref{heightconstrprop}. We then consider the locus
\[ \mathcal{G} := \bigcup_{j = 1}^{r} \underbrace{\{ g M \in (\mathbf{G}_{S}/M)(\mathbb{C}) : (g \cdot M(\mathbb{C}) \cdot P_{j}) \cap \mathcal{I} \neq \varnothing \}}_{\mathcal{G}_{j}} \]
with $P_{1}, \hdots, P_{j}$ as in \autoref{heightconstrprop}. Observe that the construction of \autoref{heightconstrprop} produces a point $g M = g_{\xi} M$ of $\mathcal{G}$ for each $\xi \in \mathcal{S}$. These points are all defined over $\overline{\mathbb{Q}}$, and in particular over a number field of degree at most some fixed constant $d$ by \autoref{heightconstrprop}(1). 

We consider subsets of $(\mathbf{G}_{S}/M)(\mathbb{C})$ definable if they are definable as subsets of $\mathbb{P}^{r}(\mathbb{C})$ after making the identification $\mathbb{C} = \mathbb{R}^2$; note that $\mathcal{G}$ and $\mathcal{G}_{j}$ are definable. Applying \autoref{heightconstrprop}(3) and our assumption (iii), there exists constants $a', b' > 0$ such that
\[ \widetilde{\theta}_{M}(g_{\xi} M) \leq a' \, [K(\xi) : K]^{b'} \]
for all $\xi \in \mathcal{S}$. Now because $\mathcal{S}$ is assumed infinite, using \autoref{heightconstrprop}(2) and (ii), one has for infinitely many positive integers $N$ that
\begin{equation}
\label{galoisorbitmakebig}
 \left| \left\{ g_{\xi} M : \widetilde{\theta}_{M}(g_{\xi} M) \leq a' N^{b'} \right\} \right| \geq \frac{1}{c} N .
\end{equation}

Fix some positive integer $\ell$. We now apply the Pila-Wilkie theorem, and in particular a strengthened form due to Pila \cite{zbMATH05680945}, to deduce that there exists 
\begin{itemize}
\item[(a)] an integer $j$;
\item[(b)] an irreducible complex algebraic subvariety $V \subset (\mathbf{G}_{S}/M)(\mathbb{C})$;
\item[(c)] an analytic open neighbourhood $\mathcal{U} \subset (\mathbf{G}_{S}/M)(\mathbb{C})$ and finitely many analytic components $V_{1}, \hdots, V_{d}$ of $V \cap \mathcal{U}$ such that $V^{\dagger} \subset \mathcal{G}_{j}$ where $V^{\dagger} = V_{1} \cup \cdots \cup V_{d}$; and
\item[(d)] distinct points 
\[ g_{\xi_{1}} M, \hdots, g_{\xi_{\ell}} M \in V^{\dagger} \subset \mathcal{G}_{j} \cap V \]
corresponding to $\xi_{1}, \hdots, \xi_{\ell} \in \mathcal{S}$.
\end{itemize}
To show this, we start by fixing $\ep > 0$, and applying \cite[Thm. 5.3]{zbMATH05680945}, which strengthens \cite[Theorem 1.6]{zbMATH05680945}, to the definable set $\mathcal{G}$ (we take $Z = \mathcal{G}$ to be a trivial family in the notation of \cite[Thm. 5.3]{zbMATH05680945}). One concludes that the points $\{ g_{\xi} M \in \mathcal{G} : \widetilde{\theta}_{M}(g_{\xi} M) \leq T \}$ lie in at most $O(T^{\ep})$ semi-algebraic ``blocks'' (in the sense of \cite[\S3.2]{zbMATH05680945}) contained in $\mathcal{G}$. Taking $\ep \ll 1/b'$ and setting $T = a' N^{b'}$, one concludes from (\ref{galoisorbitmakebig}) that there is a semi-algebraic block $B \subset \mathcal{G}$ and index $j$ such that $B \cap \mathcal{G}_{j}$ contains at least $\ell$ distinct points. Shrinking $B$ we may assume that $B \subset \mathcal{G}_{j}$. 

We now wish to show that $B$ is contained in a complex algebraic subset $V^{\dagger}$ of $\mathcal{G}_{j}$ as in (b) and (c); using the complex algebraic quotient map $q : \mathbf{G}_{S}(\mathbb{C}) \to (\mathbf{G}_{S}/M)(\mathbb{C})$, it suffices to show that each irreducible component of $q^{-1}(B)$ is contained in a complex semi-algebraic subset $\widetilde{V}^{\dagger}$ of $q^{-1}(\mathcal{G}_{j})$ with the same structure with respect to an open neighbourhood $\widetilde{\mathcal{U}} \subset \mathbf{G}_{S}(\mathbb{C})$ and some algebraic variety $\widetilde{V} \subset \mathbf{G}_{S}(\mathbb{C})$, where $\widetilde{V}^{\dagger} = \widetilde{V}_{1} \cup \cdots \cup \widetilde{V}_{\widetilde{d}}$ in an analogous way (afterwards one can compute the complex-algebraically constructible set $q(\widetilde{V})$ and work with its Zariski closure). Choose a closed algebraic embedding $\mathbf{G}_{S} \hookrightarrow \mathbb{A}^{n}$. We observe that 
\[ q^{-1}(\mathcal{G}_{j}) = \{ g : g \cdot P_{j} \in \mathcal{I} \}  \]
is a closed analytic subvariety of some definable open neighbourhood $\Omega \subset \mathbb{A}^{n}(\mathbb{C}) = \mathbb{C}^{n}$. We then apply \autoref{makecomplem} with $Z = q^{-1}(\mathcal{G}_{j})$, and $W$ an irreducible component of $q^{-1}(B) \cap \Omega$. 

The Mumford-Tate group $M$ is connected as an algebraic group (such groups are Zariski closures of connected real tori), so the map $q : \mathbf{G}_{S} \to \mathbf{G}_{S}/M$ has geometrically connected fibres. It then follows (c.f. the proof of \cite[Thm. 5.3.2]{zbMATH03708568}) that $V_{M} := q^{-1}(V)$ is irreducible. Set $\mathcal{U}_{M} := q^{-1}(\mathcal{U})$, and take $\ell > \kappa$, where $\kappa$ is the constant obtained in \autoref{finMTgroups}. Then $V_{M} \cap \mathcal{U}_{M}$ does not map into the orbit $g_{\xi_{1}} \cdot M(\mathbb{C}) \cdot P_{j}$, so there exists an irreducible algebraic curve $C \subset V_{M}$ passing through $g_{\xi_{1}}$ such that the irreducible constructible algebraic set $E := C \cdot M(\mathbb{C}) \cdot P_{j}$ is not contained in $g_{\xi_{1}} \cdot M(\mathbb{C}) \cdot P_{j}$; in particular, $\dim E = \dim O_{c} + 1$, where $O_{c} = c \cdot M(\mathbb{C}) \cdot P_{j}$. Because the intersection $O_{c} \cap \mathcal{I}$ is an isolated point in some neighbourhood of $\widetilde{\varphi}(\widetilde{\xi}_{1})$ when $c = g_{\xi_{1}}$, the intersections $O_{c} \cap \mathcal{I}$ continue to have an isolated point as a component for $c \in C \cap \mathcal{U}_{M}$, shrinking $\mathcal{U}_{M}$ if necessary; moreover, these isolated points can be taken arbitrarily close to $\widetilde{\varphi}(\widetilde{\xi}_{1})$ as $c \to g_{\xi_{1}}$. The number of such isolated points arising from a fixed translate of $M(\mathbb{C}) \cdot P_{j}$ is bounded by $\kappa$, so it follows that $E$ intersects $\mathcal{I}$ in a positive-dimensional analytic locus $F$ passing through $\widetilde{\varphi}(\widetilde{\xi}_{1})$.

Set $\xi = \xi_{1}$ and $c = g_{\xi}$. Now in particular by \autoref{heightconstrprop}(4), one has $c^{-1} M_{\widetilde{\xi}} c = M$ and $\widetilde{\varphi}(\widetilde{\xi})$ is a point of $O = c \cdot M \cdot P_{j}$. The fact that $\xi$ is atypically defined by $M_{\widetilde{\xi}}$ means that $\dim O = \delta(M_{\widetilde{\xi}}, \widetilde{\varphi}(\widetilde{\xi}))$ satisfies
\begin{equation}
\label{herpderpexpint}
\dim S + \dim O < \dim D_{S} = \dim \mathbf{H}_{S}(\mathbb{C}) \cdot \varphi(\widetilde{\xi}) ,
\end{equation}
where for the equality we use the fact that the derived subgroup of $\mathbf{G}_{S}$ is $\mathbf{H}_{S}$.

 Letting $\ch{T} = \mathbf{H}_{S}(\mathbb{C}) \cdot \varphi(\widetilde{\xi})$ one has that
\begin{align*}
\textrm{codim}_{\ch{T}} \, F &< 2 \dim \ch{T} - \dim \ch{T} \\
&\leq 2 \dim \ch{T} - (\dim S + \dim O + 1) \\
&= \textrm{codim}_{\ch{T}} \mathcal{I} + \textrm{codim}_{\ch{T}} {E} ,
\end{align*}
where we have used (\ref{herpderpexpint}), $\dim O + 1 = \dim E$, and $\dim S = \dim \mathcal{I}$. From the Ax-Schanuel Theorem for period mappings \cite{AXSCHAN}, one then learns that $\xi$, which lies inside $\varphi^{-1}(\pi_{D}(F))$ with $\pi_{D} : D \to \Gamma \backslash D$ the natural projection, lies inside a (necessarily positive dimensional) strict weakly special subvariety $W \subset S$, concluding the proof. (We note that in \cite{AXSCHAN}, the Ax-Schanuel theorem uses an inequality formulated for intersections with graphs of period mappings, but one can obtain such an inequality by pulling back an algebraic intersection with $\mathcal{I}$ to an algebraic intersection with the graph.) 
\end{proof}

\begin{lem}
\label{makecomplem}
Let $\Omega \subset \mathbb{C}^{n}$ be a definable connected open domain and let $Z \subset \Omega$ be a closed complex analytic subvariety. Let $W \subset Z$ be an irreducible real-algebraic set (i.e., the intersection with $\Omega$ of an irreducible real-algebraic semi-algebraic subset of $\mathbb{R}^{2n} = \mathbb{C}^{n}$). Then there exists an irreducible complex algebraic variety $V \subset \mathbb{C}^{n}$ and finitely many analytic components $C_{1}, \hdots, C_{d}$ of $\Omega \cap V$ such that $W \subset C_{1} \cup \cdots \cup C_{d} \subset Z$.
\end{lem}

\begin{proof}
The result is essentially shown in \cite[Lemma 4.1]{zbMATH06144656}. In particular, one may follow the proof given in \cite[Lemma 4.1]{zbMATH06144656} as written up until the very last sentence, since no maximality hypothesis on $W$ is invoked before then. In particular, the argument there shows that for every point $O \in W_{\textrm{sm}}$ in the smooth locus $W_{\textrm{sm}} \subset W$, there exists an open neighbourhood $\mathcal{N}_{O} \subset \Omega$ such that $W_{\textrm{sm}} \cap \mathcal{N}_{O} \subset i_{\mathbb{C}}(U(\mathbb{C})) \cap \mathcal{N}_{O} \subset Z$, where $V := i_{\mathbb{C}}(U(\mathbb{C}))$ is a complex algebraic subvariety of $\mathbb{C}^{n}$ independent of $O$ constructed in \cite[Lemma 4.1]{zbMATH06144656}. Since $W$ is irreducible and contained in $i_{\mathbb{C}}(U(\mathbb{C}))$, there exists a complex-algebraic irreducible component of $i_{\mathbb{C}}(U(\mathbb{C}))$ which contains $W$, and we may replace $V$ with this component.

Now let $C_{1}, \hdots, C_{d}$ be the irreducible analytic components of $V \cap \Omega$ that intersect at least one of the $\mathcal{N}_{O}$ (there are finitely many components of $V \cap \Omega$ by definability). Then necessarily each such $C_{i}$ satisfies $C_{i} \subset Z$: the intersection $C_{i} \cap Z$ is a closed complex analytic subset of $C_{i}$ of the same dimension (since $C_{i} \cap \mathcal{N}_{O} \subset C_{i} \cap Z$ and $Z$ is closed in $\Omega$), and $C_{i}$ is assumed irreducible, so $C_{i} \cap Z = C_{i}$. From the construction above, $W_{\textrm{sm}} \subset \bigcup_{i=1}^{d} C_{i}$, and then since the $C_{i}$ are closed in $\Omega$ and $W_{\textrm{sm}}$ is dense in $W$, we have $W \subset \bigcup_{i=1}^{d} C_{i} \subset Z$. 
\end{proof}

\subsection{Heights in a Siegel Set Orbit: The Case of Idempotents}
\label{idemheightsec}

\begin{notn}
For a Hodge structure $h \in D$, we write $h_{\mathbb{Q}}$ for the associated $\mathbb{Q}$-Hodge structure, and $\textrm{End}(h_{\mathbb{Q}})$ for its algebra of endomorphisms.
\end{notn}

\begin{defn}
A \emph{central Hodge idempotent} associated to $h \in D$ is a central idempotent $e$ in $\textrm{End}(h_{\mathbb{Q}})$. It is called primitive if it does not decompose as $e = e_{1} + e_{2}$, with both $e_{1}$ and $e_{2}$ non-zero central idempotents of $\textrm{End}(h_{\mathbb{Q}})$. 
\end{defn}

\begin{defn}
We say a pure polarized Hodge structure $(F^{\bullet}, W, Q)$ on the polarized $\mathbb{Q}$-vector space $(W, Q)$ is \emph{degenerate} if it is of even weight $n = 2m$ with $F^{m} = W_{\mathbb{C}}$. Otherwise it is non-degenerate.
\end{defn}

We note that in the degenerate case the associated space of polarized Hodge structures $D_{W}$ has dimension zero.

\begin{lem}
\label{genEndQ}
Suppose that $h$ is a non-degenerate Hodge structure on the polarized lattice $(W, Q)$. Then the orbit $O = \textrm{SAut}(W,Q)(\mathbb{R}) \cdot h$ contains a Hodge structure with $\mathbb{Q}$-simple special Mumford-Tate group $\textrm{SAut}(W,Q)$.
\end{lem}

\begin{proof}
That a general Hodge structure in $O$ has special Mumford-Tate group $\textrm{SAut}(W,Q)$ is then \cite[II.A.3, II.A.6]{GGK}. Note that the statements \cite[II.A.3, II.A.6]{GGK} should use $\textrm{SAut}$ instead of $\textrm{Aut}$ because the automorphism group of the polarized vector space $(W,Q)$ is not connected in general. The only situation where $\textrm{SAut}(W,Q)$ is not $\mathbb{Q}$-simple is when it is isomorphic to $\textrm{SO}_{2}$, and this case is excluded by the assumptions (see the discussion on \cite[pg. 51]{GGK}). 
\end{proof}

\begin{lem}
\label{centdecomplem}
Let $h \in D$ be a Hodge structure and let $1 = e_{1} + \cdots + e_{\ell}$ be a decomposition into associated central Hodge idempotents. Then
\begin{itemize}
\item[(i)] we have an induced decomposition $\mathbb{Q}^{m} = V_{1} \oplus \cdots \oplus V_{\ell}$ into sub-Hodge structures;
\item[(ii)] the polarization decomposes as $Q = Q_{1} + \cdots + Q_{\ell}$, with $Q_{i}$ the pullback under $e_{i}$ of a polarization on $V_{i}$;
\item[(iii)] if each $e_{i}$ is primitive and the Hodge structure on $V_{i}$ induced by $h$ is non-degenerate, then there exists a Hodge structure $h' \in D$ whose special Mumford-Tate group is 
\[ M := \textrm{SAut}(V_{1}, Q_{1}) \times \cdots \times \textrm{SAut}(V_{\ell}, Q_{\ell}) \subset  \textrm{SAut}(\mathbb{Q}^{m}, Q) ; \]
moreover, we have $h \in M(\mathbb{R}) \cdot h'$. 
\end{itemize}
\end{lem}

\begin{proof}
Part (i) is clear. For (ii) we reduce inductively to the case $\ell = 2$. From \cite[Cor. 2.12]{zbMATH05233837} we obtain a decomposition $\mathbb{Q}^{m} = V_{1} \oplus V^{\perp}_{1}$ by sub-Hodge structures. Then letting $f$ be an idempotent for which $f(V_{1}) = 0$ and projecting onto $V^{\perp}_{1}$, we have
\[ V^{\perp}_{1} = \textrm{im}(f \cdot 1) = \textrm{im}(fe_{2}) = \textrm{im}(e_{2}f) \subset V_{2} . \]
Hence $V_{2} = V^{\perp}_{1}$ by dimension counting. That $Q$ decomposes as $Q = Q_{1} + Q_{2}$ is then immediate from orthogonality.

For (iii) we have by \autoref{genEndQ} that the space of Hodge structures on $V_{i}$ polarized by $Q_{i}$ has generic special Mumford-Tate group $M_{i} := \textrm{SAut}(V_{i}, Q_{i})$. Now let $D_{\times} := D_{1} \times \cdots \times D_{\ell} \subset D$ be the product of the orbits $M_{i}(\mathbb{R}) \cdot h_{i}$ where $h = h_{1} \oplus \cdots \oplus h_{\ell}$, and let $D^{\circ}_{\times} \subset D_{\times}$ be the connected component containing $h$. Let $M' \subset M$ be the special Mumford-Tate group of a $\mathbb{Q}$-generic point $h' \in D^{\circ}_{\times}$ (outside any $\mathbb{Q}$-algebraic subvarieties of $\ch{D}$ which intersect $D^{\circ}_{\times}$ properly). We note that $h \in M(\mathbb{R}) \cdot h'$ by construction, so it suffices to show that $M' = M$. Then $M(\mathbb{Q})^{\circ} := M(\mathbb{Q}) \cap M(\mathbb{R})^{\circ}$ stabilizes $M'$ under the conjugation action, where $M(\mathbb{R})^{\circ} \subset M(\mathbb{R})$ denotes the identity component. The group $M(\mathbb{Q})^{\circ}$ is Zariski dense in $M$ (use for instance that connected $\mathbb{Q}$-algebraic groups are unirational \cite[Theorem 18.2]{lagbor}, hence $M(\mathbb{Q})$ is analytically dense in $M(\mathbb{R})$, together with the fact that $M(\mathbb{R})$ has finitely many connected components), so it follows that $M'$ is normal in $M$. 

Write $\pi_{i} : M \to M_{i}$ for the natural projection. Since the $\mathbb{Q}$-generic point $h'$ of $D^{\circ}_{\times}$ necessarily maps to a $\mathbb{Q}$-generic point of $D_{i}$ under the projection $D^{\circ}_{\times} \to D_{i}$, it follows that $\pi_{i}(M') = M_{i}$. From the characterization of normal subgroups of a product of semisimple groups in \cite[Thm. 21.51]{zbMATH06713849}, it follows that $M' = M$.
\end{proof}

\begin{lem}
\label{idemselfadjoint}
For a Hodge structure $h \in D$, the central Hodge idempotents of $h$ are self-adjoint for $\psi_{h} = Q(C_{h}(-),(-))$.
\end{lem}

\begin{proof}
Let $e$ be such an idempotent and set $e' = 1 - e$. Applying \autoref{centdecomplem}(ii) we obtain an orthogonal decomosition $\psi_{h} = \psi + \psi'$ compatible with the decomposition $\mathbb{Q}^{m} = V \oplus V'$ of \autoref{centdecomplem}(i); self-adjointness is then immediately verified.
\end{proof}

\begin{lem}
\label{finmanyorbits}
Fix a Siegel set $\mathfrak{O} = \mathfrak{G} \cdot h$ for $D$ with $\mathfrak{G} = \Omega A_{t} K_{h} \subset G(\mathbb{R})$, and $G = \textrm{Aut}(\mathbb{Q}^{m}, Q)$. Then central Hodge idempotents associated to points of $\mathfrak{O}$ lie in finitely many orbits of $\mathfrak{G}$ in $\textrm{End}(\mathbb{R}^{m})$.
\end{lem}

\begin{proof}
Given any central Hodge idempotent $e$ associated to a point of $D$, we obtain a Mumford-Tate group $M_{e} \subset G$ defined as its stabilizer. Applying \autoref{finmanyconjclasses} we may reduce to the case where all $M_{e}$ belong to a fixed $G(\mathbb{R})$-conjugacy class. Then associated to each $M_{e}$, there are finitely many $M_{e}(\mathbb{R})$-orbits $D_{e}$ of Hodge structures in $D$ with Mumford-Tate group contained in $M_{e}$ as a consequence of \cite[VI.A.2]{GGK}. Thus we may reduce to considering just those $e$ associated to a pair $(M_{e}, D_{e})$ which lies in a single $G(\mathbb{R})$-equivalence class.

Now fix a pair $(M_{e}, D_{e})$ associated to a point $x \in \mathfrak{O} \cap D_{e}$, and let $E \subset \textrm{End}(\mathbb{Q}^{m})$ be the vector subspace stabilized by $M_{e}$ under the action by conjugation. Then $E$ is the semisimple algebra associated to a generic Hodge structure in $D_{e}$. The algebra $E_{\mathbb{R}}$, being finite dimensional over the field $\mathbb{R}$, has finitely many central idempotents $e = e_{1}, \hdots, e_{\ell}$ (see \cite[Def. 1.7.3, Thm. 1.7.6, Thm. 1.7.7]{zbMATH06854373}). Pick some $r \in \mathfrak{G}$ sending $h$ to $x$, and define $f_{i} = r^{-1} \cdot e_{i}$.

Let $(M_{e'}, D_{e'})$ be an equivalent pair associated to $x' \in \mathfrak{O} \cap D_{e'}$. Then there exists $g \in G(\mathbb{R})$ which sends $(M_{e}, D_{e})$ to $(M_{e'}, D_{e'})$, and then adjusting by an element of $M_{e'}(\mathbb{R})$ we may assume $g$ sends $x$ to $x'$. If $E' \subset \textrm{End}(\mathbb{Q}^{m})$ is the subspace stabilized by $M_{e'}$, then $E_{\mathbb{R}}$ and $E'_{\mathbb{R}}$ are conjugated by $g$. This conjugation sends central idempotents to central idempotents, hence the central Hodge idempotents of $x'$ are among $g \cdot e_{1}, \hdots, g \cdot e_{\ell}$. By \autoref{invimageSiegel} the set $\mathfrak{G}$ is the preimage of $\mathfrak{O}$ in $G(\mathbb{R})$ under the orbit map, so $gr \in \mathfrak{G}$. It follows that the central Hodge idempotents associated to $x'$ are in the orbit of $\mathfrak{G}$ on $f_{1}, \hdots, f_{\ell}$.
\end{proof}

We say that a map $\alpha : h \to h'$ of Hodge structures on $\mathbb{Z}^{m}$ is an isogeny if it induces an isomorphism of $\mathbb{Q}$-Hodge structures. Its degree is the magnitude of its determinant.

\begin{lem}
\label{isogenyheightbound}
Given a Siegel set $\mathfrak{O} = \mathfrak{G} \cdot h$ for $D$, the heights of the central Hodge idempotents associated to points of $\mathfrak{O}$ are bounded by a polynomial in their isogeny height: the degree of the minimal isogeny $\alpha : h \to h'$ satisfying the property that the Hodge decomposition of $h'$ is induced by a direct sum decomposition of $\mathbb{Z}^{m}$ as a $\mathbb{Z}$-module. 
\end{lem}

\begin{proof}
Applying \autoref{finmanyorbits} we may reduce to considering just those central Hodge idempotents in a single orbit of $\mathfrak{G}$. After embedding $\mathfrak{G}$ into a finite union of Siegel sets for $\GL_{m}$ using \autoref{bettersiegelembedding}, one may apply \autoref{polyindenomlem} and \autoref{idemselfadjoint} to reduce to bounding the size of the denominators of the central Hodge idempotents associated to points of $\mathfrak{O}$. Choose $\alpha : h \to h'$ as in the statement, and let $e'$ be the idempotent of $h'$ corresponding to $e$. Then $e'$ is integral, hence so is $(\det \alpha) \alpha^{-1} e' \alpha = (\det \alpha) e$. But then $|\det \alpha|$ must be larger than the largest denominator of $e$.
\end{proof}

\subsection{Application in the Abelian Case}
\label{hypotheseschecksec}

We now consider the more concrete case of Hodge structures arising from abelian varieties; in particular, we start by verifying hypothesis (ii) of \autoref{finspptcor} in this case. 

In what follows we will use the bundle $P = P'$ introduced in \S\ref{monocompframesec} with $\mathbb{L} = \mathbb{V}_{\mathbb{C}}$ and $\mathcal{V} = \mathcal{H}$. We take the point $c$ appearing in (\ref{Pdefeq}) to be the same as the point $c$ we chose in \S\ref{Hodgetheosetupsec}. Choosing $c \in S(K)$, we may also assume that $P$ is defined over $K$ with respect to the natural $K$-structure on $\mathbb{B}(\sheafhom(\mathcal{\mathcal{H}}, \mathcal{O}_{S} \otimes_{\mathbb{C}} \mathcal{\mathcal{H}}_{c}))$. Indeed, because $\mathbf{H}_{S}$ is semisimple, it is defined by the tensors it fixes pointwise (see \cite[I, Prop. 3.1]{zbMATH03728195}). One then considers a large enough submodule $\mathcal{M} \subset \bigoplus_{a, b \geq 0} \mathcal{H}^{a,b}$ generated by global $\nabla$-flat sections so that the formula in (\ref{Pdefeq}) defines an $\mathbf{H}_{S}$-torsor. We may assume $\mathcal{M}$ is defined over $K$ because the connection $\nabla$ is \cite{katz1968}. 

We write $\widetilde{\mathcal{P}} \subset P(\mathbb{C})$ for the unique analytic continuation of $\mathcal{P} \subset P(\mathbb{C})$ to a closed analytic locus in $P^{\textrm{an}}$, as in \autoref{axschancor} and its proof.

\begin{lem}
\label{conjhypothabcase}
Suppose in the situation of \autoref{finspptcor} that $f$ is an abelian family. Then after replacing $K$ with a finite extension, hypothesis (ii) holds for the set $\mathcal{S}$ consisting of all special points defined atypically by a pair $(M_{\widetilde{\xi}}, \widetilde{\varphi}(\widetilde{\xi}))$ which is $\mathbf{G}_{S}(\mathbb{C})$-conjugate to a fixed Mumford-Tate subgroup $M \subset \mathbf{G}_{S}$.
\end{lem}

\begin{proof}
The assumption that $f$ is an abelian family will be used only through Deligne's verification of the absolute Hodge conjecture \cite{deligne} for tensors associated to the cohomology of an abelian variety. We will use this throughout.

\vspace{0.5em}

Without loss of generality we may assume that $c$ belongs to $\mathcal{S}$ (otherwise there is nothing to show), and then take $(M = M_{\widetilde{c}}, \widetilde{\varphi}(\widetilde{c}))$ to be a pair defining $c$. The absolute Hodge conjecture implies that the tensors stabilized by $M_{\widetilde{c}}$ (resp. $\mathbf{G}_{S,\widetilde{c}}$) are defined over $\overline{\mathbb{Q}}$ with respect to the $\overline{\mathbb{Q}}$-structure on the fibre $\mathcal{H}_{c}$, and hence we obtain a de Rham realization $M_{c} \subset \GL(\mathcal{H}_{c})$ (resp. $\mathbf{G}_{S,c} \subset \GL(\mathcal{H}_{c})$), which we assume is defined over $K$ after passing to a finite extension. This in particular implies that the $\mathbf{G}_{S,c}(\mathbb{C})$-conjugacy class of $M_{c}$ is defined over $K$ (i.e., stable under $\textrm{Aut}(\mathbb{C}/K)$). 

Recall that our convention, established in \autoref{MTsubgpdef} and \autoref{defaspvar}, is to view all Mumford-Tate groups associated to points $\xi \in \mathcal{S}$ as subgroups of the fixed realization of $\mathbf{G}_{S}$ inside $\GL(\mathbb{V}_{\widetilde{c}})$, which we do by first choosing a lift $\widetilde{\xi}$ and then identifying $\widetilde{\mathbb{V}}_{\widetilde{\xi}}$ and $\widetilde{\mathbb{V}}_{\widetilde{c}}$ using the unique (up to homotopy) path $\alpha = \alpha(\widetilde{\xi}, \widetilde{c})$ from $\widetilde{\xi}$ to $\widetilde{c}$. We now make this explicit rather than implicit, writing $\underline{\alpha} : \widetilde{\mathbb{V}}_{\widetilde{\xi}} \xrightarrow{\sim} \widetilde{\mathbb{V}}_{\widetilde{c}}$ for the induced linear map. More specifically, we consider a group $M_{\xi} \subset \GL(\mathbb{V}_{\xi})$ defined as the stabilizer of tensors associated to $\mathbb{V}_{\xi}$ which are Hodge for $F^{\bullet}_{\xi}$, let $M_{\widetilde{\xi}} \subset \GL(\widetilde{\mathbb{V}}_{\widetilde{\xi}})$ be its lift, and write $\underline{\alpha}(M_{\widetilde{\xi}})$ for its translate to $\GL(\widetilde{\mathbb{V}}_{\widetilde{c}})$. Then our assumption on $\xi \in \mathcal{S}$ is that there exists such a pair $(M_{\xi}, \widetilde{\xi})$ such that $(\underline{\alpha}(M_{\widetilde{\xi}}), \widetilde{\varphi}(\widetilde{\xi}))$ defines $\{ \xi \}$ atypically and $\underline{\alpha}(M_{\widetilde{\xi}})$ is $\mathbf{G}_{S,\widetilde{c}}(\mathbb{C})$-conjugate to $M_{\widetilde{c}}$.

Now take $\sigma \in \textrm{Aut}(\mathbb{C}/K)$, and consider the conjugate pair $(\xi^{\sigma}, M_{\xi^{\sigma}} := M^{\sigma}_{\xi})$, where the conjugation of Mumford-Tate groups uses the identification $\mathcal{H}_{\xi,\mathbb{C}} \xrightarrow{\sigma} \mathcal{H}_{\xi^{\sigma},\mathbb{C}}$. The absolute Hodge conjecture guarantees that $M_{\xi^{\sigma}}$ is again defined as the stabilizer of tensors associated to $\mathbb{V}_{\xi^{\sigma}}$ which are Hodge for $F^{\bullet}_{\xi^{\sigma}}$. Choose a lift $\widetilde{\xi}^{\sigma} \in \widetilde{S}$ of $\xi^{\sigma}$ with associated path $\tau$ from $\widetilde{\xi}^{\sigma}$ to $\widetilde{c}$ and linear map $\underline{\tau}$, and write $M_{\widetilde{\xi}^{\sigma}}$ for the lift of $M_{\xi^{\sigma}}$ to $\GL(\widetilde{\mathbb{V}}_{\widetilde{\xi}^{\sigma}})$. We aim to show that $(\underline{\tau}(M_{\widetilde{\xi}^{\sigma}}), \widetilde{\varphi}(\widetilde{\xi}^{\sigma}))$ defines $\{ \xi^{\sigma} \}$ atypically, and $\underline{\tau}(M_{\widetilde{\xi}^{\sigma}})$ is $\mathbf{G}_{S,\widetilde{c}}(\mathbb{C})$-conjugate to $M_{\widetilde{c}}$. 

\vspace{0.5em}

Now both $\underline{\alpha}$ and $\underline{\tau}$ induce, using the natural maps $\widetilde{\mathbb{V}}_{\mathbb{C}, \widetilde{\xi}} \xrightarrow{\sim} \mathbb{V}_{\mathbb{C}, \xi}$, $\widetilde{\mathbb{V}}_{\mathbb{C}, \widetilde{\xi}^{\sigma}} \xrightarrow{\sim} \mathbb{V}_{\mathbb{C}, \xi^{\sigma}}$ and $\widetilde{\mathbb{V}}_{\mathbb{C}, \widetilde{c}} \xrightarrow{\sim} \mathbb{V}_{\mathbb{C}, c}$, points $\underline{\alpha}_{\textrm{dR}}, \underline{\tau}_{\textrm{dR}} \in \widetilde{\mathcal{P}}$. Applying $\sigma$ and using that $P$ is defined over $K$, one obtains a point $\underline{\alpha}^{\sigma}_{\textrm{dR}} \in \Hom(\mathcal{H}_{\xi^{\sigma},\mathbb{C}}, \mathcal{H}_{c,\mathbb{C}}) \cap P(\mathbb{C})$. Using the $\mathbf{H}_{S,c}$-torsor structure of $P$ one may find an element $\beta \in \mathbf{H}_{S,c}(\mathbb{C}) \subset \GL(\mathcal{H}_{c,\mathbb{C}})$ such that $\beta \circ \underline{\alpha}_{\textrm{dR}}^{\sigma} = \underline{\tau}_{\textrm{dR}}$. Using the natural identifications, one knows that $\underline{\alpha}_{\textrm{dR}} \circ M_{\xi} \circ \underline{\alpha}_{\textrm{dR}}^{-1}$ is $\mathbf{G}_{S,c}(\mathbb{C})$-conjugate to $M_{c}$, and hence acting with $\sigma$, that $\underline{\alpha}^{\sigma}_{\textrm{dR}} \circ M_{\xi^{\sigma}} \circ (\underline{\alpha}^{\sigma}_{\textrm{dR}})^{-1}$ is $\mathbf{G}_{S,c}(\mathbb{C})$-conjugate to $M_{c}$, as the $\mathbf{G}_{S,c}(\mathbb{C})$-conjugacy class of $M_{c}$ is $\textrm{Aut}(\mathbb{C}/K)$-stable. But then since $\beta \in \mathbf{H}_{S,c}(\mathbb{C}) \subset \mathbf{G}_{S,c}(\mathbb{C})$, one learns that $\underline{\tau}_{\textrm{dR}} \circ M_{\xi^{\sigma}} \circ \underline{\tau}_{\textrm{dR}}^{-1}$ is $\mathbf{G}_{S,c}(\mathbb{C})$-conjugate to $M_{c}$. On the lifted level, this then shows that $\underline{\tau}(M_{\widetilde{\xi}^{\sigma}})$ is $\mathbf{G}_{S,\widetilde{c}}(\mathbb{C})$-conjugate to $M_{\widetilde{c}}$. 

\vspace{0.5em}

For the atypicality claim, we note that the quantity $\delta(M_{\widetilde{\xi}^{\sigma}}, h_{\widetilde{\xi}^{\sigma}})$ may also be computed from the Hodge filtration on the Lie algebra of $M^{\sigma}_{\xi}$ (instead of the Hodge direct sum decomposition). In particular, it can be computed using the de Rham realizations $M_{\xi^{\sigma}} = M^{\sigma}_{\xi}$ and $F^{\bullet}_{\xi^{\sigma}}$, and this data is $\sigma$-conjugate to the de Rham realizations of $M_{\xi}$ and $F^{\bullet}_{\xi}$. The fact that both $\{ \xi \}$ and $\{ \xi^{\sigma} \}$ are isolated zero-dimensional Hodge loci (special points) follows because $\textrm{Aut}(\mathbb{C}/K)$ preserves the positive-dimensional Hodge locus (c.f. \cite[\S1.3(**)]{voisin2010hodge} and the sentences that follow that statement). 
\end{proof}

We write $S_{\textrm{Hg},+}$ for the positive-dimensional Hodge locus: the union in $S(\mathbb{C})$ of all positive-dimensional special subvarieties. 

\begin{prop}
\label{idempotcor}
Suppose that $f : X \to S$ is an abelian family for which the derived subgroup of $\mathbf{G}_{S}$ is $\mathbb{Q}$-simple, and that $\mathbb{V} = R^{1} \an{f}_{*} \mathbb{Z}$. Let $\mathcal{S} \subset S(\mathbb{C})$ be a set of special points not lying in $S_{\textrm{Hg},+}$ which are defined by, and atypical for, their central Hodge idempotents in the sense of \autoref{defaspvar2} and \autoref{atyppointdef}. Suppose that
\[ \theta(\xi) \leq \kappa \, [K(\xi) : K]^{a} \]
for all $\xi \in \mathcal{S}$, for some $\kappa, a \in \mathbb{R}_{> 0}$ independent of $\xi$, with $\theta$ some logarithmic Weil height. Then $\mathcal{S}$ is finite.
\end{prop}

\begin{proof}
We adopt the same setup as in \autoref{finspptcor}, noting that hypothesis (ii) is now verified as a consequence of \autoref{conjhypothabcase} after possibly replacing $K$ with a finite extension and $\mathcal{S}$ with its Galois orbit. Hypothesis (i) in \autoref{finspptcor} has been assumed, and the hypothesis on $\mathbf{H}_{S}$ is a consequence of $\mathbb{Q}$-simplicity and the normality of $\mathbf{H}_{S}$ in the derived subgroup of $\mathbf{G}_{S}$ \cite{Andre1992} (note that we may assume $f$ is not iso-trivial). If we can verify hypothesis (iii) in \autoref{finspptcor}, then, assuming $\mathcal{S}$ is infinite, we will conclude that some point $\xi$ in $\mathcal{S}$ lies in a strict positive-dimensional weakly special subvariety. Then $\xi$ is therefore contained in a positive-dimensional strict weakly special subvariety. The $\mathbb{Q}$-simplicity of monodromy then implies this subvariety may be taken to be weakly non-factor in the sense of \cite[Def. 1.13]{fieldsofdef}, and hence also necessarily in the positive-dimensional Hodge locus by \cite[Lem. 2.5]{fieldsofdef}, contradicting our assumption on $\mathcal{S}$.

By hypothesis, there exists a pair $(M_{\widetilde{\xi}}, \widetilde{\varphi}(\widetilde{\xi}))$ defining $\{ \xi \}$ in the sense of \autoref{defaspvar}, where $M_{\widetilde{\xi}}$ is the stabilizer of the central Hodge idempotents of $\widetilde{\varphi}(\widetilde{\xi})$. To verify (iii), we now take the set $T_{\widetilde{\xi}}$ in \autoref{finspptcor} to consist of the span of the primitive central Hodge idempotents $e_{1}, \hdots, e_{\ell}$ of $\widetilde{\varphi}(\widetilde{\xi})$; here we are using the fixed polarization $Q$ to make the identification $\textrm{End}(\mathbb{Q}^{m}) \cong (\mathbb{Q}^{m})^{\otimes 2}$, which is compatible with the action of $\textrm{Aut}(\mathbb{Q}^{m}, Q)$. Applying \autoref{sympinvar} below, one has $_\textrm{sp} M_{\widetilde{\xi}} = M$, where $M$ is as in \autoref{centdecomplem}(iii), and we may likewise assume that the $h'$ of \autoref{defaspvar} associated to our pair $(M_{\widetilde{\xi}}, \widetilde{\varphi}(\widetilde{\xi}))$ is the $h'$ of \autoref{centdecomplem}(iii). Since $\textrm{End}(h'_{\mathbb{Q}}) = \textrm{span}_{\mathbb{Q}} \{ e_{1}, \hdots, e_{\ell} \}$ is exactly the subspace of $\textrm{End}(\mathbb{Q}^{m})$ stabilized by $M$ and $M$ does not stabilize any elements of $\mathbb{Q}^{m}$ (weight one Hodge structures have no Hodge vectors), it follows that $T_{\widetilde{\xi}}$ is minimal type.

Applying \autoref{isogenyheightbound} it now suffices to show that the isogeny heights of the associated weight one Hodge structures (i.e., abelian varieties) are polynomially bounded by the Faltings height (which differs from any logarithmic Weil height by a multiplicative constant) of the fibre above the point $\xi$ and the degree of the field of definition of $\xi$. This is just the result of Masser-Wustholtz \cite[Theorem 1.1]{zbMATH00847746}. Note that the statement of \cite[Theorem 1.1]{zbMATH00847746} is given (in our language) for $\xi$ defined over a field of bounded degree, but the constants depend polynomially on this degree as is explained at the end of \cite[pg.23]{zbMATH00847746}.
\end{proof}

\begin{lem}
\label{sympinvar}
Let $W$ be a $\mathbb{Q}$-vector space with symplectic form $Q$, and $W = W_{1} \oplus \cdots \oplus W_{\ell}$ be an orthogonal decomposition with idempotents $e_{1}, \hdots, e_{\ell}$. Then the subspace in $W \otimes W^{*}$ stabilized by $H = \textrm{Sp}(W_{1}) \times \cdots \times \textrm{Sp}(W_{\ell})$ is the span of $e_{1}, \hdots, e_{\ell}$.
\end{lem}

\begin{proof}
We can write
\begin{equation}
\label{subspdecomp}
W \otimes W^{*} = \bigoplus_{i,j} \underbrace{W_{i} \otimes W_{j}^{*}}_{W_{ij}} .
\end{equation}
When $i \neq j$, the action on the spaces $W_{i} \otimes W_{j}^{*}$ factors through $\textrm{Sp}(W_{i}) \times \textrm{Sp}(W_{j})$. If $t \in W_{ij}$ is a proposed invariant we can write $t = \sum_{k=1}^{m} w_{ik} \otimes w^{*}_{jk}$ with $m$ minimal and observe that necessarily the $w_{ik}$ (resp. the $w^{*}_{jk}$) are fixed by $\textrm{Sp}(W_{i})$ (resp. $\textrm{Sp}(W_{j})$). Since the symplectic group fixes no non-zero vectors, this implies that $t = 0$.

As the decomposition (\ref{subspdecomp}) is invariant under the $H$-action an invariant in $W \otimes W^{*}$ is a sum of invariants in each individual $W_{ij}$, so we are reduced to computing to computing the degree two tensor invariants of the symplectic groups $\textrm{Sp}(W_{i})$. This is standard (e.g., use \cite[Theorem 1B]{zbMATH04103278} and note that $S_{2}$ acts on $\theta_2$ through multiplication by $-1$). 
\end{proof}

We note that \autoref{idempotcor} gives an interpretation of \autoref{informalstatement} in the introduction.

\section{Applications}

We now return to the setting where $S$ is a curve. In what follows we will use the following fact, which is part of the proof of Theorem 2 in \cite[IX, \S4.4]{zbMATH00041964}.

\begin{thm}
\label{numberofvancycles}
Suppose that $f : X \to S$ is a semistable projective family of relative dimension $n$ with geometrically connected fibres, and let $f' : X' \to S'$ be the base-change to the smooth locus $S' \subset S$, which we assume is non-empty. Then if $s_{0} \in S(\mathbb{C}) \setminus S'(\mathbb{C})$, the vanishing cycles induced by the order-$n$ normal crossing singularities of $Y = X_{s_{0}}$ span a space of dimension $h^{n}(\Sigma_{Y})$, where $h^{n}(\Sigma_{Y})$ is the dimension of the $n$'th cohomology group of the dual graph $\Sigma_{Y}$ associated to $Y$ (see \cite[\S1]{zbMATH06290383}). \qed
\end{thm}

\subsection{Families of Curves}

We now prove \autoref{modulicurvesthm}. We first handle the case where $K$ has positive transcendence degree. We start with the following.

\begin{prop}
\label{ZardenseKtransprop}
Suppose $Y$ is a an algebraic variety defined over $\overline{\mathbb{Q}}$ and $S$ is a $\overline{\mathbb{Q}}$-Zariski dense curve in $Y_{K}$, where $K$ is an extension of $\overline{\mathbb{Q}}$ of positive transcendence degree. Then $S$ intersects the union of all codimension two $\overline{\mathbb{Q}}$-subvarieties of $Y$ in a finite set. 
\end{prop}

\begin{proof}
The $\overline{\mathbb{Q}}$-density assumption implies that $S$ does not lie in the singular locus of $Y$, so removing a closed locus we may assume $Y$ is smooth. Given a finite cover $Y = \bigcup_{i} Y_{i}$ by Zariski opens, we may reduce to the same problem for $(Y_{i}, Y_{i} \cap S)$, and therefore assume (using smoothness) that $Y$ admits an \'etale map $c : Y \to \mathbb{A}^{n}$ defined over $\overline{\mathbb{Q}}$. Reducing to the same question for the Zariski closure of $c(S)$, this is then  \cite[Theorem 1.2]{zbMATH06256440}.
\end{proof}

\noindent Now take $Y$ to be the $\overline{\mathbb{Q}}$-Zariski closure of $S$ in (the underlying coarse space of) $\mathcal{M}_{g}$. Let $\bigcup_{i=1}^{\infty} \mathcal{Z}_{i} \subset \mathcal{M}_{g}(\mathbb{C})$ be the countable union of special subvarieties of $\mathcal{M}_{g}$ which parameterize curves whose Jacobian has a CM isogeny factor. Since $Y$ has Mumford-Tate group $\textrm{GSp}_{2g}$ (recall that $S$ is Hodge-generic) and $g \geq 2$, a component $C \subset \mathcal{Z}_{i} \cap Y$ of codimension $1$ in $Y$ is atypical for the variation of Hodge structure on $Y$ in the sense of \cite[Def. 2.2]{BKU} (the special subvarieties $\mathcal{Z}_{i}$ arise as components of inverse images under the Torelli map of special subvarieties of $\mathcal{A}_g$ of codimension at least $2$). Since $K$ has positive-transcendence degree, $Y$ has dimension $\geq 2$, so such components are also positive dimensional. Thus by geometric Zilber-Pink \cite[Thm. 3.1]{BKU}, the union $U$ of all such $C$ lies in a strict Zariski closed subvariety $Z \subset Y$; note that in the notation of \cite[Thm. 3.1]{BKU} we cannot have $Z = Y$ because $\mathbf{H}_{Y} = \textrm{Sp}_{2g}$ is (absolutely) simple, where we use the Andr\'e-Deligne monodromy theorem \cite{Andre1992} to deduce that $\mathbf{H}_{S}$ is normal in $\textrm{GSp}_{2g}$. Because each component of $U$ is defined over $\overline{\mathbb{Q}}$, we may assume $Z$ is as well. Then $Z$, and hence $U$, intersects $S$ in finitely many points. The other components of $Y \cap \bigcup_{i=1}^{\infty} \mathcal{Z}_{i}$ do as well by \autoref{ZardenseKtransprop}, since they necessarily have codimension at least two in $Y$.  

\vspace{0.5em}

Now assuming $K$ is a number field, we observe that the case where $g = 2$ was already proven in \cite[Theorem 1.1]{zbMATH07481643}, albeit under different language. To reduce from our setup to the setup of Daw and Orr in \cite{zbMATH07481643}, note that the locus $\mathcal{B} \subset \overline{\mathcal{M}_{g}} \setminus \mathcal{M}_{g}$ in the case of $g = 2$ necessarily consists of singular curves $C = \mathbb{P}^{1} \cup \dots \cup \mathbb{P}^{1}$. The stability condition requires (see \cite[Def. 1.1]{zbMATH03289080}) that each $\mathbb{P}^{1}$ component intersects each other component in at least $3$ nodes. The boundary strata of $\overline{\mathcal{M}_{g}}$ are stratified by the number of nodes \cite{zbMATH03882570}, and the stratum of curves with $k$ nodes has dimension $3g - 3 - k$; in particular, the only possibility here is $C = \mathbb{P}^{1} \cup \mathbb{P}^{1}$ with the components meeting at $3$ nodes. 

Now let $f : \overline{X} \to \overline{S}$ be the family of stable curves over $\overline{S}$ coming from the moduli interpretation of $\mathcal{M}_{g}$, and let $b \in \overline{S} \cap \mathcal{B}$. Then using \cite[Thm. 4.18.1]{kleiman2005picardscheme} we can consider the relative Jacobian $\textrm{Pic}^{0}_{\overline{X}/\overline{S}}$, which is a group scheme over $\overline{S}$. Its fibre above $b$ is $\textrm{Pic}^{0}_{C}$, which by \cite[\S5.B, pg. 250]{zbMATH01179513} is isomorphic to an extension of the Jacobian of $\mathbb{P}^{1} \sqcup \mathbb{P}^{1}$, which is zero, by $\mathbb{G}^{t}_{m}$, where $t$ is the first Betti number of the dual graph associated to $C$ (as defined in loc. cit.). One easily checks $t = 2$, which puts us in the setup of \cite[Theorem 1.4]{zbMATH07481643} after polarizing $\textrm{Pic}^{0}_{X/S} \to S$, with $X = f^{-1}(S)$, by the relative Theta divisor.


\vspace{0.5em}

We now consider the statement when $g \geq 3$ and with $K$ a number field. One again we consider the associated $K$-algebraic family $f : \overline{X} \to \overline{S}$ of genus $g$ stable curves and take $\mathcal{S} \subset S(\overline{\mathbb{Q}})$ to be the set of points where the fibre of the associated variation $\mathbb{V}' = R^{1} \an{f'}_{*} \mathbb{Q}$ of Hodge structures admits a simple factor with complex multiplication; here $f'$ is the restriction of $f$ to the fibre above $S \subset \overline{S}$. Again we note that the Hodge-genericity of $\overline{S}$ implies, via the Andr\'e-Deligne monodromy theorem \cite{Andre1992}, that $\mathbf{H}_{S}$ is the full symplectic group, so in particular both $\mathbf{H}_{S}$ and the local system $\mathbb{V}'$ are absolutely simple. 

To get a polynomial bound on the heights of points of $\mathcal{S}$ in terms of the degree we now wish to apply \autoref{sndsinggivesheights} to our situation, noting that $\mathcal{S} = \mathcal{S}_{\textrm{split}}$ in our case. This means checking that integrating around the vanishing cycles corresponding to two nodal singularities of $X_{s_{0}}$ produces tuples of non-constant $G$-functions satisfying the linear independence condition in \autoref{sndsinggivesheights}. This follows from Andr\'e's result \autoref{numberofvancycles} above, as well as the remark at \cite[pg.192]{zbMATH00041964}, which says that $h^{1}(\Sigma_{X_{s_{0}}}) = g - \sum_{i} p_{g}(C_{i}) \geq 2$, where $X_{s_{0}} = C_{1} \cup \cdots \cup C_{\ell}$ is the decomposition of the singular fibre into its components. 

We are now in the situation where we have a logarithmic Weil height $\theta : S(\overline{\mathbb{Q}}) \to \mathbb{R}_{> 0}$ for which there exists constants $\kappa, a \in \mathbb{R}_{> 0}$ such that
\begin{equation}
\label{xiheightboundcurves}
\theta(\xi) \leq \kappa \, [K(\xi) : K]^{a}
\end{equation}
for all $\xi \in \mathcal{S}$. By the resolution of the Andr\'e-Oort conjecture \cite{fullandreoort}, it suffices to replace $\mathcal{S}$ with just those points admitting a CM summand $W$ but which are not CM themselves. Using \cite[Prop. 1.3.2.1]{zbMATH06244724} we may replace $W$ with the smallest isotypic summand containing it, which again has CM and is not equal to $\mathbb{V}'_{\xi}$. The idempotent $e$ corresponding to the isotypic summand $W$ is then a non-trivial central idempotent that defines $\xi$ inside $S$. After replacing the problem in question for the equivalent problem for the associated Jacobian family, the result then follows from \autoref{idempotcor} as soon as $g \geq 3$. In particular Hodge-genericity implies the monodromy assumption of \autoref{idempotcor}, and the atypicality with respect to the idempotents of the points in $\mathcal{S}$  holds under the assumption $g \geq 3$ as a special locus of $\mathcal{A}_{g}$ with a global isogeny factor has codimension at least two. Note that to verify the hypothesis in \autoref{idempotcor} that the points in question be ``defined by'' their idempotents in the sense of \autoref{defaspvar2} one may take the $h$ and $h'$ of \autoref{centdecomplem} to be the $h$ and $h'$ in \autoref{defaspvar} and \autoref{defaspvar2}.

\subsection{Families of Abelian Varieties}

We now prove \autoref{abfamthm}. The proof is the same as the case for curves above, except we do not even have to check that the vanishing cycles are independent because this has been assumed. 

\subsection{Degenerations to Hyperplanes}
\label{hyperplanesec}

Lastly, we prove \autoref{hypsufcor}. We once again apply \autoref{sndsinggivesheights} using Andr\'e's result \autoref{numberofvancycles} above. That Hodge-genericity implies the absolute simplicity of the primitive variation is a consequence of Beauville's computation of the monodromy groups of hypersurface variations in \cite[Thm. 2, Thm. 4]{beauville1986groupe}; here we note that a Hodge generic curve on which the associated primitive variation is non-constant will have maximum-possible algebraic monodromy because its algebraic monodromy group must be non-trivial and normal in the absolutely simple group of automorphisms preserving the polarization form. 

To apply \autoref{sndsinggivesheights} it therefore suffices to compute the degree $n$ cohomology of the dual graph associated to a generic hyperplane arrangement. By \cite[Thm. 1.1]{zbMATH06228506}, when $Y \subset V$ is a divisor with simple normal crossings in a smooth complete variety $V$, the cohomology of $\Sigma_{Y}$ is determined by the complement $U := V - Y$; we explain this in the concrete case when $V = \mathbb{P}^{n+1}$ and $Y$ is a union of $d$ hyperplanes. We start by considering the long exact sequence of cohomology with compact support
\begin{equation}
\label{longexactcompact}
 \cdots \to H^{i}_{c}(\mathbb{P}^{n+1},\mathbb{Q}) \to H^{i}_{c}(Y,\mathbb{Q}) \to H^{i+1}_{c}(U,\mathbb{Q}) \to H^{i+1}_{c}(\mathbb{P}^{n+1},\mathbb{Q}) \to \cdots 
\end{equation}
As explained in \cite[\S2(W3)]{zbMATH06228506}, this is canonically an exact sequence of mixed Hodge structures in such a way so that it stays exact upon applying $W_{j}$ for each $j$, and in particular $W_{0}$. Moreover for each complex algebraic variety $T$ the cup product \cite[\S B.1.2]{zbMATH05233837} maps $H^{i}(T, \mathbb{Q}) \otimes_{\mathbb{Q}} H^{j}_{c}(T,\mathbb{Q}) \to H^{i+j}_{c}(T, \mathbb{Q})$ are morphisms of mixed Hodge structures \cite[\S6.3]{zbMATH05233837}, which implies via Poincar\'e duality \cite[Cor. B.25]{zbMATH05233837} that $H^{i}_{c}(T,\mathbb{Q}) \simeq H^{2 m - i}(T,\mathbb{Q})^{*} \otimes \mathbb{Q}(-m)$ as mixed Hodge structures when $T$ is irreducible of dimension $m$. We may apply \cite[Thm. 1]{zbMATH00562914} to each $H^{i}(U,\mathbb{Q})$, noting that $U$ becomes the complement of an affine hyperplane arrangement after removing the first hyperplane from $\mathbb{P}^{n+1}$. Then $H^{i}(U,\mathbb{Q})$ is pure of weight $2i$, and hence taking $T = U$ and $i = m = n+1$ we obtain by duality that $H^{n+1}_{c}(U,\mathbb{Q})$ is pure of weight zero. Applying $W_{0}$ to (\ref{longexactcompact}) with $i = n$ one obtains $W_{0} H^{n}_{c}(Y,\mathbb{Q}) \simeq H^{n+1}_{c}(U,\mathbb{Q})$. But $H^{n}_{c}(Y,\mathbb{Q}) = H^{n}(Y,\mathbb{Q})$ since $Y$ is compact, hence applying \cite[Lem. 1.3]{zbMATH06228506} this means that $\dim_{\mathbb{C}} H^{n}(\Sigma_{Y}, \mathbb{C}) = \dim_{\mathbb{C}} H^{n+1}_{c}(U,\mathbb{C})$. 

Applying Poincar\'e duality again we reduce to computing $\dim_{\mathbb{C}} H^{n+1}(U,\mathbb{C})$. One may deduce from \cite[Lemma 3]{zbMATH03436327} that the dimension of this cohomology group is ${d-1 \choose n+1}$, which under the assumption $d > n+2$ is always at least $2$.

As is clear from the proof, a more detailed analysis of the cohomology of the dual graph associated to the singular fibre gives more results for other degeneration types as well, but we do not undertake this analysis here.

\appendix

\section{Orr's Theorem for Self-Adjoint Operators}

\begin{prop}
\label{polyindenomlem}
Suppose that $\eta \in \textrm{End}(\mathbb{R}^{m})$ is any element which is self-adjoint for a positive-definite bilinear form $\psi : \mathbb{R}^{m} \times \mathbb{R}^{m} \to \mathbb{R}$. Let $\mathfrak{S} = \Lambda B_{t} K$ be a Siegel set for $\GL_{m}$ with $K$ the orthogonal group of $\psi$. Then all rational points in the orbit $\mathfrak{S} \cdot \eta$ under the conjugation action have height uniformly bounded by a polynomial in their denominators.
\end{prop}


\begin{proof}
Applying \cite[Lem 3.8]{zbMATH06880895} there exists $\gamma \in \GL_{m}(\mathbb{Q})$ and $\sigma \in \GL_{m}(\mathbb{R})$ such that $\gamma^{-1} \mathfrak{S} \gamma \sigma$ is contained in a standard Siegel set $\mathfrak{T} = \Xi C_{u} J$ in the sense of \cite[\S2]{zbMATH06880895}. Since conjugating by a fixed rational $\gamma$ maps rational points to rational points, we may reduce to the situation where $\gamma = 1$. Then given $s \in \mathfrak{S}$ we have $s \sigma = t$ for some $t \in \mathfrak{T}$, and hence
\begin{equation*}
s \eta s^{-1} = t (\sigma^{-1} \eta \sigma) t^{-1} .
\end{equation*}
From the proof of \cite[Lem 3.8]{zbMATH06880895} we see that $\gamma \sigma$ ($= \sigma$ for us) conjugates the standard orthogonal group into $K$, so $\tau = \sigma^{-1} \eta \sigma$ is a symmetric matrix. Replacing $\mathfrak{S}$ by $\mathfrak{T}$ and $\eta$ by $\tau$, we may therefore reduce to the situation where $\mathfrak{S}$ is a standard Siegel set and $\eta$ is symmetric.

For $s \in \Lambda B_{t} K$ we now write $s = \lambda b k$ in accordance with the Siegel set decomposition. Let $\rho = k \eta k^{-1} = k \eta k^{t}$; we note that $\rho$ is symmetric and by compactness of $K$ its entries are uniformly bounded in magnitude by some constant $\kappa$. Let $b_{i}$ be the $i$'th diagonal entry of $b$, and recall \cite[\S2A(4)]{zbMATH06880895} that $b_{i} \geq t b_{i+1}$ for some constant $t$. Define a partial order on pairs $(i,j)$ by saying that $(i,j) \succcurlyeq (i',j')$ when $i' \geq i$ and $j \geq j'$; note the reversal of the first inequality. Define
\[ \beta := b \rho b^{-1} = [b_{i} \rho_{ij} b_{j}^{-1}] . \]
From the fact that $b_{i} \geq t b_{i+1}$ for all $i$ we learn that 
\[ \beta_{ij} \geq t^{(i'-i)+(j-j')} (\rho_{ij}/\rho_{i'j'}) \beta_{i'j'} \]
whenever $(i,j) \succcurlyeq (i',j')$ and both $\beta_{ij} \neq 0$ and $\beta_{i'j'} \neq 0$. We also calculate that $1/\beta_{ij} = \frac{1}{\rho^{2}_{ij}} \beta_{ji}$ using the symmetry of $\rho$. Let $\eta^{s} = s \eta s^{-1} = \lambda \beta \lambda^{-1}$ be a proposed rational matrix in the orbit, and write $d_{ij}$ for the absolute value of the denominator of $\eta^{s}_{ij}$. 

Suppose that $(i,j)$ is a minimal element for the order $\succcurlyeq$ subject to the condition that $\eta^{s}_{ij} \neq 0$. Applying \autoref{maxfliplem} below to $(i+1,j)$ and $(i,j-1)$ one learns that this is equivalent to the index $(j,i)$ being maximal for $\succcurlyeq$ subject to the condition that $\beta_{ji} \neq 0$. Calculating as in (\ref{etascalc}) below and using the $\lambda$ has diagonal entries equal to $1$ (see \cite[2A.(5)]{zbMATH06880895}) we have $\eta^{s}_{ij} = \beta_{ij}$. Then for any $(j',i') \preccurlyeq (j,i)$ with $\beta_{j'i'} \neq 0$ we have
\begin{align*}
d_{ij} &\geq \frac{1}{|\beta_{ij}|} \\ 
&= \frac{1}{\rho^{2}_{ij}} |\beta_{ji}| \\
&\geq t^{(j'-j)+(i-i')} \frac{1}{|\rho_{ij} \rho_{i'j'}|} |\beta_{j'i'}| \\
&\geq t^{(j'-j)+(i-i')} \frac{1}{\kappa^2} \beta_{j'i'} .
\end{align*}
Since every non-zero entry of $\beta$ is equal to or lies below a non-zero entry of maximal index, it follows that the entries of $\beta$ have magnitude bounded by a uniform polynomial in the denominators $d_{ij}$. Because $\Lambda$ is compact (again see \cite[2A.(5)]{zbMATH06880895}), this is also true for the entries of $\eta^{s} = \lambda \beta \lambda^{-1}$. Finally, we have
\[ \widetilde{\theta}(\eta^{s}) \leq (\textrm{max}_{ij} \{ d_{ij} \}) \cdot \textrm{max}_{ij}\{ 1, |\eta^{s}_{ij}| \} , \]
where we recall that $\widetilde{\theta}$ is our na\"ive Weil height. This completes the proof.
\end{proof}

\begin{lem}
\label{maxfliplem}
We have $\eta^{s}_{k\ell} = 0$ for all $(k,\ell) \preccurlyeq (i,j)$ if and only if $\beta_{\ell k} = \beta_{k\ell} = 0$ for all $(\ell,k) \succcurlyeq (j,i)$.
\end{lem}

\begin{proof}
Note that $(k,\ell) \preccurlyeq (i,j)$ and $(\ell,k) \succcurlyeq (j,i)$ are the same condition. We argue by upward induction on the partial order $\succcurlyeq$ on pairs, beginning with the smallest case $(i,j) = (m,1)$. Then $\eta^{s}_{m1} = \beta_{m1}$ so the result is immediate.

Next we suppose the lemma holds for $(i+1,j)$ and $(i,j-1)$ and prove the forward implication. Arguing inductively, this just requires showing that $\beta_{ij} = \beta_{ji} = 0$. Using the fact that $\lambda$ and $\lambda^{-1}$ are upper triangular, we have
\begin{align}
\eta^{s}_{ij} &= \sum_{q = 1}^{m} \sum_{r=1}^{m} \lambda_{iq} \beta_{q r} [\lambda^{-1}]_{r j} \\
\label{etascalc}
&= \sum_{q \geq i} \sum_{r \leq j} \lambda_{i q} \beta_{q r} [\lambda^{-1}]_{r j} .
\end{align}
If $\eta^{s}_{k\ell} = 0$ for all $(k, \ell) \preccurlyeq (i,j)$ then we may apply the induction hypothesis to conclude that all terms in the sum vanish except possibly the term $\lambda_{ii} \beta_{ij} [\lambda^{-1}]_{jj}$. But this only vanishes if $\beta_{ij} = 0$. Then since $\beta_{ij} = b_{i} \rho_{ij} b_{j}^{-1}$ and one has $\rho_{ij} = \rho_{ji} = 0$ by the symmetry of $\rho$, we get $\beta_{ji} = 0$ as well.

Conversely if $\beta_{\ell k} = \beta_{k \ell} = 0$ for all $(\ell, k) \succcurlyeq (j,i)$ then plainly all terms in the sum vanish, hence $\eta^{s}_{ij}$ vanishes, and this continues to hold as $(i,j)$ decreases in the order $\preccurlyeq$.
\end{proof}

\section{An Alternative Approach to \S3}

After writing this paper, the author discovered that the arguments in \S3 can be clarified using recent work in $p$-adic Hodge theory. In particular, it is possible to simplify the construction of the maps $\widehat{\gamma}^{*}_{\textrm{\'et}}$ of \S\ref{funcextsec} and verify their compatibility with the corresponding $\widehat{\gamma}^{*}_{\textrm{dR}}$ by referencing computations that have appeared in the literature since Scholze's paper \cite{zbMATH06209107}.

We consider a pair $(Y, \Delta^{a,b})$ as in \S\ref{funcextsec}, except with $\Delta^{a,b} \subset Y$ a product $\Delta^{a,b} = (\Delta^{\circ})^{a} \times \Delta^{b}$ of \emph{open} annuli and disks instead of closed ones. Working with open annuli and disks ensures that the space $\Delta^{a,b}$ is Stein (see \cite[\S3.1.1]{zbMATH07160161} for a definition) and allows us to apply the results of \cite{zbMATH07160161} to compute the pro-\'etale cohomology of $\Delta^{a,b}$ with coefficients in $\widehat{\mathbb{Q}}_{p}(a) = \widehat{\mathbb{Z}}_{p}(a)[1/p]$. Note that the construction in the proof of \autoref{bigsupergfuncthm}(ii) also produces such rigid subvarieties in the fibres $X^{\textrm{ad}}_{s_{1}}$, so the difference is immaterial for our application. 

Now, it is a consequence of \cite[Thm. 4.12]{zbMATH07160161} that:

\begin{prop}
\label{natexactsecprop}
There is a natural short exact sequence
\[ 0 \to \Omega^{a-1}(\Delta^{a,b}_{\mathbb{C}_{p}})/\ker d \to H^{a}_{\textrm{p\'et}}(\Delta^{a,b}_{\mathbb{C}_{p}}, \widehat{\mathbb{Q}}_{p}(a)) \to (H^{a}_{\textrm{HK}}(\Delta^{a,b}) \widehat{\otimes}_{F} B^{+}_{\textrm{st}})^{N = 0, \varphi=p^{a}} \to 0 . \]
The term $(H^{a}_{\textrm{HK}}(\Delta^{a,b}) \widehat{\otimes}_{F} B^{+}_{\textrm{st}})^{N = 0, \varphi=p^{a}}$ is isomorphic to $\mathbb{Q}_{p}$, and a choice of coordinate $T \in \mathcal{O}(\Delta^{\circ})$ induces a Galois-equivariant splitting 
\begin{equation}
\label{canonicsplit}
H^{a}_{\textrm{p\'et}}(\Delta^{a,b}_{\mathbb{C}_{p}}, \widehat{\mathbb{Q}}_{p}(a)) \cong (\Omega^{a-1}(\Delta^{a,b}_{\mathbb{C}_{p}})/\ker d) \oplus \mathbb{Q}_{p} 
\end{equation}
such that $1 \in \mathbb{Q}_{p}$ is the image of the Hyodo-Kato cohomology class of $dT_{1}/T_{1} \wedge \cdots \wedge dT_{a}/T_{a}$. 
\end{prop}
The object $H^{a}_{\textrm{HK}}(\Delta^{a,b})$ is the Hyodo-Kato cohomology of $\Delta^{a,b}$, which is a version of the de Rham cohomology of $\Delta^{a,b}$ but with Frobenius and monodromy operators $\varphi$ and $N$. In \cite[Thm. 4.12]{zbMATH07160161} it is defined using a choice of semistable model, but it can be also be defined intrinsically on the level of rigid varieties if one follows the approach of \cite[\S3]{arXiv:2306.06100}. The computation of $(H^{a}_{\textrm{HK}}(\Delta^{a,b}) \widehat{\otimes}_{F} B^{+}_{\textrm{st}})^{N = 0, \varphi=p^{a}}$ can be carried out analogously to the argument in \cite[\S4.3.2]{zbMATH07160161} for products of tori by also including $b$ trivial factors. The description of the section giving the splitting is given in \cite[Lem 4.30]{arXiv:2308.07712} in the $a = 1, b = 0$ case, and the section for the general case is constructed using cup product.

What is important for us is the following fact:
\begin{lem}
The summand $\Omega^{a-1}(\Delta^{a,b}_{\mathbb{C}_{p}})/\ker d$ in (\ref{canonicsplit}) is exactly the kernel of the map 
\begin{equation}
\label{QptoBdRker}
H^{a}_{\textrm{p\'et}}(\Delta^{a,b}_{\mathbb{C}_{p}}, \widehat{\mathbb{Q}}_{p}(a)) \to H^{a}_{\textrm{p\'et}}(\Delta^{a,b}_{\mathbb{C}_{p}}, \mathbb{B}_{\textrm{dR}}) .
\end{equation}
\end{lem}

\begin{proof}
This follows from Bosco's description of the exact sequence in \autoref{natexactsecprop}. More precisely, in the proof of \cite[Thm. 7.7]{arXiv:2306.06100} Bosco constructs this exact sequence in the following way. He considers the so-called ``fundamental exact sequence'' of $p$-adic Hodge theory (cf. \cite[Cor. 2.26]{arXiv:2306.06100})
\begin{equation}
\label{fundexactsec}
0 \to \widehat{\mathbb{Q}}_{p} \to \mathbb{B}_{e} \to \mathbb{B}_{\textrm{dR}} / \mathbb{B}^{+}_{\textrm{dR}} \to 0
\end{equation}
and shows that the exact sequence in \autoref{natexactsecprop} can be identified, up to twisting, with a limit of exact sequences of the form
\begin{equation}
\label{boscover}
0 \to \textrm{coker}\, \alpha_{a-1} \to H^{a}_{\textrm{p\'et}}(V^{\dagger}, \widehat{\mathbb{Q}}_{p}) \to \textrm{ker}\, \alpha_{a} \to 0
\end{equation}
where the maps $\alpha_{i}$ and the maps in (\ref{boscover}) come from the long exact sequence
\[ \cdots \to H^{i}_{\textrm{p\'et}}(V^{\dagger}, \widehat{\mathbb{Q}}_{p}) \to H^{i}_{\textrm{p\'et}}(V^{\dagger}, \mathbb{B}_{e}) \xrightarrow{\alpha_{i}} H^{i}_{\textrm{p\'et}}(V^{\dagger}, \mathbb{B}_{\textrm{dR}}/\mathbb{B}^{+}_{\textrm{dR}}) \to \cdots \]
induced by (\ref{fundexactsec}), and the $V^{\dagger} \subset \Delta^{a,b}_{\mathbb{C}_{p}}$ range over an appropriate cover consisting of dagger subvarieties (for instance a collection of affinoid products of annuli and disks of the same type, with their canonical dagger structure). The map $\widehat{\mathbb{Q}}_{p} \to \mathbb{B}_{\textrm{dR}}$ may be factored as $\widehat{\mathbb{Q}}_{p} \to \mathbb{B}_{e} \to \mathbb{B}_{dR}$, from which one can show that 
\[ \varprojlim_{V} \textrm{coker}\, \alpha_{a-1} \subset \ker \left[ H^{a}_{\textrm{p\'et}}(\Delta^{a,b}_{\mathbb{C}_{p}}, \widehat{\mathbb{Q}}_{p}) \to H^{a}_{\textrm{p\'et}}(\Delta^{a,b}_{\mathbb{C}_{p}}, \mathbb{B}_{\textrm{dR}}) \right] \]
under this identification; in particular, $\Omega^{a-1}(\Delta^{a,b}_{\mathbb{C}_{p}})/\ker d$ is contained in the kernel of (\ref{QptoBdRker}). It then suffices to show that the map (\ref{QptoBdRker}) is not zero. Using \autoref{BdRcohomologyisnaive} below this reduces, via cup product, to the $a = 1, b = 0$ case, where this is shown in the proof of \autoref{dtnonzero}. 
\end{proof}

One can also show:

\begin{lem}
\label{BdRcohomologyisnaive}
We have $H^{a}_{\textrm{p\'et}}(\Delta^{a,b}_{\mathbb{C}_{p}}, \mathbb{B}_{\textrm{dR}}) \cong B_{\textrm{dR}}$. More generally, one has $H^{i}_{\textrm{p\'et}}(\Delta^{a,b}_{\mathbb{C}_{p}}, \mathbb{B}_{\textrm{dR}}) = \bigwedge^{i} \bigoplus_{j=1}^{a} B_{\textrm{dR}}$. 
\end{lem}

\begin{proof}
The case $a = 1, b = 0$ is given explicitly toward the end of \cite[Rem. 3.21]{zbMATH07009770}. More generally, \cite[Rem. 3.21]{zbMATH07009770} explains how to compute the cohomology groups $H^{i}_{\textrm{p\'et}}(X_{\mathbb{C}_{p}}, \mathbb{B}_{\textrm{dR}})$ when $X$ is a more general Stein space. To do this, one uses the quasi-isomorphism in \cite[Rem. 3.21]{zbMATH07009770} to replace the sheaf $\mathbb{B}^{+}_{\textrm{dR}}$ by the de Rham complex in loc. cit. and performs the usual computation of the de Rham cohomology. To obtain the result for $\mathbb{B}_{\textrm{dR}}$ one takes a limit in the manner described in \cite[Rem. 3.21]{zbMATH07009770}. These calculations are compatible with cup product for products of Stein spaces.
\end{proof}

With the above considerations in mind, the eight-object diagram at the beginning of \S\ref{funcextsec} can be replaced with the simpler diagram 
\begin{equation*}
\begin{tikzcd}[center picture]
H^{a}(Y_{\overline{K_{v}}}, \hat{\mathbb{Z}}_{p}(a)) \otimes B_{\textrm{dR}} \arrow[r, "\sim"] \arrow[d] & H^{a}(Y_{\mathbb{C}_{p}}, \mathbb{B}_{\textrm{dR}}) \arrow[r, "\sim"] \arrow[d] & \arrow[l, swap, "\sim"] H^{a}(Y, \Omega^{\bullet}_{Y}) \otimes B_{\textrm{dR}} \arrow[d] \\
(\ker \oplus \mathbb{Q}_{p}) \otimes B_{\textrm{dR}} \arrow[r] & H^{a}(\Delta^{a,b}_{\mathbb{C}_{p}}, \mathbb{B}_{\textrm{dR}}) \arrow[r, "\sim"] &  \arrow[l] H^{a}(\Delta^{a,b}, \Omega^{\bullet}_{\Delta^{a,b}}) \otimes B_{\textrm{dR}}
\end{tikzcd}
\end{equation*} 
where ``$\ker$'' refers to the kernel of the map $[H^{a}(\Delta^{a,b}_{\overline{K_{v}}}, \hat{\mathbb{Z}}_{p}(a)) \to H^{a}(\Delta^{a,b}_{\mathbb{C}_{p}}, \mathbb{B}_{\textrm{dR}})]$. After modding out by ``$\ker$'' all horizontal arrows in the diagram become isomorphisms. The map $\widehat{\gamma}^{*}_{\textrm{dR}}$ is defined as before and the functional $\widehat{\gamma}^{*}_{\textrm{\'et}}$ on the \'etale side is the obvious map $H^{a}(\Delta^{a,b}_{\overline{K_{v}}}, \hat{\mathbb{Z}}_{p}(a)) \otimes \mathbb{Q}_{p} \to \mathbb{Q}_{p}$ induced by the splitting. The fact that $\widehat{\gamma}^{*}_{\textrm{dR}}$ and $\widehat{\gamma}^{*}_{\textrm{\'et}}$ are identified after tensoring with $B_{\textrm{dR}}$ --- i.e., the analogue of \autoref{indepofpathlem} in this setting --- is then immediate. 

\bibliography{hodge_theory}
\bibliographystyle{alpha}

\end{document}